\documentclass[12pt]{article}
\usepackage{amssymb,amsmath,amsthm,tikz,multirow}
\usetikzlibrary{arrows,calc}

\title{Tilings of the Sphere by Edge Congruent Pentagons}
\author{Ka Yue Cheuk, Ho Man Cheung, Min Yan\thanks{Research was supported by Hong Kong RGC General Research Fund 605610 and 606311.} \\ 
Hong Kong University of Science and Technology}

\newcommand{\sub}{\subset}
\newcommand{\pa}{\partial}
\newcommand{\mb}{\mathbf}
\newcommand{\mc}{\mathcal}
\newcommand{\bb}{\mathbb}

\newtheorem{theorem}{Theorem}
\newtheorem{lemma}[theorem]{Lemma}

\newtheorem*{theorem*}{Theorem}

\theoremstyle{definition}
\newtheorem*{definition*}{Definition}
\newtheorem*{case*}{Case}
\newtheorem*{subcase*}{Subcase}

\theoremstyle{remark}

\numberwithin{equation}{section}

\begin{document}

\maketitle

\begin{abstract}
We study edge-to-edge tilings of the sphere by edge congruent pentagons, under the assumption that there are tiles with all vertices having degree $3$. We develop the technique of neighborhood tilings and apply the technique to completely classify edge congruent earth map tilings. 
\end{abstract}

\section{Introduction}

Mathematicians have studied tilings for a long time. Among many important discoveries is the complete classification of edge-to-edge tilings of the sphere by congruent triangles \cite{so,ua}. It is an easy combinatorial fact that in an edge-to-edge tiling of the sphere by congruent polygons such that all vertices have degree $\ge 3$, the tile must be triangle, quadrilateral, or pentagon \cite{ua2}. According to \cite{ua2}, the study of quadrilateral spherical tilings is rather complicated. We believe pentagonal spherical tilings should be relatively easier to study because $5$ is an ``extreme'' among $3$, $4$, $5$.  

For triangles on the sphere, the congruence in terms of the edge length is equivalent to the congruence in terms of the angle. For pentagons, however, the two congruences are no longer equivalent. Moreover, we need to know the underlying combinatorial structure of the tiling. Therefore the study of the spherical tilings by (fully geometrically) congruent pentagons needs to be divided into the combinatorial, edge length and angle aspects. In \cite{gsy}, by independently studying the three aspects and then combining the respective classifications together, Gao, Shi and Yan completely classified the minimal case of the edge-to-edge tilings of the sphere by $12$ pentagons. 

In \cite{yan}, we attempted to study the combinatorial aspects of the spherical pentagonal tilings. We found the next simplest combinatorial structure beyond the minimal case (which is the dodecahedron). We called these earth map tilings.

In this paper, we attempt to study the edge length aspects. Two polygons are {\em edge congruent} if there is a one-to-one correspondence between the edges that preserve the lengths and the adjacency relations. Specifically, we will study edge-to-edge tilings of the sphere by edge congruent pentagons satisfying the following combinatorial conditions:
\begin{enumerate}
\item The number of tiles is $>12$.
\item All vertices have degree $\ge 3$.
\item There is a tile such that all its vertices have degree $3$.
\end{enumerate}
The first condition simply means that we are beyond the minimal case. The second condition essentially assumes that all vertices are natural, which simplifies our study. The only substantial condition is the third one. Under the second condition, it is an easy combinatorial consequence that, if $v_i$ is the number of vertices of degree $i$ (recall that the second condition means $v_1=v_2=0$), then
\[
v_3=20+\sum_{i\ge 4}(3i-10)v_i=20+2v_4+5v_5+8v_6+\cdots.
\]
This shows that degree $3$ vertices dominate. In fact, Lemma 3 of \cite{gsy} says that there must be a tile with at least four vertices of degree $3$. So we expect that the third condition is satisfied by many pentagonal tilings, and failing the condition should impose substantial constraints on the tiling. 

In a subsequent paper \cite{ccy}, we will further consider the angles and study the full geometrical congruence. We also note that Luk and Yan \cite{luk,ly} studied the numerical aspects of the angle congruence for general edge-to-edge spherical pentagonal tilings.

Now we describe the results of this paper. Let $P$ be a tile in an edge congruent tiling, such that all the vertices of $P$ have degree $3$. By Proposition 8 of \cite{gsy}, the edge lengths of the five edges of $P$ must be arranged in one of the following five ways (up to rotations and flippings): 
\begin{align*}
a^5&\colon a,a,a,a,a; \\
a^4b&\colon a,a,a,a,b; \\
a^2b^2c&\colon a,a,b,b,c; \\
a^3bc&\colon a,a,a,b,c; \\
a^3b^2&\colon a,a,a,b,b.
\end{align*}
In Section \ref{section_nd}, we develop the key technique, which is the classification of the edge congruent tilings of the neighborhood of $P$. Furthermore, we find all the ways such neighborhood tilings propagate.

In Section \ref{section_3abc}, we recall the earth map tilings introduced in \cite{yan}, which are characterized by the property of having exactly two vertices of degree $>3$. Then we apply the neighborhood technique of Section \ref{section_nd} to prove that there are no edge congruent earth map tilings for the edge length combination $a^3bc$. In Sections \ref{section_2a2bc} and \ref{section_3a2b}, we further classify edge congruent earth map tilings for the combinations $a^2b^2c$ and $a^3b^2$. In Section \ref{section_4ab}, using a more direct method instead of the neighborhood technique, we further classify edge congruent earth map tilings for the combination $a^4b$. There is clearly only one family for the combination $a^5$. So we get the complete classification of edge congruent earth map tilings.

The technique of this paper can be used to classify other edge congruent pentagonal tilings as long as we know enough about the combinatorial structures, and there are sufficiently many tiles with all vertices having degree $3$. However, the lesson we learn from our initial attempts to go beyond the minimal case is that the independent study of each aspect has its limits. For example, as edge lengths become more and more equal (i.e., the $a^4b$ case in this paper), the edge congruent classification gets too complicated to be useful and the additional angle consideration should be introduced to cut down the number of possibilities. In general, we expect that the interaction between the combinatorial, edge length and angle aspects become more and more important as we have fewer and fewer information about certain aspect of the tilings.

\section{Neighborhood Tiling}
\label{section_nd}

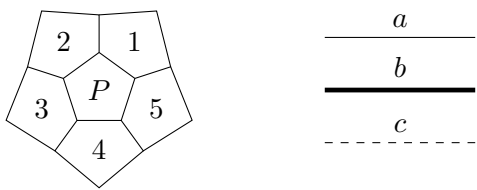
\begin{figure}[htp]
\centering
\begin{tikzpicture}[>=latex]

\foreach \x in {1,...,5}
{
\draw[rotate=72*\x]
	(90:0.5) -- (18:0.5) -- (18:1) -- (54:1.3) -- (90:1);
\node at (-18+72*\x:0.8) {\small \x};
}

\draw
	(3,0.7) -- node[above] {\small $a$} ++(2,0);
\draw[line width=1.8]
	(3,0) -- node[above] {\small $b$} ++(2,0);
\draw[dashed]
	(3,-0.7) -- node[above] {\small $c$} ++(2,0);
	
\node at (0,0) {\small $P$};

\end{tikzpicture}
\caption{Neighborhood of a tile with all vertices having degree $3$.}
\label{neighborhood}
\end{figure}

Let $P$ be a tile in an edge-to-edge pentagonal tiling. The {\em neighborhood} of $P$ is the collection of tiles with at least one common vertex with $P$. If all the vertices of $P$ have degree $3$, then the combinatorial structure of the neighborhood is given in Figure \ref{neighborhood}. The nearby tiles are labeled by $1,\dots,5$ and denoted $P_1,\dots,P_5$. The edge shared by $P_i$ and $P_j$ is denoted $E_{ij}$. We will also indicate edge lengths $a,b,c$ by normal, thick and dashed lines.

\begin{lemma}\label{edge_nd}
For distinct $a,b,c$, we have the following edge congruent tilings of the neighborhood of a tile with all vertices having degree $3$.
\begin{enumerate}
\item If the pentagon has edge length combination $a^2b^2c$, then the tiling is uniquely given by the right of Figure \ref{nd2a2bc} up to symmetry.
\item If the pentagon has edge length combination $a^3bc$, then the tiling is given by one of two in Figure \ref{nd3abc} up to symmetry.
\item If the pentagon has edge length combination $a^3b^2$, then the tiling is given by one of three in Figure \ref{nd3a2b} up to symmetry.
\item If the pentagon has edge length combination $a^4b$, then the tiling is given by one of eighteen in Figure \ref{nd4ab} up to symmetry.
\end{enumerate}
\end{lemma}

There is obviously only one edge congruent tiling of the neighborhood for the combination $a^5$.

\begin{proof}
For the edge combination $a^2b^2c$, we start with the center tile in Figure \ref{nd2a2bc}, with the edge lengths arranged in the (unique, up to the symmetries of rotations and flippings) order of $a,a,b,b,c$. The edge $E_{45}$ is adjacent to $c$ and therefore cannot be $c$. If $E_{45} = b$, then all the edges of $P_4$ are determined. The two $a$-edges of $P_3$ imply that $E_{23} = b$ or $c$. In either case, we have $E_{12} = a$. Then all the edges of $P_1$ are determined, and we find three $b$-edges in $P_5$, a contradiction. So we must have $E_{45} = a$. Then we can successively determine all the edges of $P_4$, $P_3$, $P_5$. We already know two $a$-edges of $P_2$ and two $b$-edges of $P_1$. This forces $E_{12}=c$, so that all the edges of $P_1$ and $P_2$ are determined.

\begin{figure}[htp]
\centering
\begin{tikzpicture}[>=latex]

\foreach \x in {1,...,5}
\foreach \a in {0,1}
{
\begin{scope}[xshift=3*\a cm]

\coordinate (A\x X\a) at (-54+72*\x:0.5);
\coordinate (B\x X\a) at (-54+72*\x:1);
\coordinate (C\x X\a) at (-18+72*\x:1.3);
\coordinate (P\x X\a) at (-18+72*\x:0.8);

\node at (-18+72*\x:0.8) {\small \x};

\end{scope}
}

\draw
	(C1X0) -- (B2X0) -- (A2X0) -- (A3X0) -- (A4X0) -- (B4X0) -- (C4X0) 
	(B1X1) -- (C1X1) -- (B2X1) 
	(A3X1) -- (B3X1)
	(A2X1) -- (A3X1) -- (A4X1)
	(C4X1) -- (B5X1) -- (C5X1)
	(A5X1) -- (B5X1)
	;

\draw[line width=1.8]
	(A2X0) -- (A1X0) -- (A5X0) -- (B5X0) -- (C4X0) 
	(A1X0) -- (B1X0)
	(B2X1) -- (C2X1) -- (B3X1) 
	(A2X1) -- (A1X1) -- (A5X1)
	(A1X1) -- (B1X1)
	(C3X1) -- (B4X1) -- (C4X1)
	(A4X1) -- (B4X1);

\draw[dashed]
	(B1X0) -- (C1X0) 
	(A4X0) -- (A5X0)
	(A2X1) -- (B2X1) 
	(B1X1) -- (C5X1)
	(B3X1) -- (C3X1)
	(A4X1) -- (A5X1)
	;

\draw[gray!70]
	(B2X0) -- (C2X0) -- (B3X0) -- (C3X0) -- (B4X0) 
	(A3X0) -- (B3X0)
	(B1X0) -- (C5X0) -- (B5X0)
	;

\end{tikzpicture}
\caption{Neighborhood tiling for the edge combination $a^2b^2c$.}
\label{nd2a2bc}
\end{figure}
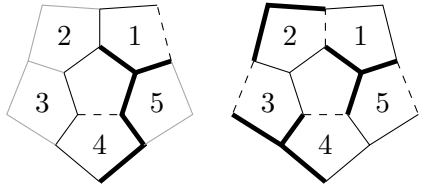

For the edge combination $a^3bc$, we start with the center tile in Figure \ref{nd3abc}, with the edge lengths arranged in the order of $a,a,a,b,c$. Since $E_{12}$ is adjacent to both $b$ and $c$, we must have $E_{12}=a$. Then we can successively determine all the edges of $P_1$, $P_2$, $P_3$, $P_5$, and we get $E_{34} = E_{45} = a$. The only ambiguity left is the location of the $b$-edge and the $c$-edge in $P_4$, which gives two possibilities.

\begin{figure}[htp]
\centering
\begin{tikzpicture}[>=latex]

\foreach \x in {1,...,5}
\foreach \a in {0,1}
{
\begin{scope}[xshift=3*\a cm]

\coordinate (A\x X\a) at (-54+72*\x:0.5);
\coordinate (B\x X\a) at (-54+72*\x:1);
\coordinate (C\x X\a) at (-18+72*\x:1.3);
\coordinate (P\x X\a) at (-18+72*\x:0.8);

\node at (-18+72*\x:0.8) {\small \x};

\end{scope}
}

\foreach \a in {0,1}
{
\draw
	(B1X\a) -- (C1X\a) -- (B2X\a) -- (C2X\a) -- (B3X\a)
	(A2X\a) -- (B2X\a)
	(A1X\a) -- (A5X\a) -- (A4X\a) -- (A3X\a)
	(A4X\a) -- (B4X\a) -- (C3X\a)
	(A5X\a) -- (B5X\a) -- (C5X\a)
	;

\draw[line width=1.8]
	(A1X\a) -- (B1X\a) 
	(A2X\a) -- (A3X\a)
	(B3X\a) -- (C3X\a)
	;

\draw[dashed]
	(A1X\a) -- (A2X\a)
	(A3X\a) -- (B3X\a)
	(B1X\a) -- (C5X\a)
	;
}

\draw[line width=1.8]
	(B4X0) -- (C4X0)
	(B5X1) -- (C4X1)
	;

\draw[dashed]
	(B5X0) -- (C4X0)
	(B4X1) -- (C4X1)
	;

\node at (0,0) {\small I};
\node at (3,0) {\small II};

\end{tikzpicture}
\caption{Neighborhood tilings for the edge combination $a^3bc$.}
\label{nd3abc}
\end{figure}
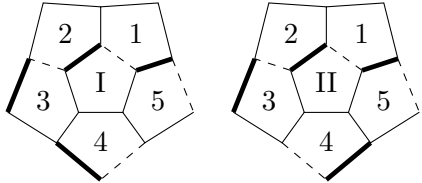

For the edge combination $a^3b^2$, we start with the center tile in Figure \ref{nd3a2b}, with edge lengths arranged in the order of $a,a,a,b,b$. The edge $E_{12}$ is either $a$ or $b$. If $E_{12} = a$, then we may successively determine all the edges of $P_1$, $P_2$, $P_3$, $P_5$, $P_4$ and get the first tiling. If $E_{12} = b$, then we may determine all the edges of $P_1$ and $P_2$. Now among the edges $E_{34}$ and $E_{45}$ of $P_4$, one must be $a$ and the other may be $a$ or $b$. The first case is $E_{34}=E_{45}=a$, and we immediately get all the edges of $P_3$, $P_4$, $P_5$ as in the second tiling. For the second case, up to symmetry, we may assume $E_{34}=b$ and $E_{45}=a$, and also immediately get all the edges of $P_3$, $P_4$, $P_5$  as in the third tiling.

\begin{figure}[htp]
\centering
\begin{tikzpicture}[>=latex]

\foreach \x in {1,...,5}
\foreach \a in {0,1,2}
{
\begin{scope}[xshift=3*\a cm]

\draw[rotate=72*\x]
	(90:0.5) -- (18:0.5) -- (18:1) -- (54:1.3) -- (90:1);
	
\node at (-18+72*\x:0.8) {\small \x};

\coordinate (A\x X\a) at (-54+72*\x:0.5);
\coordinate (B\x X\a) at (-54+72*\x:1);
\coordinate (C\x X\a) at (-18+72*\x:1.3);

\end{scope}
}

\draw[line width=1.8]
	(B1X0) -- (A1X0) -- (A2X0) -- (A3X0) -- (B3X0) 
	(B4X0) -- (C4X0) -- (B5X0)
	(B1X1) -- (C5X1) -- (B5X1) -- (C4X1) -- (B4X1) -- (C3X1) -- (B3X1) 
	(A1X1) -- (A2X1) -- (A3X1)
	(A2X1) -- (B2X1)
	(B1X2) -- (C5X2) -- (B5X2)  
	(C3X2) -- (B4X2) -- (C4X2)
	(A4X2) -- (B4X2)
	(A1X2) -- (A2X2) -- (A3X2)
	(A2X2) -- (B2X2);

\node at (0,0) {\small I};
\node at (3,0) {\small II};
\node at (6,0) {\small III};

\end{tikzpicture}
\caption{Neighborhood tiling for the edge combination $a^3b^2$.}
\label{nd3a2b}
\end{figure}
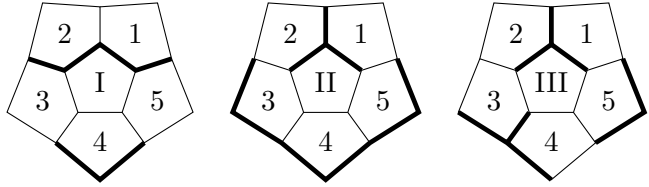

For the edge combination $a^4b$, we start with the center tile in Figure \ref{nd4ab}, with the only $b$-edge at the bottom. Then we look at the five edges out of the five vertices of the center tile. At most two of these can be the $b$-edge. If two are $b$-edges, then we get the unique first tiling. If one is $b$-edge, then this $b$-edge is either ``central'' or ``sideway''. In the first case, we get three possibilities given by the second, third and fourth tilings. In the second case, we get another three possibilities given by the fifth, sixth and seventh tilings.

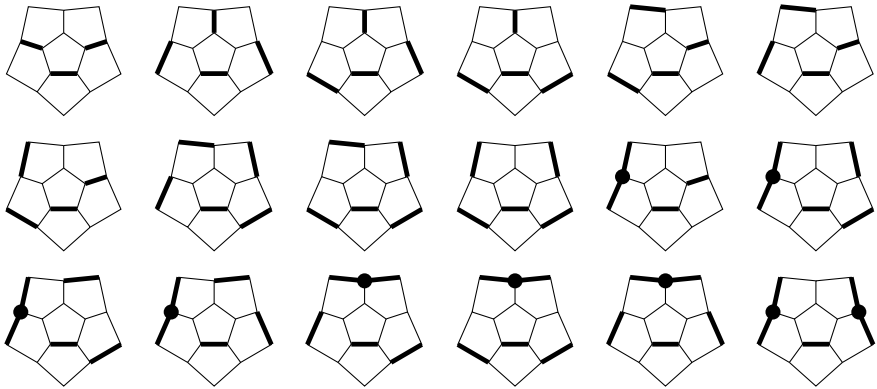
\begin{figure}[htp]
\centering
\begin{tikzpicture}[>=latex]

\foreach \x in {0,...,4}
\foreach \a in {0,...,5}
\foreach \b in {0,...,2}
{
\begin{scope}[shift={(2*\a cm, -1.8*\b cm)}]

\draw
	(-54+72*\x:0.3) -- (18+72*\x:0.3)
	(18+72*\x:0.3) -- (18+72*\x:0.6)
	(-54+72*\x:0.6) -- (-18+72*\x:0.8) -- (18+72*\x:0.6);

\coordinate  (A\x X\a\b) at (18+72*\x:0.3);
\coordinate  (B\x X\a\b) at (18+72*\x:0.6);
\coordinate  (C\x X\a\b) at (54+72*\x:0.8);
	
\end{scope}
}

\foreach \a in {0,...,5}
\foreach \b in {0,...,2}
\draw[line width=1.8]
	(A3X\a\b) -- (A4X\a\b);

\draw[line width=1.8]
	(A0X00) -- (B0X00)	
	(A2X00) -- (B2X00)
	(A1X10) -- (B1X10)	
	(B0X10) -- (C4X10)
	(B2X10) -- (C2X10)
	(A1X20) -- (B1X20)	
	(B0X20) -- (C4X20)
	(B3X20) -- (C2X20)
	(A1X30) -- (B1X30)	
	(B4X30) -- (C4X30)
	(B3X30) -- (C2X30)
	(A0X40) -- (B0X40)	
	(B1X40) -- (C1X40)
	(B3X40) -- (C2X40)
	(A0X50) -- (B0X50)	
	(B1X50) -- (C1X50)
	(B2X50) -- (C2X50)
	(A0X01) -- (B0X01)	
	(B2X01) -- (C1X01)
	(B3X01) -- (C2X01)
	(B4X11) -- (C4X11)	
	(B0X11) -- (C0X11)
	(B1X11) -- (C1X11)
	(B2X11) -- (C2X11)
	(B4X21) -- (C4X21)	
	(B0X21) -- (C0X21)
	(B1X21) -- (C1X21)
	(B3X21) -- (C2X21)
	(B4X31) -- (C4X31)	
	(B0X31) -- (C0X31)
	(B2X31) -- (C1X31)
	(B3X31) -- (C2X31)
	(A0X41) -- (B0X41)	
	(C1X41) -- (B2X41) -- (C2X41)
	(B4X51) -- (C4X51)	
	(B0X51) -- (C0X51)
	(C1X51) -- (B2X51) -- (C2X51)
	(B4X02) -- (C4X02)	
	(B1X02) -- (C0X02)
	(C1X02) -- (B2X02) -- (C2X02)
	(B0X12) -- (C4X12)	
	(B1X12) -- (C0X12)
	(C1X12) -- (B2X12) -- (C2X12)
	(B4X22) -- (C4X22)	
	(C0X22) -- (B1X22) -- (C1X22)
	(B2X22) -- (C2X22)
	(B4X32) -- (C4X32)	
	(C0X32) -- (B1X32) -- (C1X32)
	(B3X32) -- (C2X32)
	(B0X42) -- (C4X42)	
	(C0X42) -- (B1X42) -- (C1X42)
	(B2X42) -- (C2X42)
	(C0X52) -- (B0X52) -- (C4X52)	
	(C1X52) -- (B2X52) -- (C2X52)
	;

\fill
	(B2X41) circle (0.1)
	(B2X51) circle (0.1)
	(B2X02) circle (0.1)
	(B2X12) circle (0.1)
	(B1X22) circle (0.1)
	(B1X32) circle (0.1)
	(B1X42) circle (0.1)
	(B0X52) circle (0.1)
	(B2X52) circle (0.1);

\end{tikzpicture}
\caption{Neighborhood tilings for the edge combination $a^4b$.}
\label{nd4ab}
\end{figure}

It remains to consider the case that there is only one $b$-edge in the interior. We need to assign four $b$-edges on the boundary for the remaining four tiles. We get total of eleven possibilities. In seven of these eleven possibilities, the vertices indicated by the dot must have degree $>3$.
\end{proof}

By the {\em propagation} of the neighborhood tilings, we mean whether one of the nearby tiles $P_1,\dots,P_5$ may also have all its vertices having degree $3$. For the unique neighborhood tiling for the edge combination $a^2b^2c$, this is indeed the case for all the nearby tiles. An example is given by the edge congruent dodecahedron tiling in \cite{gsy}.

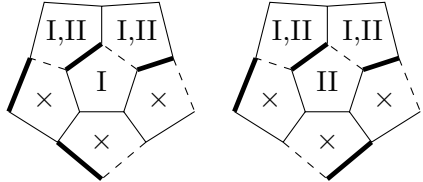
\begin{figure}[htp]
\centering
\begin{tikzpicture}[>=latex]

\foreach \a in {0,1}
{
\begin{scope}[xshift=3*\a cm]

\foreach \x in {1,...,5}
{
\coordinate (A\x X\a) at (-54+72*\x:0.5);
\coordinate (B\x X\a) at (-54+72*\x:1);
\coordinate (C\x X\a) at (-18+72*\x:1.3);
\coordinate (P\x X\a) at (-18+72*\x:0.8);
}

\node at (54:0.8) {\small I,II};
\node at (126:0.8) {\small I,II};
\node at (198:0.8) {\small $\times$};
\node at (-90:0.8) {\small $\times$};
\node at (-18:0.8) {\small $\times$};

\end{scope}
}

\foreach \a in {0,1}
{
\draw
	(B1X\a) -- (C1X\a) -- (B2X\a) -- (C2X\a) -- (B3X\a)
	(A2X\a) -- (B2X\a)
	(A1X\a) -- (A5X\a) -- (A4X\a) -- (A3X\a)
	(A4X\a) -- (B4X\a) -- (C3X\a)
	(A5X\a) -- (B5X\a) -- (C5X\a)
	;

\draw[line width=1.8]
	(A1X\a) -- (B1X\a) 
	(A2X\a) -- (A3X\a)
	(B3X\a) -- (C3X\a)
	;

\draw[dashed]
	(A1X\a) -- (A2X\a)
	(A3X\a) -- (B3X\a)
	(B1X\a) -- (C5X\a)
	;
}

\draw[line width=1.8]
	(B4X0) -- (C4X0)
	(B5X1) -- (C4X1)
	;

\draw[dashed]
	(B5X0) -- (C4X0)
	(B4X1) -- (C4X1)
	;

\node at (0,0) {\small I};
\node at (3,0) {\small II};

\end{tikzpicture}
\caption{Propagation of neighborhood for the edge combination $a^3bc$.}
\label{pnd3abc}
\end{figure}

For the edge combination $a^3bc$, we call the two neighborhood tilings in Figure \ref{nd3abc} types I and II. An immediate consequence of Lemma \ref{edge_nd} is that, in either neighborhood tiling, each of the tiles labeled $\times$ in Figure \ref{pnd3abc} must have at least one vertex of degree $>3$. The labels I and II for $P_1$ and $P_2$ mean that theses tiles can have all vertices having degree $3$, without causing contradiction.

\begin{figure}[htp]
\centering
\begin{tikzpicture}[>=latex]

\foreach \x in {1,...,5}
\foreach \a in {0,1,2}
{
\begin{scope}[xshift=3*\a cm]

\draw[rotate=72*\x]
	(90:0.5) -- (18:0.5) -- (18:1) -- (54:1.3) -- (90:1);

\coordinate (A\x X\a) at (-54+72*\x:0.5);
\coordinate (B\x X\a) at (-54+72*\x:1);
\coordinate (C\x X\a) at (-18+72*\x:1.3);
\coordinate (P\x X\a) at (-18+72*\x:0.8);

\end{scope}
}

\draw[line width=1.8]
	(B1X0) -- (A1X0) -- (A2X0) -- (A3X0) -- (B3X0) 
	(B4X0) -- (C4X0) -- (B5X0)
	(B1X1) -- (C5X1) -- (B5X1) -- (C4X1) -- (B4X1) -- (C3X1) -- (B3X1) 
	(A1X1) -- (A2X1) -- (A3X1)
	(A2X1) -- (B2X1)
	(B1X2) -- (C5X2) -- (B5X2)  
	(C3X2) -- (B4X2) -- (C4X2)
	(A4X2) -- (B4X2)
	(A1X2) -- (A2X2) -- (A3X2)
	(A2X2) -- (B2X2);

\node at (0,0) {\small I};
\node at (3,0) {\small II};
\node at (6,0) {\small III};

\node at (P1X0) {\small I};
\node at (P2X0) {\small I};
\node at (P3X0) {\small $\times$};
\node at (P4X0) {\small II};
\node at (P5X0) {\small $\times$};

\node at (P1X1) {\small III};
\node at (P2X1) {\small III};
\node at (P3X1) {\small $\times$};
\node at (P4X1) {\small I};
\node at (P5X1) {\small $\times$};

\node at (P1X2) {\small III};
\node[scale=0.8] at (P2X2) {\small II,III};
\node at (P3X2) {\small III};
\node at (P4X2) {\small III};
\node at (P5X2) {\small III};

\node[xshift=6cm,fill=white,inner sep=1,draw,shape=circle] at (-54:0.75) {\tiny $a$};
\node[xshift=6cm,fill=white,inner sep=1,draw,shape=circle] at (234:0.75) {\tiny $b$};

\end{tikzpicture}
\caption{Propagation of neighborhood for the edge combination $a^3b^2$.}
\label{pnd3a2b}
\end{figure}
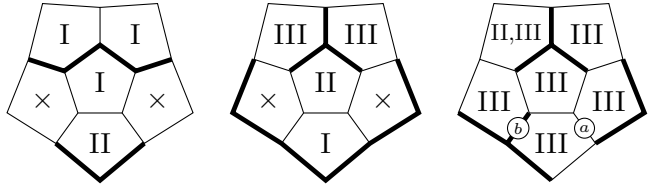

Similarly, for the edge combination $a^3b^2$, we call the three neighborhood tilings in Figure \ref{nd3a2b} types I, II, III. Again Lemma \ref{edge_nd} implies that, if a neighborhood tiling is of certain type, then each of the tiles labeled $\times$ in Figure \ref{pnd3a2b} must have at least one vertex of degree $>3$. Moreover, for those nearby tiles allowing all vertices having degree $3$, the neighborhood of that nearby tile must be of the indicated type. We also note that for the type III neighborhood, we highlight two edges with circles. These two edges will be used in future arguments.

The possible propagations for the edge combination $a^4b$ is summarized in Table \ref{pnd4ab}. The tilings in Figure \ref{nd4ab} are labeled from $1$ to $18$ instead of the Roman numerics.

\begin{table}[htp]
\begin{tabular}{|c||c|c|c|c|c|}
\hline 
Types & $P_1$ & $P_2$ & $P_3$ & $P_4$ & $P_5$ \\
\hline \hline
1 & 7,8,9,10,12 & 7,8,9,10,12 & 1,5,6,7,11 & 4,9,10,16 & 1,5,6,7,11\\
\hline 
2 & 2,3 & 2,3 & 2,3 & 2,14,17,18 & 2,3 \\
\hline 
3 & 2,3 & 6,8,14,15,17 & 5,6 & 6,11 & 2,3 \\
\hline 
4 & 6,8,14,15,17 & 6,8,14,15,17 & 5,6 & 1 & 5,6 \\
\hline 
5 & 3,4 & 8,14 & 10,11 & 5,7 & 1,5,6,7,11\\
\hline
6 &  3,4 & 6 & 3,4 & 3,8,12,13,15 & 1,5,6,7,11\\
\hline 
7 & 3,5,9,15,16 & 8,9 & 1 & 5,7 & 1,5,6,7,11\\
\hline 
8 & 7 & 5 & 3,4 & 3,8,12 & 1 \\
\hline
9 & 7 & 9,13 & 10,11 & 1 & 1 \\
\hline
10 & 10,12 & 10,12 & 1 & 1 & 1 \\
\hline
11 & 3,5,9,15,16 & $\times$ & $\times$ & 3,8,12,13,15 & 1,5,6,7,11\\
\hline 
12 & 10,12 & $\times$ & $\times$ & 6,11 & 1 \\
\hline
13 & 9,13 & $\times$ & $\times$ & 6,11 & 3,4 \\
\hline
14 & 5 & $\times$ & $\times$ & 14,15,17,18 & 3,4 \\
\hline
15 & $\times$ & $\times$ & 3,4 & 6,11 & 3,4 \\
\hline
16 & $\times$ & $\times$ & 3,4 & 1 & 3,4 \\
\hline
17 & $\times$ & $\times$ & 3,4 & 14,15,17,18 & 3,4 \\
\hline
18 & $\times$ & $\times$ & $\times$ & 14,15,17,18 & $\times$ \\
\hline
\end{tabular}
\caption{Propagation of neighborhood for the edge combination $a^4b$.}
\label{pnd4ab}
\end{table}

\section{Earth Map Tiling and Edge Length Combination $a^3bc$}
\label{section_3abc}

Earth map tilings are introduced in \cite{yan} as the pentagonal tilings of the sphere with exactly two vertices of degree $>3$, which we call {\em poles}. The concept is purely combinatorial, and ignores the edge lengths and angles.

The shortest paths between the poles are the {\em meridians}. In \cite{yan}, we showed that the {\em distance} between the two poles (i.e., the number of edges in the meridians) must be $\le 5$. Moreover, for each distance, the tiling is  given by repeating a combinatorially unique {\em timezone} bounded on both sides by meridians. The repetition will come back and be glued to the original timezone. This forms a cylinder. Then we add two poles to the top and the bottom of the cylinder to form the earth. We note that the two poles have the same degree, and for distances $4,3,2,1$, the degree is a multiple of $3$.

\begin{figure}[htp]
\centering
\begin{tikzpicture}[>=latex]


\foreach \x in {0,...,5}
{

\coordinate  (A\x) at (-0.3+1.2*\x,1.2);
\coordinate  (B\x) at (-0.3+1.2*\x,0.7);
\coordinate  (C\x) at (-0.6+1.2*\x,0.3);
\coordinate  (D\x) at (-0.6+1.2*\x,-0.3);
\coordinate  (E\x) at (-0.9+1.2*\x,-0.7);
\coordinate  (F\x) at (-0.9+1.2*\x,-1.2);
\coordinate  (G\x) at (1.2*\x,0.3);
\coordinate  (H\x) at (1.2*\x,-0.3);
	
}

\fill[gray!30]
	(B0) -- (C0) -- (D0) -- (H0) -- (E1) -- (D1) -- (C1) -- (G0) -- (B0)
	;
	
\foreach \x in {0,1}
\draw
	(A\x) -- (B\x) -- (C\x) -- (D\x) -- (E\x) -- (F\x)
	;
\draw
	(G0) -- (H0)
	(B0) -- (G0) -- (C1)
	(D0) -- (H0) -- (E1)
	;

\foreach \x in {2,3}
\draw[gray!70]
	(A\x) -- (B\x) -- (C\x) -- (D\x) -- (E\x) -- (F\x)
	;
\foreach \x in {1,2}
\draw[gray!70]
	(B\x) -- (G\x) -- (H\x) -- (D\x)
	(1.2*\x,0.3) -- ++(0.6,0)
	(1.2*\x,-0.3) -- ++(0.3,-0.4)
	;
	
\draw[<->]
	(-0.9,-1.3) -- node[below=-2] {\small timezone} ++(1.2,0);


\begin{scope}[xshift=7cm]

\fill[gray!30]
	(0,0.8) -- (0.4,0.5) -- (0.3,0) -- (0.4,-0.5) -- (0,-0.8) -- (-0.4,-0.5) -- (-0.3,0) -- (-0.4,0.5) -- (0,0.8)
	;

\foreach \x in {-1,1}
\foreach \y in {-1,1}
\draw[xscale=\x,yscale=\y]
	(0,1.2) -- (0,0.8) -- (0.4,0.5) -- (0.3,0)
	(0.4,0.5) -- (0.9,0.4) -- (1.2,0.7)
	(1.2,1.2) -- (1.2,0.7) -- (1.4,0) 
	(0.9,0.4) -- (0.9,0);
\draw
	(-0.3,0) -- (0.3,0)
	(-1.4,0) -- (-1.7,0)
	(-1.9,1.2) -- (-1.9,0.7) -- (-1.7,0) 
	(-1.9,-1.2) -- (-1.9,-0.7) -- (-1.7,0);
\draw[gray!70]
	(1.4,0) -- (1.7,0)
	(1.9,1.2) -- (1.9,0.7) -- (1.7,0) 
	(1.9,-1.2) -- (1.9,-0.7) -- (1.7,0);

\draw[<->]
	(-1.9,-1.3) -- node[below=-2] {\small timezone} ++(3.1,0);

\draw[<->]
	(-1.9,1.5) -- node[above] {\tiny meridian} node[below=-2] {\tiny part} ++(0.7,0);

\draw[<->]
	(-1.2,1.5) -- node[above] {\tiny core} node[below=-2] {\tiny part} ++(2.4,0);
		
\end{scope}


\begin{scope}[shift={(0cm,-3cm)}]

\foreach \x in {-1,1}
\fill[scale=\x,gray!30]
	(0.16,0.6) -- (-0.48,0.6) -- (-0.48,0.16) -- (-0.16,-0.16) -- (0.16,0.16)
	;

\foreach \x in {-1,1}
\draw[scale=\x]
	(0,0) -- (0.16,0.16) -- (0.48,-0.16) -- (0.8,0.16) -- (1.12,-0.16) -- (1.44,0.16) -- (1.76,-0.16)
	(0.64,1.2) -- (0.64,0.7) -- (0.16,0.6) -- (-0.48,0.6) -- (-0.96,0.7) -- (-0.96,1.2)
	(0.16,0.6) -- (0.16,0.16)
	(0.48,-0.6) -- (0.48,-0.16)
	(0.64,0.7) -- (0.8,0.16)
	(0.96,-0.7) -- (1.12,-0.16)
	(1.44,1.2) -- (1.44,0.16)
	(1.76,-1.2) -- (1.76,-0.16)
	;

\end{scope}


\begin{scope}[shift={(4.1cm,-3cm)}]

\fill[gray!30]
	(0,0.2) -- (0.4,0.4) -- (0.6,0) -- (0.4,-0.4) -- (0,-0.2) -- (0,0.2);

\foreach \y in {-1,1}
\draw[yscale=\y]
	(0,0.2) -- (-0.4,0.4) -- (-0.1,0.8) -- (0.4,0.7) -- (0.4,0.4) -- cycle
	(-0.1,0.8) -- (-0.1,1.2)
	(0.4,0.4) -- (0.6,0) -- (0.9,0) -- (1,0.5)
	(0.4,0.7) -- (1,0.5) -- (1.3,0.6) -- (1.4,1.2)
	(1.7,1.2) -- (1.8,0)
	(-0.4,0.4) -- (-0.6,0) -- (-0.9,0) -- (-1,1.2);

\draw
	(0,0.2) -- (0,-0.2)
	(1.3,0.6) -- (1.3,-0.6);

\end{scope}


\begin{scope}[shift={(8.5cm,-3cm)}]

\fill[gray!30]
	(0,0.2) -- (0.4,0.4) -- (0.6,0) -- (0.4,-0.4) -- (0,-0.2) -- (-0.4,-0.4) -- (-0.6,0) -- (-0.4,0.4) -- (0,0.2);

\foreach \x in {-1,1}
\foreach \y in {-1,1}
\draw[xscale=\x,yscale=\y]
	(0,0) -- (0,0.2) -- (0.4,0.4) -- (0.6,0) -- (0.9,0) -- (1,0.8) -- (1.1,1.2)
	(0.4,0.4) -- (0.4,0.7) 
	(0,0.7) -- (0.4,0.7) -- (1,0.8)
	(1.4,0) -- (1.4,1.2);

\end{scope}

\end{tikzpicture}
\caption{The timezone of earth map tiling and the core tile.}
\label{emt}
\end{figure}
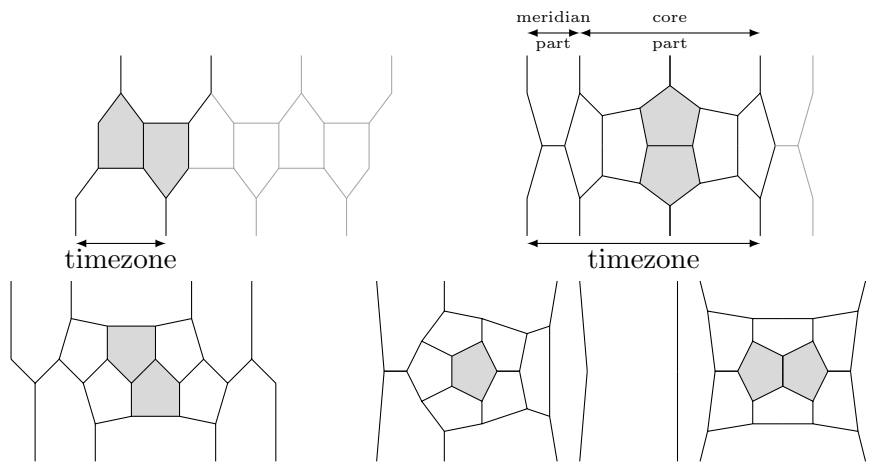

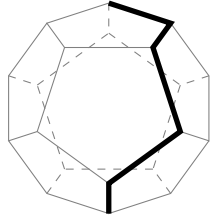
\begin{figure}[htp]
\centering
\begin{tikzpicture}[>=latex]

\foreach \x in {0,...,4}
{
\draw[gray]
	(54+72*\x:1) -- (-18+72*\x:1) -- (-18+72*\x:1.4) -- (18+72*\x:1.4) -- (54+72*\x:1.4);
\draw[gray,dashed]
	(90+72*\x:1) -- (18+72*\x:1) -- (18+72*\x:1.4);
}

\draw[line width=2]
	(90:1.4) -- (54:1.4) -- (54:1) -- (-18:1) -- (-90:1) -- (-90:1.4);

\end{tikzpicture}
\caption{Geodesic in the dodecahedron.}
\label{geodesic}
\end{figure}

The timezones are given by the black part of Figure \ref{emt}. They are obtained by cutting along any segment of the unique (up to symmetries) geodesic of length $5$ of the dodecahedron in Figure \ref{geodesic}. Note that if we cut along the whole geodesic, then we actually get three timezones glued together. Therefore we have only four tiles in the timezone for distance $5$. We also note that the timezone for distance $4$ consists of the {\em meridian part} and the {\em core part}. The meridian part can be considered as either on the left or the right of the core part.

We will apply Lemma \ref{edge_nd} to the gray tiles in Figure \ref{emt}, which we call the {\em core tiles}. For example, we observe that the core tiles cannot have three successive nearby tiles, such that each has at least one vertex of degree $>3$. Therefore the core tiles cannot have edge congruent neighborhood tilings like the ones in Figure \ref{pnd3abc}. This implies the following result.

\begin{theorem}\label{theorem3abc}
There are no edge congruent earth map tilings for the edge length combination $a^3bc$.
\end{theorem}

The earth map tilings have various symmetries. We note that the $180^{\circ}$ rotation is the composition of the vertical and horizontal flippings.
\begin{itemize}
\item For distance $5$, the timezone is symmetric with respect to the $180^{\circ}$ rotation, and the tiling is symmetric with respect to the vertical and horizontal flippings. 
\item For distance $4$, the core, the meridian part, and the tiling are all symmetric with respect to the vertical and horizontal flippings. 
\item For distance $3$, the timezone and the tiling are symmetric with respect to the $180^{\circ}$ rotation.
\item For distance $2$, the timezone and the tiling are symmetric with respect to the vertical flipping.
\item For distance $1$, the timezone and the tiling are symmetric with respect to the vertical and horizontal flippings. 
\end{itemize}
In the subsequent arguments, when we say up to symmetry, we mean up to all the symmetries above.

In general, we can always take a simple path in any spherical pentagonal tiling and cut the sphere open along the path. We get a timezone tiling with matching left and right sides. Then we may repeat this timezone tiling and get a generalized version of the earth map tiling. The difference from our specific earth map tilings is that we may have more vertices of degree $>3$ other than the poles. 

Given a particular generalized earth map tiling, we may try to list all the possible edge congruent tilings of the timezone. Then we may glue those timezone tilings with the matching boundaries. This is how we classify edge congruent earth map tilings.

\section{Earth Map Tiling and Edge Length Combination $a^2b^2c$}
\label{section_2a2bc}

In this section, we classify edge congruent earth map tilings for the edge combination $a^2b^2c$. We know that the edge lengths of each tile must be arranged in the order of $a,a,b,b,c$. However, for a specific pentagonal tile, there are ten ways of arranging the edge lengths in such an order, given by five rotations and five flippings. For distances $\ne 3$, due to the flipping symmetries of the earth map tilings, we actually only need to consider five ways of arranging the edge lengths. For distance $3$, we need to consider all ten ways. Our argument starts by considering the five (or ten for the distance $3$ case) possible ways the edge lengths of the core tile can be arranged.

The classification is summarized below. 

\begin{theorem}\label{theorem2a2bc}
The edge congruent earth map tilings for the edge length combination $a^2b^2c$ can be classified into the following numbers of families.
\begin{itemize}
\item Distance $5$: $2$ families.
\item Distance $4$: $2$ families.
\item Distance $3$: $3$ families.
\item Distance $2$: $3$ families.
\item Distance $1$: $2$ families.
\end{itemize}
\end{theorem}

The classification is up to the symmetries of earth map tilings given in Section \ref{section_3abc}, and the symmetry of exchanging $a$ and $b$. Moreover, it turns out that the families can be uniquely characterized by the distance between the poles and the edge length combinations at the poles.

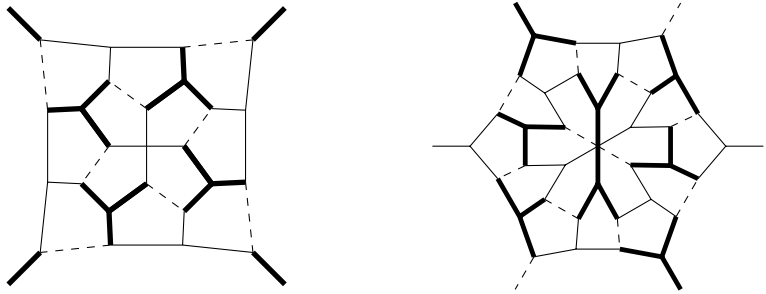
\begin{figure}[htp]
\centering
\begin{tikzpicture}[>=latex]


\foreach \x in {0,...,3}
{
\draw[rotate=90*\x]
	(0,0) -- (0.5,0)
	(20:1.4) -- (30:1)
	(45:2) -- (20:1.4) -- (-20:1.4)
	;

\draw[line width=2,rotate=90*\x]
	(0.5,0) -- (-30:1)
	(30:1) -- (60:1) -- (70:1.4)
	(0,0.5) -- (60:1)
	(45:2) -- (45:2.6)
	;

\draw[dashed,rotate=90*\x]
	(0.5,0) -- (30:1)
	(70:1.4) -- (45:2);

}


\begin{scope}[xshift=6cm]

\foreach \x in {0,...,5}
\foreach \y in {-1,1}
{

\coordinate (A\x) at (30+60*\x:0.5);
\coordinate (B\x X\y) at (15*\y+60*\x:1);
\coordinate (C\x X\y) at (18*\y+60*\x:1.4);
\coordinate (D\x) at (60*\x:1.7);
\coordinate (E\x) at (60*\x:2.2);

}

\draw
	(0,0) -- (A0)
	(0,0) -- (A3)
	(B0X1) -- (A0) -- (B1X-1)
	(B3X1) -- (A3) -- (B4X-1)
	(B2X-1) -- (B2X1) -- (A2)
	(B5X-1) -- (B5X1) -- (A5)
	(C0X-1) -- (D0) -- (C0X1)
	(C3X-1) -- (D3) -- (C3X1)
	(B1X1) -- (C1X1)
	(B4X1) -- (C4X1)
	(C2X-1) -- (C1X1) -- (D1)
	(C5X-1) -- (C4X1) -- (D4)
	(B2X1) -- (C2X1)
	(B5X1) -- (C5X1)
	(D0) -- (E0)
	(D3) -- (E3)
	;
\draw[line width=1.8]
	(0,0) -- (A1)
	(0,0) -- (A4)
	(A5) -- (B0X-1) -- (B0X1)
	(A2) -- (B3X-1) -- (B3X1)
	(B1X1) -- (A1) -- (B2X-1)
	(B4X1) -- (A4) -- (B5X-1)
	(B0X-1) -- (C0X-1)
	(B3X-1) -- (C3X-1)
	(B1X-1) -- (C1X-1)
	(B4X-1) -- (C4X-1)
	(C0X1) -- (C1X-1) -- (D1)
	(C3X1) -- (C4X-1) -- (D4)
	(C2X-1) -- (D2) -- (C2X1)
	(C5X-1) -- (D5) -- (C5X1)
	(D2) -- (E2)
	(D5) -- (E5)
	;
\draw[dashed]
	(0,0) -- (A5)
	(0,0) -- (A2)
	(B1X-1) -- (B1X1)
	(B4X-1) -- (B4X1)
	(B0X1) -- (C0X1)
	(B3X1) -- (C3X1)
	(B2X-1) -- (C2X-1)
	(B5X-1) -- (C5X-1)
	(C2X1) -- (C3X-1)
	(C5X1) -- (C0X-1)
	(D1) -- (E1)
	(D4) -- (E4)
	;

\end{scope}

\end{tikzpicture}
\caption{Earth map tilings $\binom{a}{b}^4$ and $\binom{bac}{bca}^2$, of distance $5$.}
\label{tile2a2bcB}
\end{figure}

Figure \ref{tile2a2bcB} gives one example each for the two families in distance $5$.

\subsection*{Distance $5$}

\begin{figure}[htp]
\centering
\begin{tikzpicture}[>=latex]

\foreach \x in {0,...,9}
{

\coordinate  (A\x) at (-0.3+1.2*\x,1.2);
\coordinate  (B\x) at (-0.3+1.2*\x,0.7);
\coordinate  (C\x) at (-0.6+1.2*\x,0.3);
\coordinate  (D\x) at (-0.6+1.2*\x,-0.3);
\coordinate  (E\x) at (-0.9+1.2*\x,-0.7);
\coordinate  (F\x) at (-0.9+1.2*\x,-1.2);
\coordinate  (G\x) at (1.2*\x,0.3);
\coordinate  (H\x) at (1.2*\x,-0.3);

\coordinate  (P\x) at (-0.3+1.2*\x,0.1);
\coordinate  (Q\x) at (0.3+1.2*\x,-0.1);
	
}

\node at (P0) {\small $P_1$};
\node at (Q0) {\small $P_2$};
\node[gray!70] at (P1) {\small $P_1$};
\node[gray!70] at (Q1) {\small $P_2$};
\node at (P4) {\small $P_3$};
\node at (Q4) {\small $P_4$};
\node at (P5) {\small $P_5$};
\node at (Q5) {\small $P_5$};
\node at (P6) {\small $P_4$};
\node at (Q6) {\small $P_3$};
\node[gray!70] at (P7) {\small $P_3$};
\node[gray!70] at (Q7) {\small $P_4$};
\node[gray!70] at (P8) {\small $P_5$};

\draw
	(A0) -- (B0) 
	(C0) -- (D0) -- (E0)
	(D0) -- (H0)
	(A1) -- (B1)
	(C1) -- (D1) -- (E1)
	(G4) -- (H4) 
	(D4) -- (H4) -- (E5)
	(A5) -- (B5) -- (C5)
	(B5) -- (G5)
	(H5) -- (E6) -- (F6)
	(D6) -- (E6)
	(B6) -- (G6) -- (H6)
	(G6) -- (C7)
	;

\draw[line width=2]
	(E0) -- (F0) 
	(B0) -- (G0) -- (H0)
	(G0) -- (C1)
	(E1) -- (F1)
	(A4) -- (B4) -- (C4) 
	(B4) -- (G4)
	(D4) -- (E4) -- (F4)
	(C5) -- (D5) -- (E5)
	(D5) -- (H5)
	(B6) -- (C6) -- (D6)
	(G5) -- (C6)
	(H6) -- (E7) -- (F7)
	(D7) -- (E7)
	(A7) -- (B7) -- (C7)
	;

\draw[dashed]
	(B0) -- (C0) 
	(H0) -- (E1)
	(B1) -- (C1)
	(C4) -- (D4) 
	(E5) -- (F5)
	(G4) -- (C5)
	(G5) -- (H5)
	(A6) -- (B6)
	(D6) -- (H6)
	(C7) -- (D7)
	;

\draw[gray!70]
	(D1) -- (H1) 
	(A2) -- (B2)
	(C2) -- (D2) -- (E2)
	(D7) -- (H7) -- (E8) 
	(G7) -- (H7)
	(A8) -- (B8) -- (C8)
	(B8) -- (G8)
	(H8) -- (E9) -- (F9)
	;

\draw[gray!70,line width=2]
	(B1) -- (G1) -- (H1) 
	(G1) -- (C2)
	(E2) -- (F2)
	(B7) -- (G7) 
	(C8) -- (D8) -- (E8)
	(D8) -- (H8)
	;

\draw[gray!70,dashed]
	(H1) -- (E2) 
	(B2) -- (C2)
	(G7) -- (C8) 
	(E8) -- (F8)
	(G8) -- (H8)
	;

\node[fill=white,inner sep=1,draw,shape=circle] at (A0) {\small $a$};

\node[fill=white,inner sep=0.3,draw,shape=circle] at (F0) {\small $b$};

\node[fill=white,inner sep=0.3,draw,shape=circle] at (A4) {\small $b$};
\node[fill=white,inner sep=1,draw,shape=circle] at (A5) {\small $a$};	
\node[fill=white,inner sep=1,draw,shape=circle] at (A6) {\small $c$};

\node[fill=white,inner sep=0.3,draw,shape=circle] at (F4) {\small $b$};
\node[fill=white,inner sep=1,draw,shape=circle] at (F5) {\small $c$};	
\node[fill=white,inner sep=1,draw,shape=circle] at (F6) {\small $a$};

\end{tikzpicture}
\caption{Earth map tiling, distance $5$, edge combination $a^2b^2c$.}
\label{d5_case1}
\end{figure}
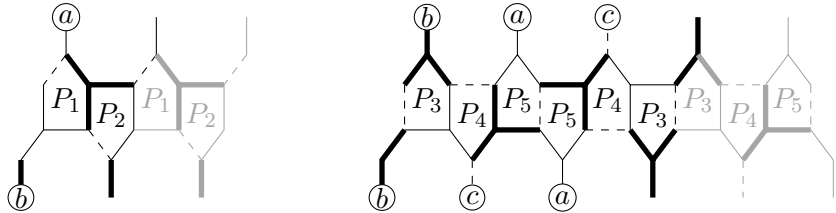

Up to symmetry, the tiles $P_1,P_2,P_3,P_4,P_5$ in Figure \ref{d5_case1} give all the possible edge length arrangements of the core tile that conform to the order of $a,a,b,b,c$. By Lemma \ref{edge_nd}, once we fix one such tile $P$, all the edges in the neighborhood of $P$ can be determined. In particular, all the edges of the core tiles $P_-$ and $P_+$ on the left and right of $P$ are determined. Then Lemma \ref{edge_nd} further determines all the edges in the neighborhoods of $P_-$ and $P_+$. The process continues, so that the edges of the starting tile uniquely determine all the other edges.

We get two repeating patterns, which give two families of edge congruent earth map tilings of distance $5$ for the edge combination $a^2b^2c$.
\begin{enumerate}
\item Repeat the black part of the first tiling in Figure \ref{d5_case1} by $k$ times. The pole vertices are $a^k$ and $b^k$.
\item Repeat the black part of the second tiling in Figure \ref{d5_case1} by $k$ times. Both pole vertices are $(abc)^k$. 
\end{enumerate}
Note that the black part of the second tiling is obtained by gluing the neighborhood tiling of $P_4$ and its $180^{\circ}$ rotation. 

We introduce a notation for the tilings of the black parts in Figure \ref{d5_case1}, by taking the circled edges at both poles and reading the edge lengths from left to right. So we denote the first black part tiling by $\binom{a}{b}$, and denote the second black part tiling by $\binom{bac}{bca}$. The two families can be described as follows.
\begin{enumerate}
\item Circular product of $\binom{a}{b}$.
\item Circular product of $\binom{bac}{bca}$.
\end{enumerate}

\subsection*{Distance $4$}
 
\begin{figure}[htp]
\centering
\begin{tikzpicture}[>=latex]

\foreach \x in {-1,1}
\foreach \y in {-1,1}
\foreach \a in {0,1,2}
{
\begin{scope}[xshift=4.8*\a cm]

\coordinate  (A\y X\a) at (0,1.2*\y);
\coordinate  (B\y X\a) at (0,0.8*\y);
\coordinate  (C\x X\a) at (0.3*\x,0);
\coordinate  (D\x\y X\a) at (0.4*\x,0.5*\y);
\coordinate  (E\x\y X\a) at (0.9*\x,0.4*\y);
\coordinate  (F\x\y X\a) at (1.2*\x,0.7*\y);
\coordinate  (G\x\y X\a) at (1.2*\x,1.2*\y);
\coordinate  (H\x X\a) at (1.4*\x,0);
\coordinate  (I\x X\a) at (1.7*\x,0);
\coordinate  (J\x\y X\a) at (1.9*\x,0.7*\y);
\coordinate  (K\x\y X\a) at (1.9*\x,1.2*\y);

\coordinate  (P\a) at (0,-1.5);
	
\end{scope}
}

\draw
	(A1X0) -- (B1X0) 
	(D11X0) -- (C1X0) -- (D1-1X0)
	(C1X0) -- (C-1X0)
	(G11X0) -- (F11X0)
	(G-11X0) -- (F-11X0)
	(H1X0) -- (F1-1X0) -- (G1-1X0)
	(E1-1X0) -- (F1-1X0)
	(E-11X0) -- (E-1-1X0) -- (D-1-1X0)
	(E-1-1X0) -- (F-1-1X0)
	(J-11X0) -- (K-11X0)
	(I-1X0) -- (J-1-1X0) -- (K-1-1X0)
	(K11X0) -- (J11X0)
	(D11X1) -- (E11X1) -- (F11X1)  
	(E11X1) -- (E1-1X1)
	(D-11X1) -- (C-1X1)
	(B1X1) -- (D-11X1) -- (E-11X1)
	(D1-1X1) -- (B-1X1) -- (D-1-1X1)
	(B-1X1) -- (A-1X1)
	(F-11X1) -- (H-1X1) -- (F-1-1X1)
	(H-1X1) -- (I-1X1)
	(H1X1) -- (I1X1)
	(J11X1) -- (I1X1) -- (J1-1X1) 
	(B1X2) -- (D-11X2) -- (E-11X2) 
	(D-11X2) -- (C-1X2)
	(H-1X2) -- (F-1-1X2) -- (G-1-1X2)
	(I-1X2) -- (J-11X2) -- (K-11X2)
	(B-1X2) -- (D1-1X2) -- (E1-1X2)
	(D1-1X2) -- (C1X2)
	(H1X2) -- (F11X2) -- (G11X2)
	(I1X2) -- (J1-1X2) -- (K1-1X2)
	(E11X2) -- (F11X2)
	(E-1-1X2) -- (F-1-1X2)
	;
	
\draw[line width=2]
	(D11X0) -- (E11X0) -- (F11X0)  
	(E11X0) -- (E1-1X0)
	(D-11X0) -- (C-1X0)
	(B1X0) -- (D-11X0) -- (E-11X0)
	(D1-1X0) -- (B-1X0) -- (D-1-1X0)
	(B-1X0) -- (A-1X0)
	(F-11X0) -- (H-1X0) -- (F-1-1X0)
	(H-1X0) -- (I-1X0)
	(H1X0) -- (I1X0)
	(J11X0) -- (I1X0) -- (J1-1X0)
	(A1X1) -- (B1X1) 
	(D11X1) -- (C1X1) -- (D1-1X1)
	(C1X1) -- (C-1X1)
	(G11X1) -- (F11X1)
	(H1X1) -- (F1-1X1) -- (G1-1X1)
	(E1-1X1) -- (F1-1X1)
	(E-11X1) -- (E-1-1X1) -- (D-1-1X1)
	(E-1-1X1) -- (F-1-1X1)
	(K11X1) -- (J11X1)
	(G-11X1) -- (F-11X1)
	(K-11X1) -- (J-11X1)
	(I-1X1) -- (J-1-1X1) -- (K-1-1X1)
	(K11X1) -- (J11X1)
	(B1X2) -- (D11X2) -- (E11X2) 
	(D11X2) -- (C1X2)
	(H1X2) -- (F1-1X2) -- (G1-1X2)
	(I1X2) -- (J11X2) -- (K11X2)
	(B-1X2) -- (D-1-1X2) -- (E-1-1X2)
	(D-1-1X2) -- (C-1X2)
	(H-1X2) -- (F-11X2) -- (G-11X2)
	(I-1X2) -- (J-1-1X2) -- (K-1-1X2)
	(E1-1X2) -- (F1-1X2)
	(E-11X2) -- (F-11X2)
	;

\draw[dashed]
	(B1X0) -- (D11X0) 
	(F11X0) -- (H1X0)
	(D1-1X0) -- (E1-1X0)
	(J1-1X0) -- (K1-1X0)
	(E-11X0) -- (F-11X0)
	(F-1-1X0) -- (G-1-1X0)
	(I-1X0) -- (J-11X0)
	(D-1-1X0) -- (C-1X0)
	(B1X1) -- (D11X1) 
	(F11X1) -- (H1X1)
	(D1-1X1) -- (E1-1X1)
	(J1-1X1) -- (K1-1X1)
	(E-11X1) -- (F-11X1)
	(F-1-1X1) -- (G-1-1X1)
	(I-1X1) -- (J-11X1)
	(D-1-1X1) -- (C-1X1)
	(C1X2) -- (C-1X2) 
	(A1X2) -- (B1X2)
	(E11X2) -- (E1-1X2)
	(H1X2) -- (I1X2)
	(A-1X2) -- (B-1X2)
	(E-11X2) -- (E-1-1X2)
	(H-1X2) -- (I-1X2)
	;

\node[fill=white,inner sep=1,draw,shape=circle] at (G-11X0) {\small $a$};	
\node[fill=white,inner sep=1,draw,shape=circle] at (A1X0) {\small $a$};
\node[fill=white,inner sep=1,draw,shape=circle] at (G11X0) {\small $a$};

\node[fill=white,inner sep=1,draw,shape=circle] at (G-1-1X0) {\small $c$};	
\node[fill=white,inner sep=0.3,draw,shape=circle] at (A-1X0) {\small $b$};
\node[fill=white,inner sep=1,draw,shape=circle] at (G1-1X0) {\small $a$};

\end{tikzpicture}
\caption{Timezone tiling, distance $4$, edge combination $a^2b^2c$.}
\label{d4_case1}
\end{figure}
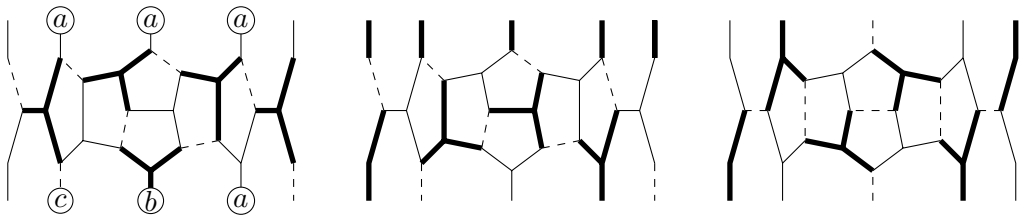

Up to symmetry, the core tiles in Figure \ref{d4_case1} give all five possible edge length arrangements. Then we may apply Lemma \ref{edge_nd} to get all the edges in the neighborhood of one core tile. If all vertices of any tile in this neighborhood have degree $3$, then we may further apply Lemma \ref{edge_nd} to this tile and get the edges of more tiles. Keep going, we see that the edges of any core tile uniquely determine all the other edges in the timezone (actually also the meridian parts on both sides).

By comparing the edge length arrangements of the meridians (or matching the meridian parts), we see that each timezone tiling can only repeat itself. Moreover, if a symmetry of the earth map tiling fixes the boundary meridians, then it also fixes the whole timezone tiling. Therefore the symmetries do not produce new timezone tilings. On the other hand, exchanging $a$ and $b$ means exchanging the first and the second tilings. After all the considerations, we get two families of edge congruent earth map tilings of distance $4$ for the edge combination $a^2b^2c$.
\begin{enumerate}
\item Circular product of $\binom{aaa}{cba}$.
\item Circular product of $\binom{bca}{acb}$.
\end{enumerate}
Note that the notation takes the ``middle three'' edges at the poles and reads from left to right. The choice respects the symmetry of the earth map tiling.

\subsection*{Distance $3$}

\begin{figure}[htp]
\centering
\begin{tikzpicture}[>=latex]

\foreach \a in {0,1,2}
\foreach \b in {0,1}
\foreach \x in {-1,1}
\foreach \y in {-1,1}
{
\begin{scope}[shift={(4.8*\a cm, -3*\b cm)}]

\coordinate  (A\x\y X\a\b) at (0.32*\x-0.16*\y,0.16*\y);
\coordinate  (B\x\y X\a\b) at (0.32*\x-0.16*\y,0.6*\y);
\coordinate  (C\x\y X\a\b) at (0.8*\x-0.16*\y,1.2*\y);
\coordinate  (D\x\y X\a\b) at (0.8*\x-0.16*\y,0.7*\y);
\coordinate  (E\x\y X\a\b) at (0.96*\x-0.16*\y,0.16*\y);
\coordinate  (F\x\y X\a\b) at (1.6*\x-0.16*\y,0.16*\y);
\coordinate  (G\x\y X\a\b) at (1.6*\x-0.16*\y,1.2*\y);
	
\end{scope}
}

\draw
	(A11X00) -- (A1-1X00) -- (E11X00) 
	(A1-1X00) -- (B1-1X00)
	(A-1-1X00) -- (A-11X00) -- (E-1-1X00)
	(A-11X00) -- (B-11X00)
	(C11X00) -- (D11X00) 
	(C1-1X00) -- (D1-1X00) 
	(C-11X00) -- (D-11X00) 
	(C-1-1X00) -- (D-1-1X00)
	(F11X00) -- (G11X00) 
	(F1-1X00) -- (G1-1X00) 
	(F-11X00) -- (G-11X00) 
	(F-1-1X00) -- (G-1-1X00)
	(A11X10) -- (A1-1X10) -- (E11X10) 
	(A1-1X10) -- (B1-1X10)
	(E1-1X10) -- (F11X10) -- (F1-1X10)
	(F11X10) -- (G11X10) 
	(B11X10) -- (B-11X10) -- (D-11X10) 
	(B-11X10) -- (A-11X10)
	(C-1-1X10) -- (D-1-1X10) -- (E-1-1X10)
	(D-1-1X10) -- (B-1-1X10)
	(G-11X10) -- (F-11X10) -- (F-1-1X10)
	(B-11X20) -- (B11X20) -- (D11X20) 
	(A11X20) -- (B11X20) 
	(B-1-1X20) -- (B1-1X20) -- (D1-1X20)
	(A1-1X20) -- (B1-1X20) 
	(A-11X20) -- (E-1-1X20) -- (E-11X20) 
	(E-1-1X20) -- (D-1-1X20) 
	(G11X20) -- (F11X20) -- (F1-1X20) 
	(E1-1X20) -- (F11X20)
	(G-11X20) -- (F-11X20) -- (F-1-1X20)	
	(B-11X01) -- (B11X01) -- (D11X01) 
	(B1-1X01) -- (B-1-1X01) -- (D-1-1X01)
	(A11X01) -- (B11X01) 
	(A-1-1X01) -- (B-1-1X01)
	(E11X01) -- (E1-1X01) -- (F11X01) 
	(E-1-1X01) -- (E-11X01) -- (F-1-1X01) 
	(D1-1X01) -- (E1-1X01)
	(D-11X01) -- (E-11X01)
	(C11X11) -- (D11X11) -- (E11X11) 
	(B11X11) -- (D11X11)
	(C1-1X11) -- (D1-1X11) -- (E1-1X11)
	(B1-1X11) -- (D1-1X11)
	(A-11X11) -- (A-1-1X11) -- (A11X11)
	(A-1-1X11) -- (B-1-1X11)
	(E-1-1X11) -- (E-11X11) -- (F-1-1X11)
	(D-11X11) -- (E-11X11)
	(A1-1X21) -- (E11X21) -- (E1-1X21) 
	(E11X21) -- (D11X21)
	(A-11X21) -- (A-1-1X21) -- (A11X21)
	(A-1-1X21) -- (B-1-1X21)
	(C-11X21) -- (D-11X21) -- (E-11X21)
	(B-11X21) -- (D-11X21)
	(C1-1X21) -- (D1-1X21)
	(C-1-1X21) -- (D-1-1X21)
	(G1-1X21) -- (F1-1X21)
	(G-1-1X21) -- (F-1-1X21)
	;

\draw[line width=2]
	(B-11X00) -- (B11X00) -- (D11X00) 
	(B1-1X00) -- (B-1-1X00) -- (D-1-1X00)
	(A11X00) -- (B11X00) 
	(A-1-1X00) -- (B-1-1X00)
	(E11X00) -- (E1-1X00) -- (F11X00) 
	(E-1-1X00) -- (E-11X00) -- (F-1-1X00) 
	(D1-1X00) -- (E1-1X00)
	(D-11X00) -- (E-11X00)
	(C11X10) -- (D11X10) -- (E11X10) 
	(B11X10) -- (D11X10)
	(C1-1X10) -- (D1-1X10) -- (E1-1X10)
	(B1-1X10) -- (D1-1X10)
	(A-11X10) -- (A-1-1X10) -- (A11X10)
	(A-1-1X10) -- (B-1-1X10)
	(E-1-1X10) -- (E-11X10) -- (F-1-1X10)
	(D-11X10) -- (E-11X10)
	(A1-1X20) -- (E11X20) -- (E1-1X20) 
	(E11X20) -- (D11X20)
	(A-11X20) -- (A-1-1X20) -- (A11X20)
	(A-1-1X20) -- (B-1-1X20)
	(C-11X20) -- (D-11X20) -- (E-11X20)
	(B-11X20) -- (D-11X20)
	(C1-1X20) -- (D1-1X20)
	(C-1-1X20) -- (D-1-1X20)
	(G1-1X20) -- (F1-1X20)
	(G-1-1X20) -- (F-1-1X20)
	(A11X01) -- (A1-1X01) -- (E11X01) 
	(A1-1X01) -- (B1-1X01)
	(A-1-1X01) -- (A-11X01) -- (E-1-1X01)
	(A-11X01) -- (B-11X01)
	(C11X01) -- (D11X01) 
	(C1-1X01) -- (D1-1X01) 
	(C-11X01) -- (D-11X01) 
	(C-1-1X01) -- (D-1-1X01)
	(F11X01) -- (G11X01) 
	(F1-1X01) -- (G1-1X01) 
	(F-11X01) -- (G-11X01) 
	(F-1-1X01) -- (G-1-1X01)
	(A11X11) -- (A1-1X11) -- (E11X11) 
	(A1-1X11) -- (B1-1X11)
	(E1-1X11) -- (F11X11) -- (F1-1X11)
	(F11X11) -- (G11X11) 
	(B11X11) -- (B-11X11) -- (D-11X11) 
	(B-11X11) -- (A-11X11)
	(C-1-1X11) -- (D-1-1X11) -- (E-1-1X11)
	(D-1-1X11) -- (B-1-1X11)
	(G-11X11) -- (F-11X11) -- (F-1-1X11)
	(B-11X21) -- (B11X21) -- (D11X21) 
	(A11X21) -- (B11X21) 
	(B-1-1X21) -- (B1-1X21) -- (D1-1X21)
	(A1-1X21) -- (B1-1X21) 
	(A-11X21) -- (E-1-1X21) -- (E-11X21) 
	(E-1-1X21) -- (D-1-1X21) 
	(G11X21) -- (F11X21) -- (F1-1X21) 
	(E1-1X21) -- (F11X21)
	(G-11X21) -- (F-11X21) -- (F-1-1X21)
	;

\foreach \a in {0,1}
\draw[dashed]
	(A-1-1X0\a) -- (A11X0\a) 
	(B-11X0\a) -- (D-11X0\a)
	(B1-1X0\a) -- (D1-1X0\a)
	(D11X0\a) -- (E11X0\a)
	(D-1-1X0\a) -- (E-1-1X0\a)
	(F11X0\a) -- (F1-1X0\a)
	(F-11X0\a) -- (F-1-1X0\a) 
	(E11X1\a) -- (E1-1X1\a) 
	(F1-1X1\a) -- (G1-1X1\a)
	(F-1-1X1\a) -- (G-1-1X1\a)
	(C-11X1\a) -- (D-11X1\a)
	(A11X1\a) -- (B11X1\a)
	(B-1-1X1\a) -- (B1-1X1\a)
	(A-11X1\a) -- (E-1-1X1\a)
	(A11X2\a) -- (A1-1X2\a) 
	(C11X2\a) -- (D11X2\a)
	(D1-1X2\a) -- (E1-1X2\a) 
	(A-11X2\a) -- (B-11X2\a) 
	(B-1-1X2\a) -- (D-1-1X2\a)
	(E-11X2\a) -- (F-1-1X2\a)
	;

\node[fill=white,inner sep=1,draw,shape=circle] at (G-11X00) {\small $a$};	
\node[fill=white,inner sep=1,draw,shape=circle] at (C-11X00) {\small $a$};
\node[fill=white,inner sep=1,draw,shape=circle] at (C11X00) {\small $a$};

\node[fill=white,inner sep=1,draw,shape=circle] at (C-1-1X00) {\small $a$};	
\node[fill=white,inner sep=1,draw,shape=circle] at (C1-1X00) {\small $a$};
\node[fill=white,inner sep=1,draw,shape=circle] at (G1-1X00) {\small $a$};

\end{tikzpicture}
\caption{Timezone tiling, distance $3$, edge combination $a^2b^2c$.}
\label{d3_case1}
\end{figure}
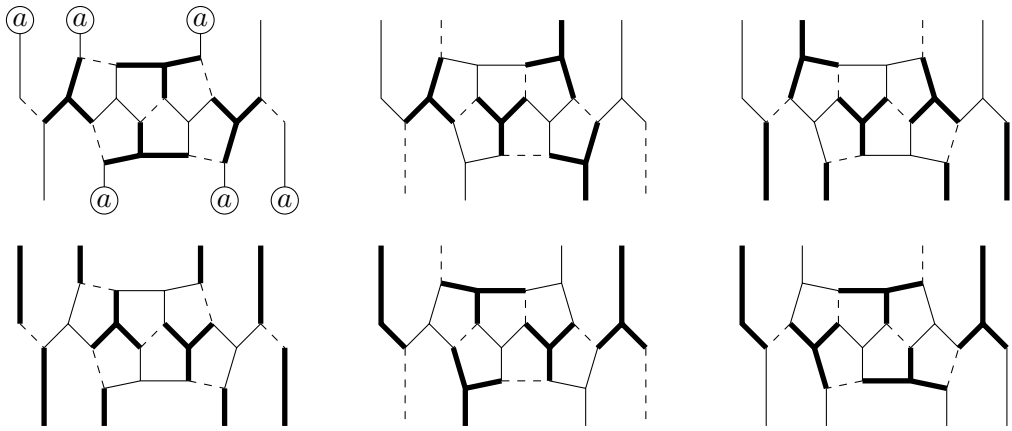

Up to symmetry, the core tiles in Figure \ref{d3_case1} give all ten possible edge length arrangements. Then Lemma \ref{edge_nd} determines all the edges in the timezone.

By comparing the edge length arrangements of the meridians, we see that each timezone tiling can only repeat itself. Moreover, the $180^{\circ}$ rotation does not produce new timezone tilings. On the other hand, exchanging $a$ and $b$ simply exchanges a tiling in the first row with the corresponding tiling in the second row. After all the considerations, we only need to repeat the three timezone tilings in the first rwo. This gives three families of edge congruent earth map tilings of distance $3$ for the edge combination $a^2b^2c$.
\begin{enumerate}
\item Circular product of $\binom{aaa}{aaa}$.
\item Circular product of $\binom{acb}{abc}$.
\item Circular product of $\binom{abc}{bbb}$.
\end{enumerate}
The choice of notation respects the symmetry of $180^{\circ}$ rotation.

\subsection*{Distance $2$}

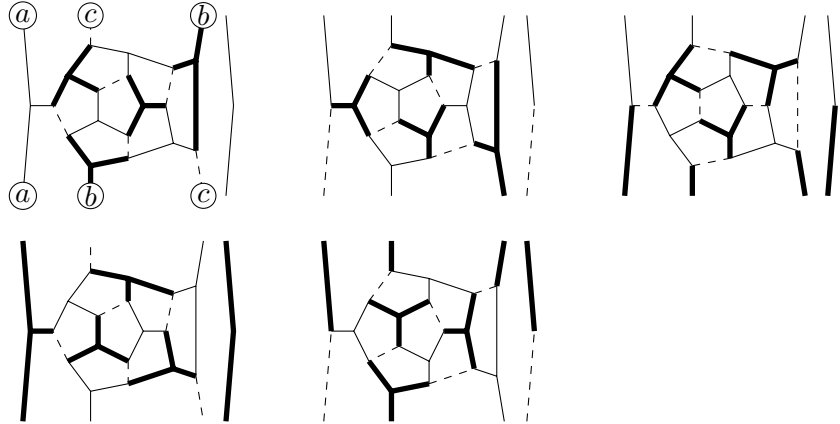
\begin{figure}[htp]
\centering
\begin{tikzpicture}[>=latex]

\foreach \a in {0,1,2}
\foreach \b in {0,1}
\foreach \y in {-1,1}
{
\begin{scope}[shift={(4*\a cm, -3*\b cm)}]

\coordinate  (A\y X\a\b) at (0,0.2*\y);
\coordinate  (B\y X\a\b) at (-0.4,0.4*\y);
\coordinate  (C\y X\a\b) at (0.4,0.4*\y);
\coordinate  (D\y X\a\b) at (-0.1,0.8*\y);
\coordinate  (E\y X\a\b) at (0.4,0.7*\y);
\coordinate  (F\y X\a\b) at (-0.1,1.2*\y);
\coordinate  (GX\a\b) at (0.6,0);
\coordinate  (HX\a\b) at (0.9,0);
\coordinate  (I\y X\a\b) at (1,0.5*\y);
\coordinate  (J\y X\a\b) at (1.3,0.6*\y);
\coordinate  (K\y X\a\b) at (1.4,1.2*\y);
\coordinate  (L\y X\a\b) at (1.7,1.2*\y);
\coordinate  (MX\a\b) at (1.8,0);
\coordinate  (NX\a\b) at (-0.6,0);
\coordinate  (OX\a\b) at (-0.9,0);
\coordinate  (P\y X\a\b) at (-1,1.2*\y);

\end{scope}
}

\draw
	(B-1X00) -- (A-1X00) -- (C-1X00)  
	(A1X00) -- (A-1X00)
	(D1X00) -- (E1X00) -- (I1X00) 
	(C1X00) -- (E1X00) 
	(E-1X00) -- (I-1X00) -- (J-1X00) 
	(HX00) -- (I-1X00)
	(L1X00) -- (MX00) -- (L-1X00)
	(P1X00) -- (OX00) -- (P-1X00)
	(OX00) -- (NX00)
	(B1X10) -- (A1X10) -- (C1X10)  
	(A1X10) -- (A-1X10)
	(B-1X10) -- (D-1X10) -- (E-1X10)
	(D-1X10) -- (F-1X10) 
	(I1X10) -- (HX10) -- (I-1X10) 
	(GX10) -- (HX10)
	(OX10) -- (P1X10)
	(D1X10) -- (F1X10) 
	(J1X10) -- (K1X10)
	(MX10) -- (L1X10)
	(A-1X20) -- (B-1X20) -- (D-1X20) 
	(B-1X20) -- (NX20)
	(A1X20) -- (C1X20) -- (GX20) 
	(C1X20) -- (E1X20) 
	(E-1X20) -- (I-1X20) -- (J-1X20) 
	(HX20) -- (I-1X20)
	(OX20) -- (P1X20)
	(D1X20) -- (F1X20) 
	(J1X20) -- (K1X20)
	(MX20) -- (L1X20)
	(A1X01) -- (B1X01) -- (D1X01) 
	(B1X01) -- (NX01)
	(C1X01) -- (GX01) -- (C-1X01) 
	(GX01) -- (HX01) 
	(B-1X01) -- (D-1X01) -- (E-1X01)
	(D-1X01) -- (F-1X01)
	(I1X01) -- (J1X01) -- (K1X01) 
	(J1X01) -- (J-1X01)
	(B1X11) -- (NX11) -- (B-1X11) 
	(OX11) -- (NX11)
	(A-1X11) -- (C-1X11) -- (GX11) 
	(C-1X11) -- (E-1X11) 
	(D1X11) -- (E1X11) -- (I1X11) 
	(C1X11) -- (E1X11)
	(J1X11) -- (J-1X11) -- (K-1X11)
	(I-1X11) -- (J-1X11)
	;

\draw[line width=2]
	(A1X00) -- (B1X00) -- (D1X00) 
	(B1X00) -- (NX00)
	(C1X00) -- (GX00) -- (C-1X00) 
	(GX00) -- (HX00) 
	(B-1X00) -- (D-1X00) -- (E-1X00)
	(D-1X00) -- (F-1X00)
	(I1X00) -- (J1X00) -- (K1X00) 
	(J1X00) -- (J-1X00)
	(B1X10) -- (NX10) -- (B-1X10) 
	(OX10) -- (NX10)
	(A-1X10) -- (C-1X10) -- (GX10) 
	(C-1X10) -- (E-1X10) 
	(D1X10) -- (E1X10) -- (I1X10) 
	(C1X10) -- (E1X10)
	(J1X10) -- (J-1X10) -- (K-1X10)
	(I-1X10) -- (J-1X10)
	(A1X20) -- (B1X20) -- (D1X20) 
	(B1X20) -- (NX20)
	(A-1X20) -- (C-1X20) -- (GX20) 
	(C-1X20) -- (E-1X20) 
	(E1X20) -- (I1X20) -- (J1X20) 
	(HX20) -- (I1X20)
	(OX20) -- (P-1X20)
	(D-1X20) -- (F-1X20) 
	(J-1X20) -- (K-1X20)
	(MX20) -- (L-1X20)
	(B-1X01) -- (A-1X01) -- (C-1X01)  
	(A1X01) -- (A-1X01)
	(D1X01) -- (E1X01) -- (I1X01) 
	(C1X01) -- (E1X01) 
	(E-1X01) -- (I-1X01) -- (J-1X01) 
	(HX01) -- (I-1X01)
	(L1X01) -- (MX01) -- (L-1X01)
	(P1X01) -- (OX01) -- (P-1X01)
	(OX01) -- (NX01)
	(B1X11) -- (A1X11) -- (C1X11)  
	(A1X11) -- (A-1X11)
	(B-1X11) -- (D-1X11) -- (E-1X11)
	(D-1X11) -- (F-1X11) 
	(I1X11) -- (HX11) -- (I-1X11) 
	(GX11) -- (HX11)
	(OX11) -- (P1X11)
	(D1X11) -- (F1X11) 
	(J1X11) -- (K1X11)
	(MX11) -- (L1X11)
	;

\foreach \a in {0,1}
\draw[dashed]
	(A1X0\a) -- (C1X0\a) 
	(D1X0\a) -- (F1X0\a)
	(B-1X0\a) -- (NX0\a)
	(C-1X0\a) -- (E-1X0\a)
	(HX0\a) -- (I1X0\a) 
	(J-1X0\a) -- (K-1X0\a)
	(B1X1\a) -- (D1X1\a) 
	(A-1X1\a) -- (B-1X1\a) 
	(P-1X1\a) -- (OX1\a) 
	(C1X1\a) -- (GX1\a) 
	(I1X1\a) -- (J1X1\a)
	(E-1X1\a) -- (I-1X1\a)
	(L-1X1\a) -- (MX1\a)
	(A1X2\a) -- (A-1X2\a) 
	(OX2\a) -- (NX2\a) 
	(D1X2\a) -- (E1X2\a) 
	(D-1X2\a) -- (E-1X2\a) 
	(GX2\a) -- (HX2\a) 
	(J1X2\a) -- (J-1X2\a)
	;
	
\fill[white]
	(6.9,-4.25) rectangle (9.9,-1.75);

\node[fill=white,inner sep=1,draw,shape=circle] at (P1X00) {\small $a$};	
\node[fill=white,inner sep=1,draw,shape=circle] at (F1X00) {\small $c$};
\node[fill=white,inner sep=0.3,draw,shape=circle] at (K1X00) {\small $b$};

\node[fill=white,inner sep=1,draw,shape=circle] at (P-1X00) {\small $a$};	
\node[fill=white,inner sep=0.3,draw,shape=circle] at (F-1X00) {\small $b$};
\node[fill=white,inner sep=1,draw,shape=circle] at (K-1X00) {\small $c$};

\end{tikzpicture}
\caption{Timezone tiling, distance $2$, edge combination $a^2b^2c$.}
\label{d2_case1}
\end{figure}

Up to symmetry, the core tiles in Figure \ref{d2_case1} give all five possible edge lengths arrangements. Then Lemma 1 determines all the edges in the timezone.

By comparing the edge length arrangements of the meridians, we see that each timezone tiling can only repeat itself. Moreover, the vertical flipping produces new timezone tilings when applied to the first and the fourth (i.e., the left two) tilings. However, exchanging $a$ and $b$ simply exchanges the tiling in the first row with the corresponding tiling in the second row, and leaves the third tiling invariant (up to symmetry). After all the considerations, we only need to repeat the three timezone tilings in the first row and the vertical flipping of the first tiling. This gives three families of edge congruent earth map tilings of distance $2$ for the edge combination $a^2b^2c$.
\begin{enumerate}
\item Circular product of $\binom{acb}{abc}$, $\binom{abc}{acb}$.
\item Circular product of $\binom{aaa}{cab}$.
\item Circular product of $\binom{aaa}{bbb}$.
\end{enumerate}
The vertical flipping of $\binom{acb}{abc}$ is $\binom{abc}{acb}$. The choice of notation respects the symmetry of vertical flipping.

\subsection*{Distance $1$}

Up to symmetry, the core tiles in Figure \ref{d2_case1} give all five possible edge length arrangements. Then Lemma 1 determines all the edges in the timezone. 

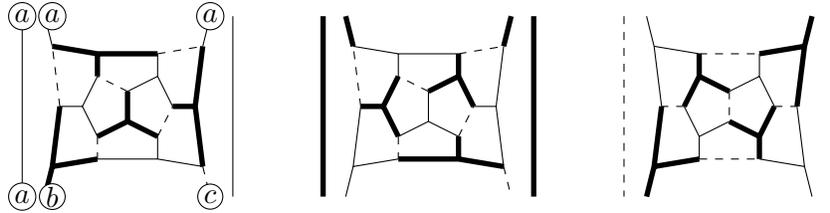
\begin{figure}[htp]
\centering
\begin{tikzpicture}[>=latex]

\foreach \x in {-1,1}
\foreach \y in {-1,1}
\foreach \a in {0,1,2}
{
\begin{scope}[xshift=4*\a cm]

\coordinate  (A\y X\a) at (0,0.2*\y);
\coordinate  (B\x\y X\a) at (0.4*\x,0.4*\y);
\coordinate  (C\x\y X\a) at (0.4*\x,0.7*\y);
\coordinate  (D\x X\a) at (0.6*\x,0);
\coordinate  (E\x X\a) at (0.9*\x,0);
\coordinate  (F\x\y X\a) at (1*\x,0.8*\y);
\coordinate  (G\x\y X\a) at (1.1*\x,1.2*\y);
\coordinate  (H\x\y X\a) at (1.4*\x,1.2*\y);
	
\end{scope}
}

\draw	
	(C11X0) -- (B11X0) -- (D1X0) 
	(A1X0) -- (B11X0)
	(C-1-1X0) -- (C1-1X0) -- (F1-1X0) 
	(B1-1X0) -- (C1-1X0)  
	(B-11X0) -- (D-1X0) -- (B-1-1X0) 
	(D-1X0) -- (E-1X0)
	(F11X0) -- (G11X0)
	(F-11X0) -- (G-11X0)
	(H11X0) -- (H1-1X0)
	(H-11X0) -- (H-1-1X0)
	(D1X1) -- (E1X1) 
	(F11X1) -- (E1X1) -- (F1-1X1)
	(B1-1X1) -- (A-1X1) -- (B-1-1X1)
	(A1X1) -- (A-1X1) 
	(C11X1) -- (C-11X1) -- (F-11X1) 
	(B-11X1) -- (C-11X1) 
	(E-1X1) -- (F-1-1X1) -- (G-1-1X1)
	(C-1-1X1) -- (F-1-1X1)
	(A-1X2) -- (B-1-1X2) 
	(C-1-1X2) -- (B-1-1X2) -- (D-1X2) 
	(A1X2) -- (B11X2)
	(C11X2) -- (B11X2) -- (D1X2) 
	(E1X2) -- (F1-1X2) -- (G1-1X2)
	(C1-1X2) -- (F1-1X2)
	(E-1X2) -- (F-11X2) -- (G-11X2)
	(C-11X2) -- (F-11X2)
	;

\draw[dashed]
	(C11X0) -- (F11X0) 
	(B1-1X0) -- (D1X0) 
	(F1-1X0) -- (G1-1X0)
	(A1X0) -- (B-11X0) 
	(B-1-1X0) -- (C-1-1X0) 
	(E-1X0) -- (F-11X0)
	(C11X1) -- (F11X1) 
	(B1-1X1) -- (D1X1) 
	(F1-1X1) -- (G1-1X1)
	(A1X1) -- (B-11X1) 
	(B-1-1X1) -- (C-1-1X1) 
	(E-1X1) -- (F-11X1)	
	(A1X2) -- (A-1X2) 
	(C-11X2) -- (C11X2) 
	(C-1-1X2) -- (C1-1X2) 
	(D1X2) -- (E1X2)
	(D-1X2) -- (E-1X2)
	(H11X2) -- (H1-1X2)
	(H-11X2) -- (H-1-1X2)
	;

\draw[line width=2]
	(D1X0) -- (E1X0) 
	(F11X0) -- (E1X0) -- (F1-1X0)
	(B1-1X0) -- (A-1X0) -- (B-1-1X0)
	(A1X0) -- (A-1X0) 
	(C11X0) -- (C-11X0) -- (F-11X0) 
	(B-11X0) -- (C-11X0) 
	(E-1X0) -- (F-1-1X0) -- (G-1-1X0)
	(C-1-1X0) -- (F-1-1X0)	
	(C11X1) -- (B11X1) -- (D1X1) 
	(A1X1) -- (B11X1)
	(C-1-1X1) -- (C1-1X1) -- (F1-1X1) 
	(B1-1X1) -- (C1-1X1)  
	(B-11X1) -- (D-1X1) -- (B-1-1X1) 
	(D-1X1) -- (E-1X1)
	(F11X1) -- (G11X1)
	(F-11X1) -- (G-11X1)
	(H11X1) -- (H1-1X1)
	(H-11X1) -- (H-1-1X1)
	(A1X2) -- (B-11X2) 
	(C-11X2) -- (B-11X2) -- (D-1X2) 
	(A-1X2) -- (B1-1X2)
	(C1-1X2) -- (B1-1X2) -- (D1X2) 
	(E1X2) -- (F11X2) -- (G11X2)
	(C11X2) -- (F11X2)
	(E-1X2) -- (F-1-1X2) -- (G-1-1X2)
	(C-1-1X2) -- (F-1-1X2)
	;
	
\node[fill=white,inner sep=1,draw,shape=circle] at (H-11X0) {\small $a$};	
\node[xshift=0.1cm,fill=white,inner sep=1,draw,shape=circle] at (G-11X0) {\small $a$};
\node[fill=white,inner sep=1,draw,shape=circle] at (G11X0) {\small $a$};

\node[fill=white,inner sep=1,draw,shape=circle] at (H-1-1X0) {\small $a$};	
\node[xshift=0.1cm,fill=white,inner sep=0.3,draw,shape=circle] at (G-1-1X0) {\small $b$};
\node[fill=white,inner sep=1,draw,shape=circle] at (G1-1X0) {\small $c$};

\end{tikzpicture}
\caption{Timezone tiling, distance $1$, edge combination $a^2b^2c$.}
\label{d1_case1}
\end{figure}

By comparing the edge lengths of the meridians, we see that each timezone tiling can only repeat itself. Moreover, applying the symmetries may produce new timezone tilings, and exchanging $a$ and $b$ simply exchanges the first and the second tilings and leave the third tiling invariant. After all the considerations, we get two families of edge congruent earth map tilings of distance $1$ for the edge combination $a^2b^2c$.
\begin{enumerate}
\item Circular product of $\binom{aaa}{abc}$, $\binom{abc}{aaa}$, $\binom{acb}{aaa}$, $\binom{aaa}{acb}$.
\item Circular product of $\binom{cab}{cba}$, $\binom{cba}{cab}$.
\end{enumerate}
In the first family, $\binom{aaa}{abc}$ is the first timezone tiling, $\binom{abc}{aaa}$ is its vertical flipping, $\binom{acb}{aaa}$ is its horizontal flipping, $\binom{aaa}{acb}$ is its $180^{\circ}$ rotation. In the second family, $\binom{cab}{cba}$ is the third timezone tiling, and $\binom{cba}{cab}$ is its vertical flipping. The choice of notation respects the symmetry of vertical flipping.

\section{Earth Map Tiling and Edge Length Combination $a^3b^2$}
\label{section_3a2b}

The edge combination $a^3b^2$ is more complicated than $a^2b^2c$ due to three possible types of neighborhood tilings. This can be dealt with by using the propagations in Figure \ref{pnd3a2b}. Another difference is that we no longer have the symmetry of exchanging $a$ and $b$. Examples can be obtained by taking $a=c$ in the classification for $a^2b^2c$. In fact, it turns out that all the timezone tilings are obtained by taking $a=c$ in the timezone tilings for $a^2b^2c$. Due to new ways of glueing timezone tilings together, however, there are edge congruent tilings for $a^3b^2$ that are not degenerated from edge congruent tilings for $a^2b^2c$.

\begin{theorem}\label{theorem3a2b}
The edge congruent earth map tilings for the edge length combination $a^3b^2$ can be classified into the following numbers of families.
\begin{itemize}
\item Distance $5$: $2$ families.
\item Distance $4$: $3$ families.
\item Distance $3$: $4$ families.
\item Distance $2$: $3$ families.
\item Distance $1$: $2$ families.
\end{itemize}
\end{theorem}

The classification is up to the symmetries of earth map tilings given in Section \ref{section_3abc}, but without exchanging $a$ and $b$. Unlike the edge combination $a^2b^2c$, the distance between the poles and the edge length combinations at the poles do not uniquely characterize the families.

\subsection*{Distance $5$}

\begin{figure}[htp]
\centering
\begin{tikzpicture}[>=latex]

\foreach \x in {0,...,9}
{

\coordinate  (A\x) at (-0.3+1.2*\x,1.2);
\coordinate  (B\x) at (-0.3+1.2*\x,0.7);
\coordinate  (C\x) at (-0.6+1.2*\x,0.3);
\coordinate  (D\x) at (-0.6+1.2*\x,-0.3);
\coordinate  (E\x) at (-0.9+1.2*\x,-0.7);
\coordinate  (F\x) at (-0.9+1.2*\x,-1.2);
\coordinate  (G\x) at (1.2*\x,0.3);
\coordinate  (H\x) at (1.2*\x,-0.3);

\coordinate  (P\x) at (-0.3+1.2*\x,0.1);
\coordinate  (Q\x) at (0.3+1.2*\x,-0.1);
	
}

\node at (P0) {\small III};
\node at (Q0) {\small III};
\node at (P1) {\small III};
\node at (Q1) {\small II};
\node at (P2) {\small I};
\node at (Q2) {\small I};
\node at (P3) {\small II};
\node at (Q3) {\small III};
\node at (P4) {\small III};
\node at (Q4) {\small III};

\draw
	(E0) -- (F0)
	(E1) -- (F1)
	(E2) -- (F2)
	(A3) -- (B3)
	(A4) -- (B4)
	(A5) -- (B5)
	(B0) -- (C0) -- (D0) -- (E0)
	(B1) -- (C1) -- (D1) -- (E1)
	(B2) -- (C2) -- (D2) -- (E2)
	(B3) -- (C3) -- (D3) -- (E3)
	(B4) -- (C4) -- (D4) -- (E4)
	(B5) -- (C5) -- (D5) -- (E5)
	(B0) -- (G0) -- (C1)
	(B1) -- (G1) -- (C2)
	(D2) -- (H2)
	(G2) -- (C3)
	(D3) -- (H3) -- (E4)
	(D4) -- (H4) -- (E5)
	;

\draw[line width=2]
	(A0) -- (B0)
	(A1) -- (B1)
	(A2) -- (B2)
	(E3) -- (F3)
	(E4) -- (F4)
	(E5) -- (F5)
	(G0) -- (H0) -- (E1)
	(D0) -- (H0)
	(G1) -- (H1) -- (E2)
	(D1) -- (H1)
	(B2) -- (G2) -- (H2) -- (E3)
	(B3) -- (G3) -- (H3)
	(G3) -- (C4)
	(B4) -- (G4) -- (H4)
	(G4) -- (C5)
	;

\draw (Q2) circle (0.2);

\node at (2.7,0.8) {\small $\times$};
\node at (3.3,-0.8) {\small $\times$};

\end{tikzpicture}
\caption{No neighborhood of type I or II.}
\label{d5_case2A}
\end{figure}
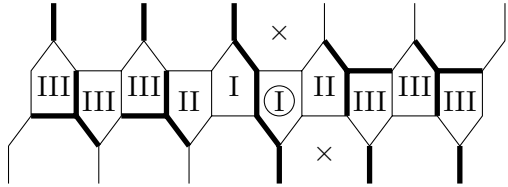

Consider a core tile in the earth map tiling. By the propagation in Figure \ref{pnd3a2b}, if the neighborhood of this tile is of type I, as indicated by the circle in Figure \ref{d5_case2A}, then the two nearby tiles labeled $\times$ must have poles as vertices. Up to symmetry, we may assume that these two tiles are located as in Figure \ref{d5_case2A}. Then we get all the edges in the neighborhood. Moreover, the propagation tells us the types of the two nearby neighborhoods. This determines all the edges in the two nearby neighborhoods. By using the propagation again, we can successively determine all the edges. By the time we get type III neighborhoods, we find that the pattern repeats. Moreover, the edges at the north poles should have length $b$ to the left of the initial tile, and should have length $a$ to the right of the initial tile. However, these two directions will eventually meet, causing a contradiction. We conclude that type I neighborhood tiling is actually forbidden. 

If the neighborhood of a core tile is of type II, then the propagation in Figure \ref{pnd3a2b} similarly implies that a nearby neighborhood is of type I. Since we have argued that there is no type I neighborhood, there is also no type II neighborhood.

\begin{figure}[htp]
\centering
\begin{tikzpicture}[>=latex]

\foreach \x in {0,...,9}
{

\coordinate  (A\x) at (-0.3+1.2*\x,1.2);
\coordinate  (B\x) at (-0.3+1.2*\x,0.7);
\coordinate  (C\x) at (-0.6+1.2*\x,0.3);
\coordinate  (D\x) at (-0.6+1.2*\x,-0.3);
\coordinate  (E\x) at (-0.9+1.2*\x,-0.7);
\coordinate  (F\x) at (-0.9+1.2*\x,-1.2);
\coordinate  (G\x) at (1.2*\x,0.3);
\coordinate  (H\x) at (1.2*\x,-0.3);

\coordinate  (P\x) at (-0.3+1.2*\x,0.1);
\coordinate  (Q\x) at (0.3+1.2*\x,-0.1);
	
}

\node at (P0) {\small $P_1$};
\node at (Q0) {\small $P_2$};
\node[gray!70] at (P1) {\small $P_1$};
\node[gray!70] at (Q1) {\small $P_2$};
\node at (P4) {\small $P_3$};
\node at (Q4) {\small $P_4$};
\node at (P5) {\small $P_5$};
\node at (Q5) {\small $P_5$};
\node at (P6) {\small $P_4$};
\node at (Q6) {\small $P_3$};
\node[gray!70] at (P7) {\small $P_3$};
\node[gray!70] at (Q7) {\small $P_4$};
\node[gray!70] at (P8) {\small $P_5$};

\draw
	(A0) -- (B0) 
	(C0) -- (D0) -- (E0)
	(D0) -- (H0)
	(A1) -- (B1)
	(C1) -- (D1) -- (E1)
	(G4) -- (H4) 
	(D4) -- (H4) -- (E5)
	(A5) -- (B5) -- (C5)
	(B5) -- (G5)
	(H5) -- (E6) -- (F6)
	(D6) -- (E6)
	(B6) -- (G6) -- (H6)
	(G6) -- (C7)
	;

\draw[line width=2]
	(E0) -- (F0) 
	(B0) -- (G0) -- (H0)
	(G0) -- (C1)
	(E1) -- (F1)
	(A4) -- (B4) -- (C4) 
	(B4) -- (G4)
	(D4) -- (E4) -- (F4)
	(C5) -- (D5) -- (E5)
	(D5) -- (H5)
	(B6) -- (C6) -- (D6)
	(G5) -- (C6)
	(H6) -- (E7) -- (F7)
	(D7) -- (E7)
	(A7) -- (B7) -- (C7)
	;

\draw
	(B0) -- (C0) 
	(H0) -- (E1)
	(B1) -- (C1)
	(C4) -- (D4) 
	(E5) -- (F5)
	(G4) -- (C5)
	(G5) -- (H5)
	(A6) -- (B6)
	(D6) -- (H6)
	(C7) -- (D7)
	;

\draw[gray!70]
	(D1) -- (H1) 
	(A2) -- (B2)
	(C2) -- (D2) -- (E2)
	(D7) -- (H7) -- (E8) 
	(G7) -- (H7)
	(A8) -- (B8) -- (C8)
	(B8) -- (G8)
	(H8) -- (E9) -- (F9)
	;

\draw[gray!70,line width=2]
	(B1) -- (G1) -- (H1) 
	(G1) -- (C2)
	(E2) -- (F2)
	(B7) -- (G7) 
	(C8) -- (D8) -- (E8)
	(D8) -- (H8)
	;

\draw[gray!70]
	(H1) -- (E2) 
	(B2) -- (C2)
	(G7) -- (C8) 
	(E8) -- (F8)
	(G8) -- (H8)
	;

\end{tikzpicture}
\caption{Earth map tiling, distance $5$, edge length combination $a^3b^2$.}
\label{d5_case2}
\end{figure}
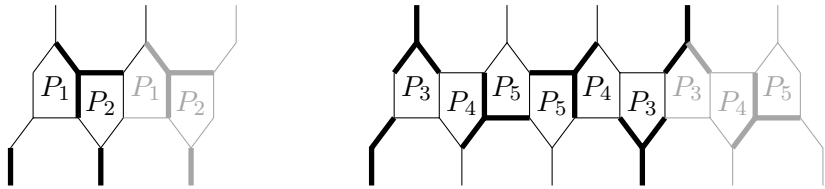

So the core tiles can only have type III neighborhood. Observe that the edges in a type III neighborhood are completely determined by the edges of the center tile and the two ``legs'' circled in Figure \ref{pnd3a2b}. Guided by this observation, we find that, up to symmetry, the neighborhoods of the tiles $P_1,P_2,P_3,P_4,P_5$ in Figure \ref{d5_case2} give all the possible ways a type III neighborhood fits into the tiling. Once we fix one such neighborhood, the edges of the core tiles $P_-$ and $P_+$ on the left and right of $P$ are determined. Since the neighborhoods of $P_-$ and $P_+$ must be of type III, all the edges in the neighborhoods of $P_-$ and $P_+$ can be further determined. The process continues, so that the edges of the initial tile uniquely determine all the other edges.

We get two repeating patterns, which give two families of edge congruent earth map tilings of distance $5$ for the edge combination $a^3b^2$.
\begin{enumerate}
\item Circular product of $\binom{a}{b}$.
\item Circular product of $\binom{baa}{baa}$.
\end{enumerate}
The notation follows the same choice as the case of $a^2b^2c$. These are actually obtained by taking $c=a$ in the two families for the combination $a^2b^2c$.

\subsection*{Distance $4$}

In the earth map tiling of distance $4$, a core tile has exactly two nearby tiles with vertices of degree $>3$, and these two tiles are always next to each other. Among the three types of neighborhood tilings in Figure \ref{pnd3a2b}, only type III has this property. Therefore the other two types are forbidden. 

\begin{figure}[htp]
\centering
\begin{tikzpicture}[>=latex]

\foreach \x in {-1,1}
\foreach \y in {-1,1}
\foreach \a in {0,1,2}
{
\begin{scope}[xshift=4.8*\a cm]

\coordinate  (A\y X\a) at (0,1.2*\y);
\coordinate  (B\y X\a) at (0,0.8*\y);
\coordinate  (C\x X\a) at (0.3*\x,0);
\coordinate  (D\x\y X\a) at (0.4*\x,0.5*\y);
\coordinate  (E\x\y X\a) at (0.9*\x,0.4*\y);
\coordinate  (F\x\y X\a) at (1.2*\x,0.7*\y);
\coordinate  (G\x\y X\a) at (1.2*\x,1.2*\y);
\coordinate  (H\x X\a) at (1.4*\x,0);
\coordinate  (I\x X\a) at (1.7*\x,0);
\coordinate  (J\x\y X\a) at (1.9*\x,0.7*\y);
\coordinate  (K\x\y X\a) at (1.9*\x,1.2*\y);

\coordinate  (P\a) at (0,-1.5);
	
\end{scope}
}

\foreach \x in {-1,1}
\foreach \y in {-1,1}
\foreach \a in {0,1,2}
\draw
	(A\y X\a) -- (B\y X\a) -- (D\x\y X\a) -- (C\x X\a) 
	(D\x\y X\a) -- (E\x\y X\a) -- (F\x\y X\a)
	(G\x\y X\a) -- (F\x\y X\a) -- (H\x X\a) -- (I\x X\a) -- (J\x\y X\a) -- (K\x\y X\a) 
	(E\x 1X\a) -- (E\x -1X\a);

\foreach \a in {0,1,2}
\draw
	(C1X\a) -- (C-1X\a);
	
\draw[line width=2]
	(D11X0) -- (E11X0) -- (F11X0)  
	(E11X0) -- (E1-1X0)
	(B1X0) -- (D-11X0) -- (E-11X0)
	(D1-1X0) -- (B-1X0) -- (D-1-1X0)
	(B-1X0) -- (A-1X0)
	(F-11X0) -- (H-1X0) -- (F-1-1X0)
	(H-1X0) -- (I-1X0)
	(A1X1) -- (B1X1) 
	(D11X1) -- (C1X1) -- (D1-1X1)
	(C1X1) -- (C-1X1)
	(G11X1) -- (F11X1)
	(H1X1) -- (F1-1X1) -- (G1-1X1)
	(E1-1X1) -- (F1-1X1)
	(E-11X1) -- (E-1-1X1) -- (D-1-1X1)
	(E-1-1X1) -- (F-1-1X1)
	(K11X1) -- (J11X1)
	(G-11X1) -- (F-11X1)
	(K-11X1) -- (J-11X1)
	(I-1X1) -- (J-1-1X1) -- (K-1-1X1)
	(B1X2) -- (D11X2) -- (E11X2) 
	(D11X2) -- (C1X2)
	(H1X2) -- (F1-1X2) -- (G1-1X2)
	(I1X2) -- (J11X2) -- (K11X2)
	(B-1X2) -- (D-1-1X2) -- (E-1-1X2)
	(D-1-1X2) -- (C-1X2)
	(H-1X2) -- (F-11X2) -- (G-11X2)
	(I-1X2) -- (J-1-1X2) -- (K-1-1X2)
	;

\draw[gray!70]
	(K11X0) -- (J11X0) 
	(J1-1X0) -- (K1-1X0);
\draw[gray!70,line width=2]
	(H1X0) -- (I1X0)
	(J11X0) -- (I1X0) -- (J1-1X0)
	;

\end{tikzpicture}
\caption{Timezone tiling, distance $4$, edge length combination $a^3b^2$.}
\label{d4_case2}
\end{figure}
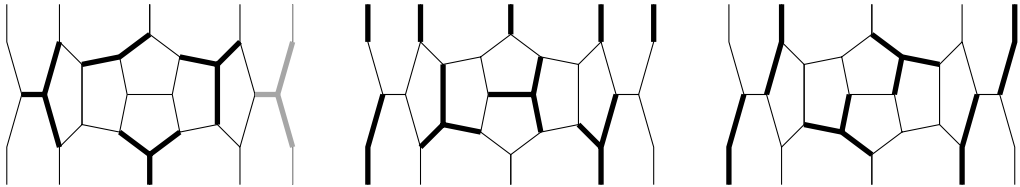

Up to symmetry, the core tiles in Figure \ref{d4_case2} give all the possible ways that a type III neighborhood fits into the tiling. Then it is easy to see that each way determines all the other edges in the timezone. The gray part of the first timezone tiling is not yet uniquely determined (but will be).

By comparing the edge length arrangements of the meridians, we see that the second and third timezone tilings can only repeat themselves. Therefore the first timezone tiling cannot be glued to the second or the third. This actually determines the gray part of the first tiling, and we conclude that the first timezone tiling can repeat either itself or its vertical flipping. Note also that the three tilings are actually obtained by taking $c=a$ in Figure \ref{d4_case1}.

So we get three families of edge congruent earth map tilings of distance $4$ for the edge combination $a^3b^2$. 
\begin{enumerate}
\item Circular product of $\binom{aaa}{aba}$, $\binom{aba}{aaa}$.
\item Circular product of $\binom{bbb}{aab}$.
\item Circular product of $\binom{baa}{aab}$.
\end{enumerate}
The notation follows the same choice as the case of $a^2b^2c$. Moreover, $\binom{aba}{aaa}$ is the vertical flipping of $\binom{aaa}{aba}$.

\subsection*{Distance $3$}

\begin{figure}[htp]
\centering
\begin{tikzpicture}[>=latex]

\foreach \a in {0,1,2}
\foreach \b in {0,1}
{
\begin{scope}[shift={(4.8*\a cm, -3*\b cm)}]

\foreach \x in {-1,1}
\foreach \y in {-1,1}
{
\coordinate  (A\x\y X\a\b) at (0.32*\x-0.16*\y,0.16*\y);
\coordinate  (B\x\y X\a\b) at (0.32*\x-0.16*\y,0.6*\y);
\coordinate  (C\x\y X\a\b) at (0.8*\x-0.16*\y,1.2*\y);
\coordinate  (D\x\y X\a\b) at (0.8*\x-0.16*\y,0.7*\y);
\coordinate  (E\x\y X\a\b) at (0.96*\x-0.16*\y,0.16*\y);
\coordinate  (F\x\y X\a\b) at (1.6*\x-0.16*\y,0.16*\y);
\coordinate  (G\x\y X\a\b) at (1.6*\x-0.16*\y,1.2*\y);

\draw
	(A\x\y X\a\b) -- (B\x\y X\a\b)
	(C\x\y X\a\b) -- (D\x\y X\a\b) -- (B\x\y X\a\b)
	(D\x\y X\a\b) -- (E\x\y X\a\b)
	(F\x\y X\a\b) -- (G\x\y X\a\b);
}

\draw
	(F1-1X\a\b) -- (F11X\a\b) -- (E1-1X\a\b) -- (E11X\a\b) -- (A1-1X\a\b) -- (A11X\a\b) -- (A-1-1X\a\b) -- (A-11X\a\b) -- (E-1-1X\a\b) -- (E-11X\a\b) -- (F-1-1X\a\b) -- (F-11X\a\b)
	(B11X\a\b) -- (B-11X\a\b)
	(B1-1X\a\b) -- (B-1-1X\a\b);
	
\end{scope}
}

\draw[line width=2]
	(B-11X00) -- (B11X00) -- (D11X00) 
	(B1-1X00) -- (B-1-1X00) -- (D-1-1X00)
	(A11X00) -- (B11X00) 
	(A-1-1X00) -- (B-1-1X00)
	(E11X00) -- (E1-1X00) -- (F11X00) 
	(E-1-1X00) -- (E-11X00) -- (F-1-1X00) 
	(D1-1X00) -- (E1-1X00)
	(D-11X00) -- (E-11X00)
	(C11X10) -- (D11X10) -- (E11X10) 
	(B11X10) -- (D11X10)
	(C1-1X10) -- (D1-1X10) -- (E1-1X10)
	(B1-1X10) -- (D1-1X10)
	(A-11X10) -- (A-1-1X10) -- (A11X10)
	(A-1-1X10) -- (B-1-1X10)
	(E-1-1X10) -- (E-11X10) -- (F-1-1X10)
	(D-11X10) -- (E-11X10)
	(A11X20) -- (A1-1X20) -- (E11X20) 
	(A-1-1X20) -- (A-11X20) -- (E-1-1X20)
	(A1-1X20) -- (B1-1X20)
	(A-11X20) -- (B-11X20)
	(C11X20) -- (D11X20)
	(C1-1X20) -- (D1-1X20)
	(C-11X20) -- (D-11X20)
	(C-1-1X20) -- (D-1-1X20)
	(G11X20) -- (F11X20)
	(G1-1X20) -- (F1-1X20)
	(G-11X20) -- (F-11X20)
	(G-1-1X20) -- (F-1-1X20)
	(A1-1X01) -- (E11X01) -- (E1-1X01) 
	(E11X01) -- (D11X01)
	(F11X01) -- (F1-1X01) -- (G1-1X01)
	(B11X01) -- (B-11X01) -- (D-11X01)
	(B-11X01) -- (A-11X01)
	(B1-1X01) -- (B-1-1X01) -- (D-1-1X01)
	(B-1-1X01) -- (A-1-1X01)
	(F-11X01) -- (F-1-1X01) -- (G-1-1X01)
	(F-1-1X01) -- (E-11X01)
	(C11X11) -- (D11X11) -- (E11X11) 
	(B11X11) -- (D11X11)
	(A-1-1X11) -- (A-11X11) -- (E-1-1X11)
	(A-11X11) -- (B-11X11)
	(B-1-1X11) -- (B1-1X11) -- (D1-1X11)
	(B1-1X11) -- (A1-1X11)
	(F11X11) -- (F1-1X11) -- (G1-1X11)
	(F-11X11) -- (F-1-1X11) -- (G-1-1X11)
	(F-1-1X11) -- (E-11X11)
	(A1-1X21) -- (E11X21) -- (E1-1X21) 
	(E11X21) -- (D11X21)
	(A-11X21) -- (A-1-1X21) -- (A11X21)
	(A-1-1X21) -- (B-1-1X21)
	(C-11X21) -- (D-11X21) -- (E-11X21)
	(B-11X21) -- (D-11X21)
	(C1-1X21) -- (D1-1X21)
	(C-1-1X21) -- (D-1-1X21)
	(G1-1X21) -- (F1-1X21)
	(G-1-1X21) -- (F-1-1X21)
	;

\end{tikzpicture}
\caption{Timezone tiling, distance $3$, edge combination $a^3b^2$.}
\label{d3_case2}
\end{figure}
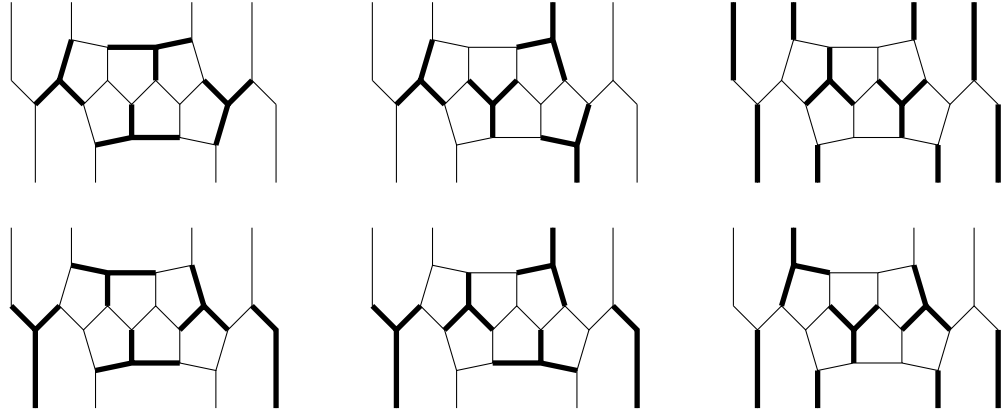

Similar to the distance $4$ case, the location of the poles implies that the neighborhood of a core tile must be of type III. Up to symmetry, the core tiles in Figure \ref{d3_case2} give all the possible ways that a type III neighborhood fits into the tiling. Then each way determines all the other edges in the timezone. The six timezone tilings are actually obtained by taking $c=a$ in Figure \ref{d3_case1}, but arranged in different order.

The first and the second tilings can be glued together. The fourth and the fifth tilings can also be glued together. The third and the sixth tilings can only repeat themselves. Moreover, only the $180^{\circ}$ rotation of the middle tiling in the first row gives a new timezone tiling. So we get four families of edge congruent earth map tilings of distance $3$ for the edge combination $a^3b^2$. 
\begin{enumerate}
\item Circular product of $\binom{aaa}{aaa}$, $\binom{aab}{aba}$, $\binom{aba}{baa}$.
\item Circular product of $\binom{bbb}{bbb}$.
\item Circular product of $\binom{aaa}{aab}$, $\binom{aab}{aab}$.
\item Circular product of $\binom{aba}{bbb}$.
\end{enumerate}
In the first family, $\binom{aba}{baa}$ is the $180^{\circ}$ rotation of $\binom{aab}{aba}$.

\subsection*{Distance $2$}

\begin{figure}[htp]
\centering
\begin{tikzpicture}[>=latex]

\foreach \a in {0,1,2}
\foreach \b in {0,1}
{
\begin{scope}[shift={(4*\a cm, -3*\b cm)}]

\foreach \y in {-1,1}
{
\coordinate  (A\y X\a\b) at (0,0.2*\y);
\coordinate  (B\y X\a\b) at (-0.4,0.4*\y);
\coordinate  (C\y X\a\b) at (0.4,0.4*\y);
\coordinate  (D\y X\a\b) at (-0.1,0.8*\y);
\coordinate  (E\y X\a\b) at (0.4,0.7*\y);
\coordinate  (F\y X\a\b) at (-0.1,1.2*\y);
\coordinate  (GX\a\b) at (0.6,0);
\coordinate  (HX\a\b) at (0.9,0);
\coordinate  (I\y X\a\b) at (1,0.5*\y);
\coordinate  (J\y X\a\b) at (1.3,0.6*\y);
\coordinate  (K\y X\a\b) at (1.4,1.2*\y);
\coordinate  (L\y X\a\b) at (1.7,1.2*\y);
\coordinate  (MX\a\b) at (1.8,0);
\coordinate  (NX\a\b) at (-0.6,0);
\coordinate  (OX\a\b) at (-0.9,0);
\coordinate  (P\y X\a\b) at (-1,1.2*\y);

\draw
	(A\y X\a\b) -- (B\y X\a\b) -- (D\y X\a\b) -- (E\y X\a\b) -- (C\y X\a\b) -- cycle
	(D\y X\a\b) -- (F\y X\a\b)
	(C\y X\a\b) -- (GX\a\b) -- (HX\a\b) -- (I\y X\a\b)
	(E\y X\a\b) -- (I\y X\a\b) -- (J\y X\a\b) -- (K\y X\a\b)
	(L\y X\a\b) -- (MX\a\b)
	(B\y X\a\b) -- (NX\a\b) -- (OX\a\b) -- (P\y X\a\b);
}

\end{scope}

\draw
	(A1X\a\b) -- (A-1X\a\b)
	(J1X\a\b) -- (J-1X\a\b);
}

\draw[line width=2]
	(B1X00) -- (NX00) -- (B-1X00) 
	(OX00) -- (NX00)
	(A-1X00) -- (C-1X00) -- (GX00) 
	(C-1X00) -- (E-1X00) 
	(D1X00) -- (E1X00) -- (I1X00) 
	(C1X00) -- (E1X00)
	(J1X00) -- (J-1X00) -- (K-1X00)
	(I-1X00) -- (J-1X00)
	(A1X10) -- (B1X10) -- (D1X10) 
	(B1X10) -- (NX10)
	(C1X10) -- (GX10) -- (C-1X10) 
	(GX10) -- (HX10) 
	(B-1X10) -- (D-1X10) -- (E-1X10)
	(D-1X10) -- (F-1X10)
	(I1X10) -- (J1X10) -- (K1X10) 
	(J1X10) -- (J-1X10)
	(D1X20) -- (E1X20) -- (I1X20) 
	(C1X20) -- (E1X20)
	(B-1X20) -- (A-1X20) -- (C-1X20) 
	(A1X20) -- (A-1X20) 
	(E-1X20) -- (I-1X20) -- (J-1X20) 
	(HX20) -- (I-1X20)
	(L1X20) -- (MX20) -- (L-1X20)
	(P1X20) -- (OX20) -- (P-1X20)
	(OX20) -- (NX20)
	(A1X01) -- (B1X01) -- (D1X01) 
	(B1X01) -- (NX01)
	(A-1X01) -- (C-1X01) -- (GX01) 
	(C-1X01) -- (E-1X01) 
	(E1X01) -- (I1X01) -- (J1X01) 
	(HX01) -- (I1X01)
	(OX01) -- (P-1X01)
	(D-1X01) -- (F-1X01) 
	(J-1X01) -- (K-1X01)
	(MX01) -- (L-1X01)
	(B1X11) -- (D1X11) -- (E1X11) 
	(D1X11) -- (F1X11)
	(B-1X11) -- (A-1X11) -- (C-1X11) 
	(A1X11) -- (A-1X11) 
	(I1X11) -- (HX11) -- (I-1X11) 
	(GX11) -- (HX11)
	(OX11) -- (P-1X11)
	(D-1X11) -- (F-1X11) 
	(J-1X11) -- (K-1X11)
	(MX11) -- (L-1X11)
	;

\fill[white]
	(6.9,-4.25) rectangle (9.9,-1.75);

\end{tikzpicture}
\caption{Timezone tiling, distance $2$, edge combination $a^3b^2$.}
\label{d2_case2}
\end{figure}
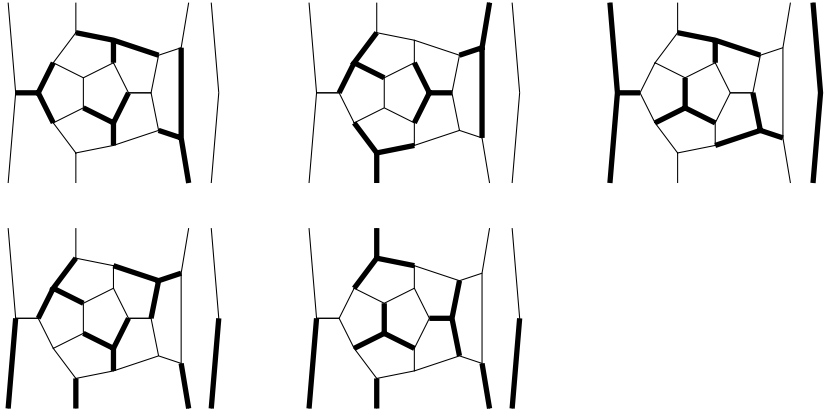

The same argument as the distances $3$ and $4$ cases give all the timezone tilings in Figure \ref{d2_case2}, which happens to be again obtained by taking $c=a$ in Figure \ref{d2_case2}, but arranged in different order. We have three families of edge congruent earth map tilings of distance $2$ for the edge combination $a^3b^2$.
\begin{enumerate}
\item Circular product of $\binom{aaa}{aab}$, $\binom{aab}{aba}$, $\binom{aab}{aaa}$, $\binom{aba}{aab}$.
\item Circular product of $\binom{aaa}{bbb}$, $\binom{aba}{bbb}$.
\item Circular product of $\binom{baa}{baa}$, $\binom{baa}{baa}^{\text{v.flip}}$.
\end{enumerate}
The first family is obtained by glueing the first and second tilings $\binom{aaa}{aab}$, $\binom{aab}{aba}$ and their vertical flippings $\binom{aab}{aaa}$, $\binom{aba}{aab}$. The second family is obtained by glueing the fourth and fifth tilings $\binom{aaa}{bbb}$, $\binom{aba}{bbb}$. The third family is obtained by glueing the third tiling $\binom{baa}{baa}$ and its vertical flipping $\binom{baa}{baa}^{\text{v.flip}}$. 

\subsection*{Distance $1$}

\begin{figure}[htp]
\centering
\begin{tikzpicture}[>=latex]

\foreach \x in {-1,1}
\foreach \y in {-1,1}
\foreach \a in {0,1,2}
{
\begin{scope}[xshift=4*\a cm]

\coordinate  (A\y X\a) at (0,0.2*\y);
\coordinate  (B\x\y X\a) at (0.4*\x,0.4*\y);
\coordinate  (C\x\y X\a) at (0.4*\x,0.7*\y);
\coordinate  (D\x X\a) at (0.6*\x,0);
\coordinate  (E\x X\a) at (0.9*\x,0);
\coordinate  (F\x\y X\a) at (1*\x,0.8*\y);
\coordinate  (G\x\y X\a) at (1.1*\x,1.2*\y);
\coordinate  (H\x\y X\a) at (1.4*\x,1.2*\y);

\draw
	(A\y X\a) -- (B\x\y X\a) -- (D\x X\a) -- (E\x X\a)  -- (F\x\y X\a) 
	(B\x\y X\a) -- (C\x\y X\a) -- (F\x\y X\a)  -- (G\x\y X\a);
	
\end{scope}
}

\foreach \x in {-1,1}
\foreach \a in {0,1,2}
\draw
	(A1X\a) -- (A-1X\a)
	(C1\x X\a) -- (C-1\x X\a)
	(H\x 1X\a) -- (H\x -1X\a);

\draw[line width=2]
	(D1X0) -- (E1X0) 
	(F11X0) -- (E1X0) -- (F1-1X0)
	(B1-1X0) -- (A-1X0) -- (B-1-1X0)
	(A1X0) -- (A-1X0) 
	(C11X0) -- (C-11X0) -- (F-11X0) 
	(B-11X0) -- (C-11X0) 
	(E-1X0) -- (F-1-1X0) -- (G-1-1X0)
	(C-1-1X0) -- (F-1-1X0)
	(A1X1) -- (B-11X1) 
	(C-11X1) -- (B-11X1) -- (D-1X1) 
	(A-1X1) -- (B1-1X1)
	(C1-1X1) -- (B1-1X1) -- (D1X1) 
	(E1X1) -- (F11X1) -- (G11X1)
	(C11X1) -- (F11X1)
	(E-1X1) -- (F-1-1X1) -- (G-1-1X1)
	(C-1-1X1) -- (F-1-1X1)
	(A-1X2) -- (B1-1X2) 
	(C1-1X2) -- (B1-1X2) -- (D1X2) 
	(C-11X2) -- (C11X2) -- (F11X2) 
	(B11X2) -- (C11X2) 
	(B-11X2) -- (D-1X2) -- (B-1-1X2) 
	(D-1X2) -- (E-1X2) 
	(F1-1X2) -- (G1-1X2)
	(F-1-1X2) -- (G-1-1X2)
	(H11X2) -- (H1-1X2)
	(H-11X2) -- (H-1-1X2)
	;

\end{tikzpicture}
\caption{Timezone tiling, distance $1$, edge combination $a^3b^2$.}
\label{d1_case2}
\end{figure}
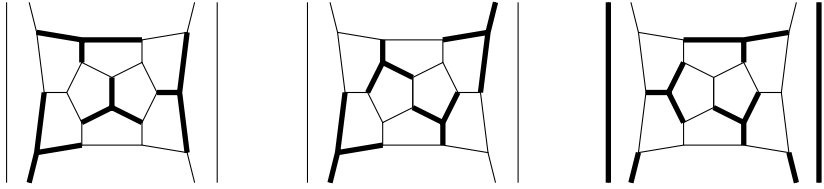

The argument for the case of distance 1 is the same as before. Figure \ref{d1_case2} gives all the timezone tilings, which happens to be obtained by taking $c=a$ in Figure \ref{d1_case1}. We have two families of edge congruent earth map tilings of distance $1$ for the edge combination $a^3b^2$.
\begin{enumerate}
\item Circular product of $\binom{aaa}{aba}$, $\binom{aba}{aaa}$, $\binom{aaa}{aab}$, $\binom{aab}{aaa}$, $\binom{aab}{aba}$, $\binom{aba}{aab}$.
\item Circular product of $\binom{baa}{bbb}$, $\binom{bbb}{baa}$, $\binom{baa}{bbb}^{\text{h.flip}}$, $\binom{bbb}{baa}^{\text{h.flip}}$.
\end{enumerate}
The first family is obtained by glueing the first tiling $\binom{aaa}{aba}$ and its vertical flipping $\binom{aba}{aaa}$, its horizontal flipping $\binom{aaa}{aab}$, and its $180^{\circ}$ rotation $\binom{aab}{aaa}$, together with the second tiling $\binom{aab}{aba}$ and its vertical flipping $\binom{aba}{aab}$. The second family is obtained by glueing the third tiling $\binom{baa}{bbb}$ and its vertical flipping $\binom{bbb}{baa}$, its horizontal flipping $\binom{baa}{bbb}^{\text{h.flip}}$, and its $180^{\circ}$ rotation $\binom{bbb}{baa}^{\text{h.flip}}$ (which is the composition of horizontal and vertical flippings). 

\section{Earth Map Tiling and Edge Length Combination $a^4b$}
\label{section_4ab}

The complication for the edge combination $a^4b$ lies in the big number of edge congruent tilings of the timezone. Such complication already appeared for the edge combination $a^3b^2$, especially when the distance between the two poles is small. Due to the complication, the families of tilings will be described by cycles in certain directed graphs. A node $m$ of the graph is an edge combination of the meridian, and a directed edge (arrow) from $m_1$ to $m_2$ is an edge congruent timezone tiling with $m_1$ and $m_2$ as the left and right boundaries.

We also remark that, since each $b$-edge is shared by exactly two tiles, and the number of tiles in a timezone is even ($4$ for distance $5$ and $12$ for distances $4,3,2,1$), the numbers of $b$-edges in $m_1$ and $m_2$ must have the same parity. This gives even and odd families of edge congruent tilings for each distance.

\begin{theorem}\label{theorem4ab}
The edge congruent earth map tilings for the edge length combination $a^4b$ can be classified into two families for each distance.
\end{theorem}

The proof will not make use of the neighborhood tilings in Section \ref{section_nd}. We will simply list all the timezone tilings.

\subsection*{Distance $5$}

\begin{figure}[htp]
\centering
\begin{tikzpicture}[>=latex]

\foreach \x in {0,1}
\foreach \a in {0,...,7}
\foreach \b in {0,1}
{
\begin{scope}[shift={(1.6*\a cm, -1.6*\b cm)}]

\coordinate  (A\x X\a\b) at (-0.2+0.8*\x,0.7);
\coordinate  (B\x X\a\b) at (-0.2+0.8*\x,0.35);
\coordinate  (C\x X\a\b) at (-0.4+0.8*\x,0.2);
\coordinate  (D\x X\a\b) at (-0.4+0.8*\x,-0.2);
\coordinate  (E\x X\a\b) at (-0.6+0.8*\x,-0.35);
\coordinate  (F\x X\a\b) at (-0.6+0.8*\x,-0.7);

\coordinate  (GX\a\b) at (0,0.2);
\coordinate  (HX\a\b) at (0,-0.2);
	
\end{scope}
}

\foreach \a in {0,...,7}
\foreach \b in {0,1}
{

\foreach \x in {0,1}
\draw
	(A\x X\a\b) -- (B\x X\a\b) -- (C\x X\a\b) -- (D\x X\a\b) -- (E\x X\a\b) -- (F\x X\a\b);

\draw
	(B0X\a\b) -- (GX\a\b) -- (C1X\a\b)
	(D0X\a\b) -- (HX\a\b) -- (E1X\a\b)
	(GX\a\b) -- (HX\a\b);
}
	
\draw[line width=1.8]
	(B0X00) -- (GX00)	
	(E1X00) -- (HX00)	
	(C1X10) -- (GX10)	
	(D0X10) -- (HX10)
	(E1X20) -- (F1X20)	
	(C1X20) -- (D1X20)
	(B0X20) -- (GX20)
	(E1X30) -- (F1X30)	
	(A1X30) -- (B1X30)
	(GX30) -- (HX30)
	(D0X40) -- (E0X40)	
	(B0X40) -- (C0X40)
	(D1X40) -- (E1X40)
	(B1X40) -- (C1X40)
	(D0X50) -- (E0X50)	
	(B0X50) -- (C0X50)
	(D1X50) -- (E1X50)
	(A1X50) -- (B1X50)
	(D0X60) -- (E0X60)	
	(B0X60) -- (C0X60)
	(C1X60) -- (D1X60)
	(A1X60) -- (B1X60)
	(D0X70) -- (E0X70)	
	(B0X70) -- (C0X70)
	(GX70) -- (C1X70)
	(C0X01) -- (D0X01)	
	(A1X01) -- (B1X01)
	(C1X01) -- (D1X01)
	(E1X01) -- (F1X01)
	(C1X11) -- (D1X11)	
	(A0X11) -- (B0X11)
	(C0X11) -- (D0X11)
	(E0X11) -- (F0X11)
	(D1X21) -- (E1X21)	
	(A0X21) -- (B0X21)
	(C0X21) -- (D0X21)
	(E0X21) -- (F0X21)
	(A0X31) -- (B0X31)	
	(E1X31) -- (F1X31)
	(GX31) -- (HX31)
	(D0X41) -- (E0X41)	
	(D1X41) -- (E1X41)
	(B0X41) -- (GX41)
	(D0X51) -- (E0X51)	
	(C1X51) -- (D1X51)
	(B0X51) -- (GX51)
	(D0X61) -- (E0X61)	
	(B1X61) -- (C1X61)
	(GX61) -- (HX61)
	(D0X71) -- (E0X71)	
	(A1X71) -- (B1X71)
	(GX71) -- (HX71)
	;
	
\end{tikzpicture}
\caption{Timezone tiling (partial), distance $5$, edge combination $a^4b$.}
\label{timezone4ab_d5}
\end{figure}
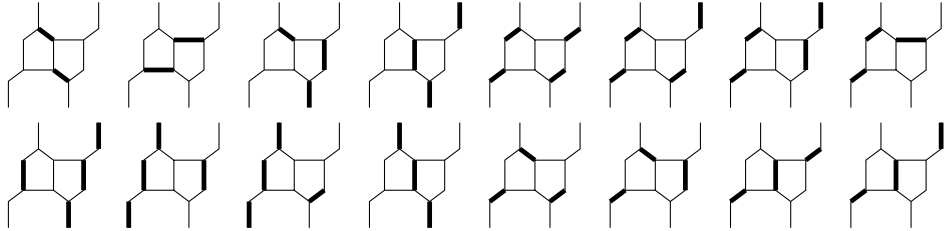

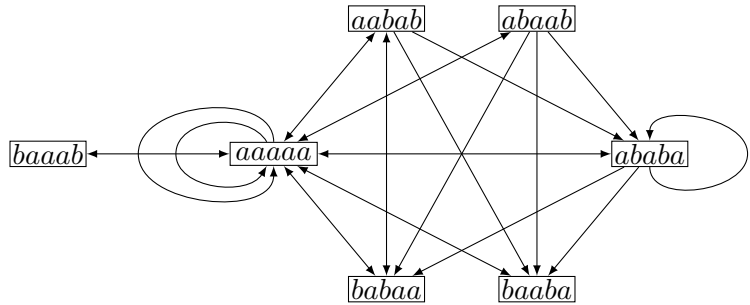
\begin{figure}[htp]
\centering
\begin{tikzpicture}[>=latex]

\node[rectangle,draw,inner sep=2] (A) at (0,0) {\small $aaaaa$};
\node[rectangle,draw,inner sep=1] (B1) at (1.5,1.8) {\small $aabab$};
\node[rectangle,draw,inner sep=1] (B-1) at (1.5,-1.8) {\small $babaa$};
\node[rectangle,draw,inner sep=1] (C1) at (3.5,1.8) {\small $abaab$};
\node[rectangle,draw,inner sep=1] (C-1) at (3.5,-1.8) {\small $baaba$};
\node[rectangle,draw,inner sep=1] (D) at (5,0) {\small $ababa$};
\node[rectangle,draw,inner sep=1] (E) at (-3,0) {\small $baaab$};	

\draw
	(A) edge[<->] (B1)
	(A) edge[<->] (B-1)
	(A) edge[<->] (C1)
	(A) edge[<->] (C-1)
	(A) edge[<->] (D)
	(A) edge[<->] (E)
	(B1) edge[<->] (B-1)
	(B1) edge[->] (C-1)
	(B1) edge[->] (D)
	(C1) edge[->] (B-1)
	(C1) edge[->] (C-1)
	(C1) edge[->] (D)
	(D) edge[->] (B-1)
	(D) edge[->] (C-1);

\draw[->]
	(A) to[out=120,in=90] (-1.3,0) to[out=270,in=240] (A);
\draw[->]
	(A) to[out=90,in=90] (-1.8,0) to[out=270,in=270] (A);
\draw[->]
	(D) to[out=270,in=270] (6.3,0) to[out=90,in=90] (D);
		
\end{tikzpicture}
\caption{Even family of earth map tiling, distance $5$, edge combination $a^4b$.}
\label{4abeven_d5}
\end{figure}

Figure \ref{timezone4ab_d5} shows some examples of timezone tilings. The tilings in the first row belong to the even family, which means that there are $0$ or $2$ $b$-edges in each boundary meridian. Instead of listing all the timezone tilings in the even family, we only represent all of them as arrows in Figure \ref{4abeven_d5}. For example, the first and second tilings in the first row of Figure \ref{timezone4ab_d5} are the two loops at the meridian $aaaaa$. The third is the arrow from $aaaaa$ to $aabab$. The fourth is the arrow from $aaaaa$ to $baaab$. The remaining four are the four arrows out of $ababa$ (including the loop).

With the only exception of two timezone tilings with no $b$-edge in the boundary meridians, all the other timezone tilings in the even family are uniquely determined by the two boundary meridians. This means that, with the exception of the two loops at the node $aaaaa$, there is at most one arrow from one node to another. We note that a double arrow between $m_1$ and $m_2$ means one timezone tiling with $m_1$ and $m_2$ as left and right meridians, and another timezone tiling with $m_2$ and $m_1$ as left and right meridians.

The symmetry of $180^{\circ}$ rotation is manifested as the vertical flipping of Figure \ref{4abeven_d5}, followed by reversing the arrow direction. An edge congruent earth map tiling is given by a directed cycle of the graph.

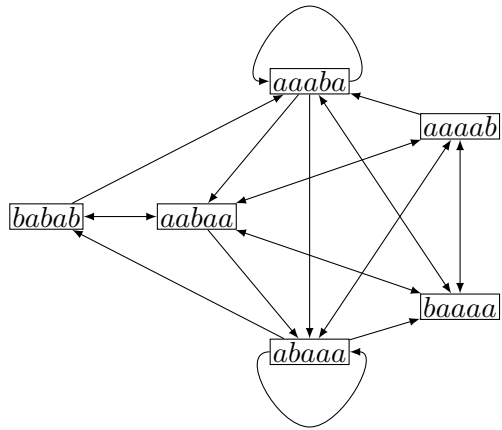
\begin{figure}[htp]
\centering
\begin{tikzpicture}[>=latex]

\node[rectangle,draw,inner sep=1] (A) at (0,0) {\small $aabaa$};
\node[rectangle,draw,inner sep=1] (B1) at (1.5,1.8) {\small $aaaba$};
\node[rectangle,draw,inner sep=1] (B-1) at (1.5,-1.8) {\small $abaaa$};
\node[rectangle,draw,inner sep=1] (C1) at (3.5,1.2) {\small $aaaab$};
\node[rectangle,draw,inner sep=1] (C-1) at (3.5,-1.2) {\small $baaaa$};
\node[rectangle,draw,inner sep=1] (E) at (-2,0) {\small $babab$};

\draw
	(A) edge[->] (B-1)
	(A) edge[<->] (C1)
	(A) edge[<->] (C-1)
	(B1) edge[->] (A)
	(B1) edge[->] (B-1)
	(B1) edge[<->] (C-1)
	(B-1) edge[->] (C-1)
	(B-1) edge[<->] (C1)
	(C1) edge[->] (B1)
	(C1) edge[<->] (C-1)
	(E) edge[<->] (A)
	(E) edge[->] (B1)
	(B-1) edge[->] (E);

\draw[<-]
	(B1) to[out=180,in=180] (1.5,2.8) to[out=0,in=0] (B1);
\draw[<-]
	(B-1) to[out=0,in=0] (1.5,-2.8) to[out=180,in=180] (B-1);
	
\end{tikzpicture}
\caption{Odd family of earth map tiling, distance $5$, edge combination $a^4b$.}
\label{4abodd_d5}
\end{figure}

Similar description applies to the odd family, which means $1$ or $3$ $b$-edges in each boundary meridian. The second row of Figure \ref{timezone4ab_d5} shows some such examples. The first three tilings in the second row give three arrows at $babab$ in Figure \ref{4abodd_d5}. The fourth tiling is the arrow from $baaaa$ to $aaaab$ (the reverse arrow is not related to this tiling). The last four tilings are all the arrows out of $aaaba$.

The symmetry of $180^{\circ}$ rotation is again the vertical flipping of Figure \ref{4abodd_d5}, followed by reversing the arrow direction. An edge congruent earth map tiling is also given by a directed cycle of the graph.

\subsection*{Distance $4$}

We need to consider meridian and core parts, both having the symmetries of vertical and horizontal flippings, and $180^{\circ}$ rotation. The tiling is obtained by alternatively gluing the meridian and core parts together. 

\begin{figure}[htp]
\centering
\begin{tikzpicture}[>=latex]

\foreach \x in {-1,1}
\foreach \y in {-1,1}
\foreach \a in {0,...,5}
{
\begin{scope}[xshift=1.2*\a cm]

\coordinate  (A\x\y X\a) at (0.25*\x,0.7*\y);
\coordinate  (B\x\y X\a) at (0.25*\x,0.35*\y);
\coordinate  (C\x X\a) at (0.15*\x,0);

\draw
	(A\x\y X\a) -- (B\x\y X\a) -- (C\x X\a);

\end{scope}
}

\foreach \a in {0,...,5}
\draw 
	(C1X\a) -- (C-1X\a);

\draw[dotted, xshift=3 cm]
	(0,-0.7) -- (0,0.7);

\draw[line width=1.8]
	(C1X0) -- (C-1X0)
	(B11X1) -- (C1X1)
	(A1-1X1) -- (B1-1X1)
	(A11X2) -- (B11X2)
	(A1-1X2) -- (B1-1X2)
	(B11X3) -- (C1X3)
	(B-1-1X3) -- (C-1X3)
	(B11X4) -- (C1X4)
	(A-1-1X4) -- (B-1-1X4)
	(A11X5) -- (B11X5)
	(A-1-1X5) -- (B-1-1X5);

\end{tikzpicture}
\caption{Meridian part, distance $4$, edge combination $a^4b$.}
\label{timezone4ab_d4m}
\end{figure}
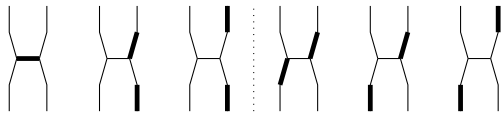

\begin{figure}[htp]
\centering
\begin{tikzpicture}[>=latex]


\node[rectangle,draw,inner sep=2] (A) at (-1.5,0) {\small $aaaa$};
\node[rectangle,draw,inner sep=1] (B1) at (0,1.5) {\small $abab$};
\node[rectangle,draw,inner sep=1] (B-1) at (0,-1.5) {\small $baba$};
\node[rectangle,draw,inner sep=1] (C) at (1.5,0) {\small $baab$};	

\draw[<->]
	(B1) to[out=210,in=60] (A);
\node[draw,shape=circle,fill=white,inner sep=0] at (-0.9,0.9) {\small 10};	
\draw[<->,dashed]
	(B1) to[out=240,in=30] (A);
	
\draw[<->]
	(B-1) to[out=150,in=-60] (A);
\node[draw,shape=circle,fill=white,inner sep=0] at (-0.9,-0.9) {\small 10};	
\draw[<->,dashed]
	(B-1) to[out=120,in=-30] (A);

\draw[<->]
	(B1) to[out=-30,in=120] (C);
\node[draw,shape=circle,fill=white,inner sep=1] at (0.9,0.9) {\small 2};	

\draw[<->]
	(B-1) to[out=30,in=240] (C);
\node[draw,shape=circle,fill=white,inner sep=1] at (0.9,-0.9) {\small 2};	

\draw[<->]
	(A) to[out=10,in=170] (C);
\node[draw,shape=circle,fill=white,inner sep=1] at (0.3,0.2) {\small 5};	
\draw[<->,dashed]
	(A) to[out=-10,in=190] (C);

\draw[<->]
	(B1) -- node[draw,shape=circle,fill=white,inner sep=1,near start] {\small 3} (B-1);

\draw[->]
	(A) to[out=100,in=90] (-2.8,0) to[out=270,in=260] (A);
\node[draw,shape=circle,fill=white,inner sep=0] at (-2.8,0) {\small 29};	
\draw[->,dashed]
	(A) to[out=120,in=90] (-2.4,0) to[out=270,in=240] (A);
		
\draw[->]
	(C) to[out=270,in=270] (2.8,0) to[out=90,in=80] (C);
\node[draw,shape=circle,fill=white,inner sep=1] at (2.8,0) {\small 1};		
\draw[->]
	(B1) to[out=0,in=0] (0,2.5) to[out=180,in=180] (B1);	
\node[draw,shape=circle,fill=white,inner sep=1] at (0,2.5) {\small 5};	

\draw[<-]
	(B-1) to[out=180,in=180] (0,-2.5) to[out=0,in=0] (B-1);
\node[draw,shape=circle,fill=white,inner sep=1] at (0,-2.5) {\small 5};	


\begin{scope}[xshift=7cm]

\node[rectangle,draw,inner sep=1] (A) at (-1.5,0) {\small $abaa$};
\node[rectangle,draw,inner sep=1] (B1) at (0,1.5) {\small $aaab$};
\node[rectangle,draw,inner sep=1] (B-1) at (0,-1.5) {\small $baaa$};
\node[rectangle,draw,inner sep=1] (C) at (1.5,0) {\small $aaba$};	

\draw[<->]
	(B1) to[out=210,in=60] (A);
\node[draw,shape=circle,fill=white,inner sep=1] at (-0.9,0.9) {\small 5};	
\draw[<->,dashed]
	(B1) to[out=240,in=30] (A);
	
\draw[<->]
	(B-1) to[out=150,in=-60] (A);
\node[draw,shape=circle,fill=white,inner sep=1] at (-0.9,-0.9) {\small 5};	

\draw[<->]
	(B1) to[out=-30,in=120] (C);
\node[draw,shape=circle,fill=white,inner sep=1] at (0.9,0.9) {\small 5};	

\draw[<->]
	(B-1) to[out=30,in=240] (C);
\node[draw,shape=circle,fill=white,inner sep=1] at (0.9,-0.9) {\small 5};
\draw[<->,dashed]
	(B-1) to[out=60,in=210] (C);	

\draw[<->]
	(A) to[out=10,in=170] (C);
\node[draw,shape=circle,fill=white,inner sep=0] at (0.5,0.2) {\small 10};	
\draw[<->,dashed]
	(A) to[out=-10,in=190] (C);
	
\draw[<->]
	(B1) to[out=260,in=100] (B-1);
\node[draw,shape=circle,fill=white,inner sep=1] at (-0.2,-0.5) {\small 5};	
\draw[<->,dashed]
	(B1) to[out=280,in=80] (B-1);

\draw[->]
	(A) to[out=100,in=90] (-2.8,0) to[out=270,in=260] (A);
\node[draw,shape=circle,fill=white,inner sep=0] at (-2.8,0) {\small 10};	
		
\draw[->]
	(C) to[out=270,in=270] (2.8,0) to[out=90,in=80] (C);
\node[draw,shape=circle,fill=white,inner sep=0] at (2.8,0) {\small 10};		
\draw[->]
	(B1) to[out=0,in=0] (0,2.5) to[out=180,in=180] (B1);	
\node[draw,shape=circle,fill=white,inner sep=1] at (0,2.5) {\small 5};	

\draw[<-]
	(B-1) to[out=180,in=180] (0,-2.5) to[out=0,in=0] (B-1);
\node[draw,shape=circle,fill=white,inner sep=1] at (0,-2.5) {\small 5};

\end{scope}
		
\end{tikzpicture}
\caption{Earth map tiling, distance $4$, edge combination $a^4b$.}
\label{4ab_d4}
\end{figure}
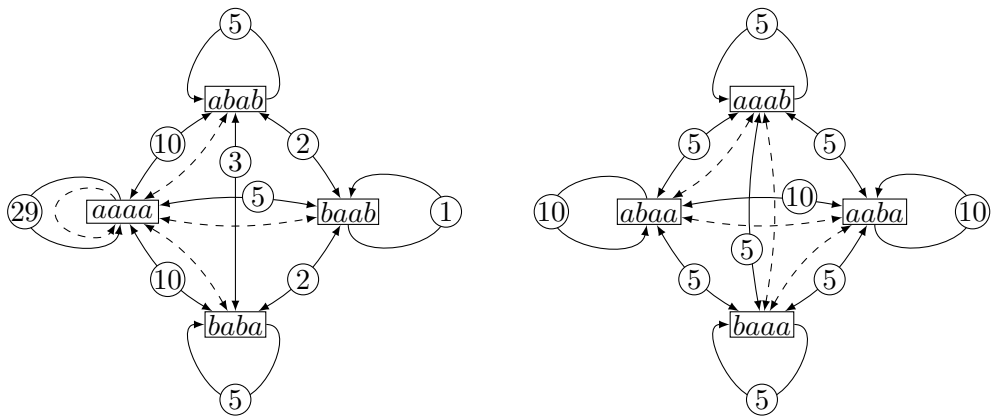

Up to symmetry, all the edge congruent tilings of the meridian part are given by Figure \ref{timezone4ab_d4m}. The first three form the even family because of the even number (i.e., $0$ or $2$) of $b$-edges in each boundary meridian. The three tilings and their symmetry transformations are represented by the dashed arrows on the left of Figure \ref{4ab_d4}. For example, the second meridian part tiling is the dashed arrow from $aaaa$ to $abab$, and its vertical flipping is the dashed arrow from $aaaa$ to $baba$. 

Similarly, the last three tilings in Figure \ref{timezone4ab_d4m} form the odd family because there is only one $b$-edge in each boundary meridian. The three tilings and their symmetry transformations are represented by the dashed arrows on the right of Figure \ref{4ab_d4}. 

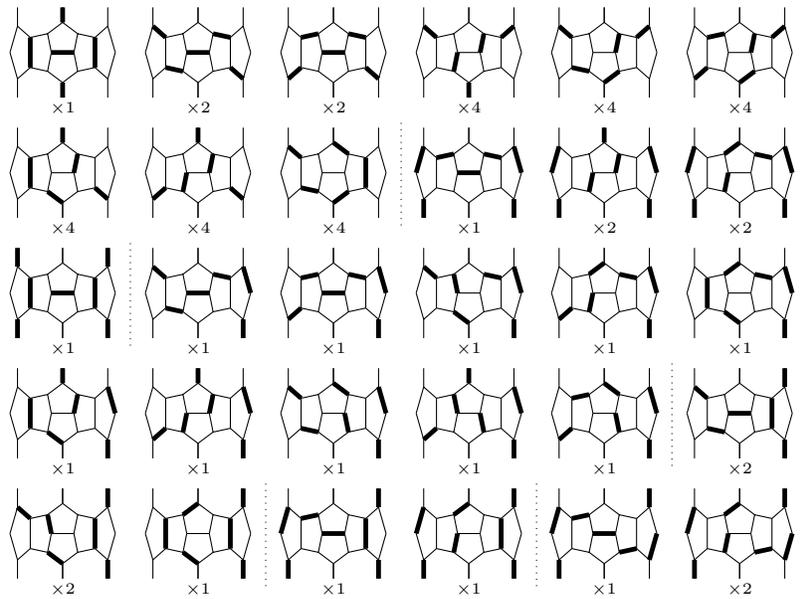
\begin{figure}[htp]
\centering
\begin{tikzpicture}[>=latex]

\foreach \x in {-1,1}
\foreach \y in {-1,1}
\foreach \a in {0,...,7}
\foreach \b in {0,...,8}
{
\begin{scope}[shift={(1.8*\a cm, -1.6*\b cm)}]

\coordinate  (A\y X\a\b) at (0,0.6*\y);
\coordinate  (B\y X\a\b) at (0,0.4*\y);
\coordinate  (C\x X\a\b) at (0.15*\x,0);
\coordinate  (D\x\y X\a\b) at (0.2*\x,0.25*\y);
\coordinate  (E\x\y X\a\b) at (0.43*\x,0.2*\y);
\coordinate  (F\x\y X\a\b) at (0.6*\x,0.35*\y);
\coordinate  (G\x\y X\a\b) at (0.6*\x,0.6*\y);
\coordinate  (H\x X\a\b) at (0.7*\x,0);

\coordinate  (P\a\b) at (0,-0.75);
	
\end{scope}
}

\foreach \a in {0,...,5}
\foreach \b in {0,...,4}
{

\foreach \x in {-1,1}
{

\foreach \y in {-1,1}
\draw
	(A\y X\a\b) -- (B\y X\a\b) -- (D\x\y X\a\b) -- (C\x X\a\b) 
	(D\x\y X\a\b) -- (E\x\y X\a\b) -- (F\x\y X\a\b)
	(G\x\y X\a\b) -- (F\x\y X\a\b) -- (H\x X\a\b);

\draw
	(E\x 1X\a\b) -- (E\x -1X\a\b);
	
}

\draw
	(C1X\a\b) -- (C-1X\a\b);
	
}


\foreach \u / \v in {3/1,6/1,1/2,5/3,2/4,4/4}
\draw[dotted, shift={(1.8*\u cm,-1.6*\v cm)}]
	(-0.9,-0.7) -- (-0.9,0.7);


\node at (P00) {\tiny $\times 1$};
\node at (P10) {\tiny $\times 2$};
\node at (P20) {\tiny $\times 2$};
\node at (P30) {\tiny $\times 4$};
\node at (P40) {\tiny $\times 4$};
\node at (P50) {\tiny $\times 4$};
\node at (P01) {\tiny $\times 4$};
\node at (P11) {\tiny $\times 4$};
\node at (P21) {\tiny $\times 4$};
\node at (P31) {\tiny $\times 1$};
\node at (P41) {\tiny $\times 2$};
\node at (P51) {\tiny $\times 2$};
\node at (P02) {\tiny $\times 1$};
\node at (P12) {\tiny $\times 1$};
\node at (P22) {\tiny $\times 1$};
\node at (P32) {\tiny $\times 1$};
\node at (P42) {\tiny $\times 1$};
\node at (P52) {\tiny $\times 1$};
\node at (P03) {\tiny $\times 1$};
\node at (P13) {\tiny $\times 1$};
\node at (P23) {\tiny $\times 1$};
\node at (P33) {\tiny $\times 1$};
\node at (P43) {\tiny $\times 1$};
\node at (P53) {\tiny $\times 2$};
\node at (P04) {\tiny $\times 2$};
\node at (P14) {\tiny $\times 1$};
\node at (P24) {\tiny $\times 1$};
\node at (P34) {\tiny $\times 1$};
\node at (P44) {\tiny $\times 1$};
\node at (P54) {\tiny $\times 2$};


\draw[line width=1.8]
	(C1X00) -- (C-1X00)	
	(A1X00) -- (B1X00)
	(A-1X00) -- (B-1X00)
	(E11X00) -- (E1-1X00)
	(E-11X00) -- (E-1-1X00)
	(C1X10) -- (C-1X10)	
	(D11X10) -- (E11X10)
	(D-1-1X10) -- (E-1-1X10)
	(E1-1X10) -- (F1-1X10)
	(E-11X10) -- (F-11X10)
	(C1X20) -- (C-1X20)	
	(D11X20) -- (E11X20)
	(E1-1X20) -- (F1-1X20)
	(E-1-1X20) -- (F-1-1X20)
	(D-11X20) -- (E-11X20)
	(C1X30) -- (D11X30)	
	(E11X30) -- (F11X30)
	(A-1X30) -- (B-1X30)
	(C-1X30) -- (D-1-1X30)
	(E-11X30) -- (F-11X30)
	(C1X40) -- (D11X40)	
	(E11X40) -- (F11X40)
	(B-1X40) -- (D1-1X40)
	(D-1-1X40) -- (E-1-1X40)
	(E-11X40) -- (F-11X40)
	(C1X50) -- (D11X50)	
	(E11X50) -- (F11X50)
	(B-1X50) -- (D1-1X50)
	(E-1-1X50) -- (F-1-1X50)
	(D-11X50) -- (E-11X50)
	(C1X01) -- (D11X01)	
	(E1-1X01) -- (F1-1X01)
	(A1X01) -- (B1X01)
	(E-11X01) -- (E-1-1X01)
	(B-1X01) -- (D-1-1X01)
	(C1X11) -- (D11X11)	
	(E1-1X11) -- (F1-1X11)
	(A1X11) -- (B1X11)
	(C-1X11) -- (D-1-1X11)
	(E-1-1X11) -- (F-1-1X11)
	(B1X21) -- (D11X21)	
	(E11X21) -- (E1-1X21)
	(B-1X21) -- (D1-1X21)
	(D-1-1X21) -- (E-1-1X21)
	(E-11X21) -- (F-11X21)
	;


\draw[line width=1.8]
	(F11X31) -- (H1X31)	
	(F-11X31) -- (H-1X31)
	(F1-1X31) -- (G1-1X31)
	(F-1-1X31) -- (G-1-1X31)
	(C1X31) -- (C-1X31)
	(D11X31) -- (E11X31)
	(D-11X31) -- (E-11X31)
	(F11X41) -- (H1X41)	
	(F-11X41) -- (H-1X41)
	(F1-1X41) -- (G1-1X41)
	(F-1-1X41) -- (G-1-1X41)
	(D-1-1X41) -- (C-1X41)
	(D11X41) -- (C1X41)
	(A1X41) -- (B1X41)
	(F11X51) -- (H1X51)	
	(F-11X51) -- (H-1X51)
	(F1-1X51) -- (G1-1X51)
	(F-1-1X51) -- (G-1-1X51)
	(D-1-1X51) -- (C-1X51)
	(D11X51) -- (E11X51)
	(B1X51) -- (D-11X51)
	;


\draw[line width=1.8]
	(F11X02) -- (G11X02)	
	(F1-1X02) -- (G1-1X02)
	(F-11X02) -- (G-11X02)
	(F-1-1X02) -- (G-1-1X02)
	(E11X02) -- (E1-1X02)
	(E-11X02) -- (E-1-1X02)
	(C1X02) -- (C-1X02)
	;	
	

\draw[line width=1.8]
	(F11X12) -- (H1X12)	
	(F1-1X12) -- (G1-1X12)
	(C1X12) -- (C-1X12)
	(D11X12) -- (E11X12)
	(D-1-1X12) -- (E-1-1X12)
	(E-11X12) -- (F-11X12)
	(F11X22) -- (H1X22)	
	(F1-1X22) -- (G1-1X22)
	(C1X22) -- (C-1X22)
	(D11X22) -- (E11X22)
	(E-1-1X22) -- (F-1-1X22)
	(D-11X22) -- (E-11X22)
	(F11X32) -- (H1X32)	
	(F1-1X32) -- (G1-1X32)
	(D11X32) -- (E11X32)
	(C-1X32) -- (D-11X32)
	(B-1X32) -- (D-1-1X32)
	(E-11X32) -- (F-11X32)
	(F11X42) -- (H1X42)	
	(F1-1X42) -- (G1-1X42)
	(D11X42) -- (E11X42)
	(B1X42) -- (D-11X42)
	(C-1X42) -- (D-1-1X42)
	(E-1-1X42) -- (F-1-1X42)
	(F11X52) -- (H1X52)	
	(F1-1X52) -- (G1-1X52)
	(D11X52) -- (E11X52)
	(B1X52) -- (D-11X52)
	(B-1X52) -- (D-1-1X52)
	(E-11X52) -- (E-1-1X52)
	(F11X03) -- (H1X03)	
	(F1-1X03) -- (G1-1X03)
	(C1X03) -- (D11X03)
	(A1X03) -- (B1X03)
	(E-11X03) -- (E-1-1X03)
	(B-1X03) -- (D-1-1X03)
	(F11X13) -- (H1X13)	
	(F1-1X13) -- (G1-1X13)
	(C1X13) -- (D11X13)
	(A1X13) -- (B1X13)
	(E-1-1X13) -- (F-1-1X13)
	(C-1X13) -- (D-1-1X13)
	(F11X23) -- (H1X23)	
	(F1-1X23) -- (G1-1X23)
	(C1X23) -- (D1-1X23)
	(D-1-1X23) -- (E-1-1X23)
	(E-11X23) -- (F-11X23)
	(B1X23) -- (D11X23)
	(F11X33) -- (H1X33)	
	(F1-1X33) -- (G1-1X33)
	(C1X33) -- (D1-1X33)
	(E-1-1X33) -- (F-1-1X33)
	(A1X33) -- (B1X33)
	(C-1X33) -- (D-11X33)
	(F11X43) -- (H1X43)	
	(F1-1X43) -- (G1-1X43)
	(C1X43) -- (D1-1X43)
	(E-1-1X43) -- (F-1-1X43)
	(B1X43) -- (D11X43)
	(D-11X43) -- (E-11X43)
	;


\draw[line width=1.8]
	(F11X53) -- (G11X53)	
	(F1-1X53) -- (G1-1X53)
	(E11X53) -- (E1-1X53)
	(C1X53) -- (C-1X53)
	(D-1-1X53) -- (E-1-1X53)
	(E-11X53) -- (F-11X53)
	(F11X04) -- (G11X04)	
	(F1-1X04) -- (G1-1X04)
	(E11X04) -- (E1-1X04)
	(C-1X04) -- (D-11X04)
	(B-1X04) -- (D-1-1X04)
	(E-11X04) -- (F-11X04)
	(F11X14) -- (G11X14)	
	(F1-1X14) -- (G1-1X14)
	(E11X14) -- (E1-1X14)
	(E-11X14) -- (E-1-1X14)
	(B1X14) -- (D-11X14)
	(B-1X14) -- (D-1-1X14)
	;


\draw[line width=1.8]
	(F11X24) -- (G11X24)	
	(F-11X24) -- (H-1X24)
	(F1-1X24) -- (G1-1X24)
	(F-1-1X24) -- (G-1-1X24)
	(E11X24) -- (E1-1X24)
	(C1X24) -- (C-1X24)
	(D-11X24) -- (E-11X24)
	(F11X34) -- (G11X34)	
	(F-11X34) -- (H-1X34)
	(F1-1X34) -- (G1-1X34)
	(F-1-1X34) -- (G-1-1X34)
	(E11X34) -- (E1-1X34)
	(C-1X34) -- (D-1-1X34)
	(B1X34) -- (D-11X34)
	;


\draw[line width=1.8]
	(F11X44) -- (G11X44)	
	(F-11X44) -- (H-1X44)
	(F1-1X44) -- (H1X44)
	(F-1-1X44) -- (G-1-1X44)
	(C1X44) -- (C-1X44)
	(D1-1X44) -- (E1-1X44)
	(D-11X44) -- (E-11X44)
	(F11X54) -- (G11X54)	
	(F-11X54) -- (H-1X54)
	(F1-1X54) -- (H1X54)
	(F-1-1X54) -- (G-1-1X54)
	(C-1X54) -- (D-1-1X54)
	(D1-1X54) -- (E1-1X54)
	(B1X54) -- (D-11X54)
	;


\end{tikzpicture}
\caption{Core part of even family, distance $4$, edge combination $a^4b$.}
\label{timezone4ab_d4even}
\end{figure}

Now consider the tilings of the core part. Up to symmetry, the complete list for the even family (i.e., $0$ or $2$ $b$-edges in each boundary meridian) is given by Figure \ref{timezone4ab_d4even}. Below each tiling, $\times m$ means that, if we apply symmetries that do not change the boundary meridians, then we get $m$ tilings. Such ``multiplicity'' already appeared before. For example, for the first family of the earth map tiling of distance $1$ with edge combination $a^3b^2$, the tiling $\binom{aaa}{aba}$ has multiplicity $4$ because applying the symmetries also gives $\binom{aba}{aaa}$, $\binom{aaa}{aab}$, $\binom{aab}{aaa}$.
\begin{itemize}
\item Nos.1-9 have the meridian $aaaa$ as both boundaries. There are total of $29$ core part tilings.
\item Nos.10-12 have the meridian $abab$ as both boundaries. There are total of $5$ core part tilings. Their vertical flippings give $5$ core part tilings with $baba$ as both boundaries.
\item No.13 has the meridian $baab$ as both boundaries. There is only $1$ core part tiling. 
\item Nos.14-23 have the meridians $aaaa$ and $abab$ as the left and right boundaries. There are total of $10$ core part tilings. Their vertical flippings give $10$ tilings with $aaaa$ and $baba$ as boundaries. Their horizontal flippings give $10$ tilings with $abab$ and $aaaa$ as boundaries. Their $180^{\circ}$ rotations give $10$ tilings with $baba$ and $aaaa$ as boundaries.
\item Nos.24-26 have the meridians $aaaa$ and $baab$ as the left and right boundaries. There are total of $5$ core part tilings. Their horizontal flippings give $5$ tilings with $baab$ and $aaaa$ as boundaries.
\item Nos.27-28 have the meridians $abab$ and $baab$ as the left and right boundaries. There are total of $2$ core part tilings. Their vertical flippings, horizontal flippings, and $180^{\circ}$ rotations give $2$ tilings each for the boundary combinations $baba$ and $baab$, $baab$ and $abab$, $baab$ and $baba$.
\item Nos.29-30 have the meridians $abab$ and $baba$ as the left and right boundaries. There are total of $3$ core part tilings. Their vertical flippings give $2$ tilings with $baba$ and $abab$ as boundaries.
\end{itemize}

Like the meridian part, we represent each core part tiling (and any symmetric transformation) in Figure \ref{timezone4ab_d4even} by a solid arrow on the left of Figure \ref{4ab_d4}. Since we often have too many arrows between two given nodes, we simply draw one arrow with a number indicating the actual number of arrows. For example, we have $10$ arrows from $aaaa$ and $abab$ representing the tilings Nos.14-23 in Figure \ref{timezone4ab_d4even}. The horizontal flippings of these tilings give $10$ arrows from $abab$ and $aaaa$.

The horizontal flipping converts a tiling with $m_1$ and $m_2$ as the left and right boundary meridians and to a tiling with $m_2$ and $m_1$ as the left and right boundary meridians. Since the horizontal flipping means reversing the arrow direction in Figure \ref{4ab_d4}, all the arrows (dashed as well as solid) are double arrows. We also note that the vertical flipping of the the meridian and core part tilings is manifested as the vertical flipping of the graph. 

An edge congruent earth map tiling in the even family is then a directed cycle in the graph on the left of Figure \ref{4ab_d4}, such that the solid and dashed arrows appear alternatively.

\begin{figure}[htp]
\centering
\begin{tikzpicture}[>=latex]

\foreach \x in {-1,1}
\foreach \y in {-1,1}
\foreach \a in {0,...,7}
\foreach \b in {0,...,8}
{
\begin{scope}[shift={(1.8*\a cm, -1.6*\b cm)}]

\coordinate  (A\y X\a\b) at (0,0.6*\y);
\coordinate  (B\y X\a\b) at (0,0.4*\y);
\coordinate  (C\x X\a\b) at (0.15*\x,0);
\coordinate  (D\x\y X\a\b) at (0.2*\x,0.25*\y);
\coordinate  (E\x\y X\a\b) at (0.43*\x,0.2*\y);
\coordinate  (F\x\y X\a\b) at (0.6*\x,0.35*\y);
\coordinate  (G\x\y X\a\b) at (0.6*\x,0.6*\y);
\coordinate  (H\x X\a\b) at (0.7*\x,0);

\coordinate  (P\a\b) at (0,-0.75);
	
\end{scope}
}

\foreach \x in {-1,1}
\foreach \y in {-1,1}
\foreach \a in {0,...,5}
\foreach \b in {0,...,4}
\draw
	(A\y X\a\b) -- (B\y X\a\b) -- (D\x\y X\a\b) -- (C\x X\a\b) 
	(D\x\y X\a\b) -- (E\x\y X\a\b) -- (F\x\y X\a\b)
	(G\x\y X\a\b) -- (F\x\y X\a\b) -- (H\x X\a\b)
	(C1X\a\b) -- (C-1X\a\b)
	(E\x 1X\a\b) -- (E\x -1X\a\b);


\foreach \u / \v in {3/0,3/1,2/2,1/3,5/3}
\draw[dotted, shift={(1.8*\u cm,-1.6*\v cm)}]
	(-0.9,-0.7) -- (-0.9,0.7);


\node at (P00) {\tiny $\times 1$};
\node at (P10) {\tiny $\times 2$};
\node at (P20) {\tiny $\times 2$};
\node at (P30) {\tiny $\times 1$};
\node at (P40) {\tiny $\times 1$};
\node at (P50) {\tiny $\times 2$};
\node at (P01) {\tiny $\times 2$};
\node at (P11) {\tiny $\times 2$};
\node at (P21) {\tiny $\times 2$};
\node at (P31) {\tiny $\times 1$};
\node at (P41) {\tiny $\times 1$};
\node at (P51) {\tiny $\times 1$};
\node at (P02) {\tiny $\times 1$};
\node at (P12) {\tiny $\times 1$};
\node at (P22) {\tiny $\times 1$};
\node at (P32) {\tiny $\times 1$};
\node at (P42) {\tiny $\times 1$};
\node at (P52) {\tiny $\times 1$};
\node at (P03) {\tiny $\times 1$};
\node at (P13) {\tiny $\times 1$};
\node at (P23) {\tiny $\times 1$};
\node at (P33) {\tiny $\times 2$};
\node at (P43) {\tiny $\times 1$};
\node at (P53) {\tiny $\times 2$};
\node at (P04) {\tiny $\times 1$};
\node at (P14) {\tiny $\times 2$};
\node at (P24) {\tiny $\times 1$};
\node at (P34) {\tiny $\times 1$};
\node at (P44) {\tiny $\times 2$};
\node at (P54) {\tiny $\times 1$};


\draw[line width=1.8]
	(F-1-1X00) -- (G-1-1X00)	
	(F1-1X00) -- (G1-1X00)
	(C-1X00) -- (C1X00)
	(E-11X00) -- (E-1-1X00)
	(E11X00) -- (E1-1X00)
	(A1X00) -- (B1X00)	
	(F-1-1X10) -- (G-1-1X10)	
	(F1-1X10) -- (G1-1X10)
	(C-1X10) -- (D-1-1X10)
	(E-11X10) -- (F-11X10)
	(B1X10) -- (D11X10)
	(E11X10) -- (E1-1X10)	
	(F-1-1X20) -- (G-1-1X20)	
	(F1-1X20) -- (G1-1X20)
	(C-1X20) -- (D-1-1X20)
	(E-11X20) -- (F-11X20)
	(C1X20) -- (D11X20)
	(E11X20) -- (F11X20);


\draw[line width=1.8]	
	(H-1X30) -- (F-1-1X30)	
	(H1X30) -- (F1-1X30)
	(C-1X30) -- (C1X30)
	(A1X30) -- (B1X30)
	(E-1-1X30) -- (D-1-1X30)
	(E1-1X30) -- (D1-1X30)
	(H-1X40) -- (F-1-1X40)	
	(H1X40) -- (F1-1X40)
	(C-1X40) -- (C1X40)
	(D-11X40) -- (E-11X40)	
	(D11X40) -- (E11X40)
	(A-1X40) -- (B-1X40)
	(H-1X50) -- (F-1-1X50)	
	(H1X50) -- (F1-1X50)
	(C-1X50) -- (D-1-1X50)
	(A-1X50) -- (B-1X50)
	(A1X50) -- (B1X50)
	(C1X50) -- (D11X50)	
	(H-1X01) -- (F-1-1X01)	
	(H1X01) -- (F1-1X01)
	(C-1X01) -- (D-1-1X01)
	(A-1X01) -- (B-1X01)	
	(D-11X01) -- (B1X01)
	(D11X01) -- (E11X01)
	(H-1X11) -- (F-1-1X11)	
	(H1X11) -- (F1-1X11)
	(B-1X11) -- (D-1-1X11)
	(D1-1X11) -- (E1-1X11)	
	(D-11X11) -- (E-11X11)
	(B1X11) -- (D11X11)
	(H-1X21) -- (F-1-1X21)	
	(H1X21) -- (F1-1X21)
	(B-1X21) -- (D-1-1X21)
	(D1-1X21) -- (E1-1X21)
	(C-1X21) -- (D-11X21)
	(A1X21) -- (B1X21)
	;


\draw[line width=1.8]
	(F-1-1X31) -- (G-1-1X31)	
	(F1-1X31) -- (H1X31)
	(C1X31) -- (C-1X31)	
	(D1-1X31) -- (E1-1X31)
	(A1X31) -- (B1X31)
	(E-11X31) -- (E-1-1X31)
	(F-1-1X41) -- (G-1-1X41)	
	(F1-1X41) -- (H1X41)
	(D1-1X41) -- (E1-1X41)
	(D-1-1X41) -- (C-1X41)	
	(F-11X41) -- (E-11X41)
	(D11X41) -- (B1X41)
	(F-1-1X51) -- (G-1-1X51)	
	(F1-1X51) -- (H1X51)
	(D1-1X51) -- (B-1X51)
	(D11X51) -- (C1X51)	
	(A1X51) -- (B1X51)
	(E-11X51) -- (E-1-1X51)
	(F-1-1X02) -- (G-1-1X02)	
	(F1-1X02) -- (H1X02)
	(D1-1X02) -- (B-1X02)	
	(D11X02) -- (E11X02)	
	(E-11X02) -- (F-11X02)
	(C-1X02) -- (D-11X02)
	(F-1-1X12) -- (G-1-1X12)	
	(F1-1X12) -- (H1X12)
	(D1-1X12) -- (B-1X12)
	(D11X12) -- (E11X12)
	(D-11X12) -- (B1X12)
	(E-11X12) -- (E-1-1X12)
	;


\draw[line width=1.8]

	(F-1-1X22) -- (G-1-1X22)  
	(F11X22) -- (H1X22)	
	(C1X22) -- (C-1X22)
	(D1-1X22) -- (E1-1X22)
	(A1X22) -- (B1X22)
	(E-11X22) -- (E-1-1X22)
	(F-1-1X32) -- (G-1-1X32)  
	(F11X32) -- (H1X32)	
	(D1-1X32) -- (E1-1X32)
	(C-1X32) -- (D-1-1X32)
	(F-11X32) -- (E-11X32)
	(B1X32) -- (D11X32)
	(F-1-1X42) -- (G-1-1X42)  
	(F11X42) -- (H1X42)	
	(D1-1X42) -- (B-1X42)
	(C1X42) -- (D11X42)
	(E-11X42) -- (E-1-1X42)
	(A1X42) -- (B1X42)
	(F-1-1X52) -- (G-1-1X52)  
	(F11X52) -- (H1X52)	
	(B-1X52) -- (D1-1X52)
	(D11X52) -- (E11X52)
	(F-11X52) -- (E-11X52)
	(D-11X52) -- (C-1X52)
	(F-1-1X03) -- (G-1-1X03)  
	(F11X03) -- (H1X03)	
	(B-1X03) -- (D1-1X03)
	(D11X03) -- (E11X03)
	(E-11X03) -- (E-1-1X03)
	(B1X03) -- (D-11X03)
	;
	

\draw[line width=1.8]	
	
	(F-1-1X13) -- (G-1-1X13)	
	(F11X13) -- (G11X13)
	(C-1X13) -- (D-1-1X13)
	(F-11X13) -- (E-11X13)
	(C1X13) -- (D11X13)	
	(E1-1X13) -- (F1-1X13)

	(F-1-1X23) -- (G-1-1X23)	
	(F11X23) -- (G11X23)
	(C1X23) -- (D1-1X23)
	(C-1X23) -- (D-11X23)	
	(E1-1X23) -- (F1-1X23)	
	(E-11X23) -- (F-11X23)
	(F-1-1X33) -- (G-1-1X33)	
	(F11X33) -- (G11X33)
	(C1X33) -- (D1-1X33)
	(E1-1X33) -- (F1-1X33)
	(D-11X33) -- (B1X33)
	(E-11X33) -- (E-1-1X33)
	(F-1-1X43) -- (G-1-1X43)	
	(F11X43) -- (G11X43)
	(D-11X43) -- (B1X43)
	(E-11X43) -- (E-1-1X43)
	(D1-1X43) -- (B-1X43)
	(E11X43) -- (E1-1X43)	
	;


\draw[line width=1.8]
	(F-1-1X53) -- (H-1X53)	
	(F11X53) -- (H1X53)
	(C1X53) -- (C-1X53)
	(A1X53) -- (B1X53)
	(D-1-1X53) -- (E-1-1X53)
	(D1-1X53) -- (E1-1X53)	
	(F-1-1X04) -- (H-1X04)	
	(F11X04) -- (H1X04)
	(C1X04) -- (D1-1X04)
	(C-1X04) -- (D-11X04)	
	(A1X04) -- (B1X04)
	(A-1X04) -- (B-1X04)
	(F-1-1X14) -- (H-1X14)	
	(F11X14) -- (H1X14)
	(C-1X14) -- (D-11X14)	
	(A1X14) -- (B1X14)
	(D1-1X14) -- (E1-1X14)
	(B-1X14) -- (D-1-1X14)
	(F-1-1X24) -- (H-1X24)	
	(F11X24) -- (H1X24)
	(E-11X24) -- (D-11X24)
	(E1-1X24) -- (D1-1X24)
	(B1X24) -- (D11X24)
	(B-1X24) -- (D-1-1X24)
	(F-1-1X34) -- (H-1X34)	
	(F11X34) -- (H1X34)
	(C1X34) -- (D11X34)
	(C-1X34) -- (D-1-1X34)
	(A1X34) -- (B1X34)
	(A-1X34) -- (B-1X34)
	(F-1-1X44) -- (H-1X44)	
	(F11X44) -- (H1X44)
	(C-1X44) -- (D-1-1X44)	
	(A-1X44) -- (B-1X44)
	(D11X44) -- (E11X44)
	(B1X44) -- (D-11X44)
	(F-1-1X54) -- (H-1X54)	
	(F11X54) -- (H1X54)
	(E-1-1X54) -- (D-1-1X54)
	(E11X54) -- (D11X54)
	(B-1X54) -- (D1-1X54)
	(B1X54) -- (D-11X54)
	;


\end{tikzpicture}
\caption{Core part of odd family, distance $4$, edge combination $a^4b$.}
\label{timezone4ab_d4odd}
\end{figure}
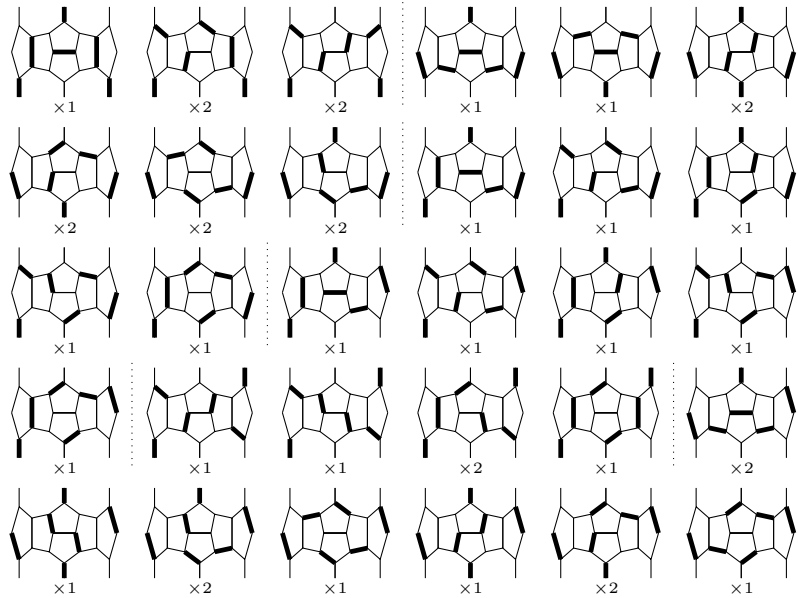

Similarly, up to symmetry, the complete list for the odd family (i.e., $1$ $b$-edge in each boundary meridian) is given by Figure \ref{timezone4ab_d4odd}. 
\begin{itemize}
\item Nos.1-3 have the meridian $aaab$ as both boundaries. There are total of $5$ core part tilings. Their vertical flippings give $5$ tilings with $baaa$ as both boundaries.
\item Nos.4-9 have the meridian $aaba$ as both boundaries. There are total of $10$ core part tilings. Their vertical flippings give $10$ tilings with $abaa$ as both boundaries.
\item Nos.10-14 has the meridian $aaab$ and $aaba$ as the left and right boundaries. There are $5$ core part tilings. Their vertical flippings, horizontal flippings, and $180^{\circ}$ rotations give $5$ tilings each for the boundary combinations $baaa$ and $abaa$, $aaba$ and $aaab$, $abaa$ and $baaa$. 
\item Nos.15-19 have the meridians $aaab$ and $abaa$ as the left and right boundaries. There are total of $5$ core part tilings. Their vertical flippings, horizontal flippings, and $180^{\circ}$ rotations give $5$ tilings each for the boundary combinations $baaa$ and $aaba$, $abaa$ and $aaab$, $aaba$ and $baaa$.
\item Nos.20-23 have the meridians $aaab$ and $baaa$ as the left and right boundaries. There are total of $5$ core part tilings. Their horizontal flippings give $5$ tilings with $baaa$ and $aaab$ as boundaries.
\item Nos.24-30 have the meridians $aaba$ and $abaa$ as the left and right boundaries. There are total of $10$ core part tilings. Their horizontal flippings give $10$ tilings with $abaa$ and $aaba$ as boundaries.
\end{itemize}

Like the even family, the data for the odd family is summarized on the right of Figure \ref{4ab_d4}. An edge congruent earth map tiling in the odd family is a directed cycle in this graph, such that the solid and dashed arrows appear alternatively.

\subsection*{Distance 3}

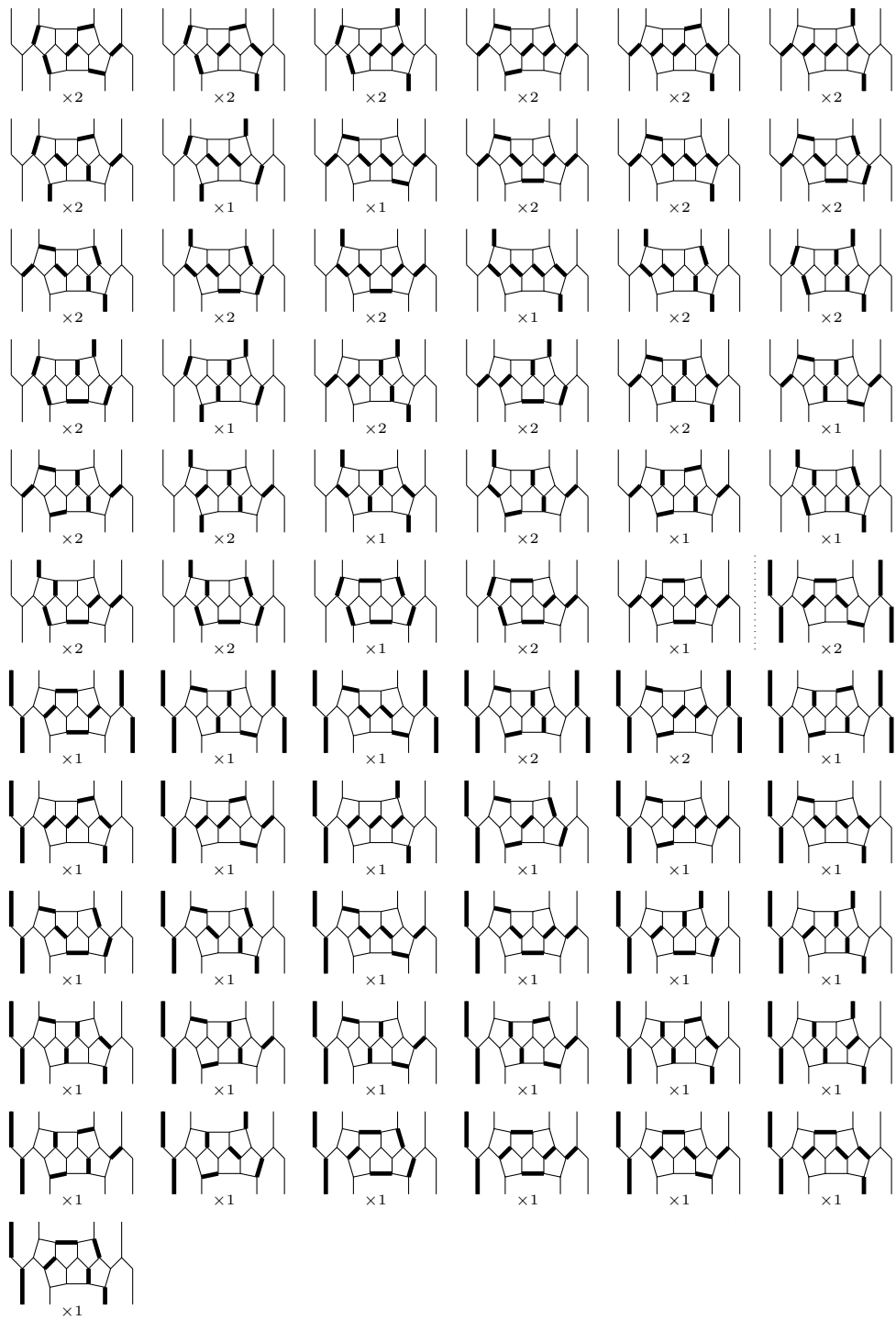
\begin{figure}[htp]
\centering
\begin{tikzpicture}[>=latex]

\foreach \x in {-1,1}
\foreach \y in {-1,1}
\foreach \a in {0,...,5}
\foreach \b in {0,...,11}
{
\begin{scope}[shift={(2.2*\a cm, -1.6*\b cm)}]

\coordinate  (A\x\y X\a\b) at (0.16*\x-0.08*\y,0.08*\y);
\coordinate  (B\x\y X\a\b) at (0.16*\x-0.08*\y,0.3*\y);
\coordinate  (C\x\y X\a\b) at (0.40*\x-0.08*\y,0.6*\y);
\coordinate  (D\x\y X\a\b) at (0.40*\x-0.08*\y,0.35*\y);
\coordinate  (E\x\y X\a\b) at (0.48*\x-0.08*\y,0.08*\y);
\coordinate  (F\x\y X\a\b) at (0.80*\x-0.08*\y,0.08*\y);
\coordinate  (G\x\y X\a\b) at (0.80*\x-0.08*\y,0.6*\y);

\coordinate  (P\a\b) at (0,-0.7);

\draw
	(A\x\y X\a\b) -- (B\x\y X\a\b)
	(C\x\y X\a\b) -- (D\x\y X\a\b) -- (B\x\y X\a\b)
	(D\x\y X\a\b) -- (E\x\y X\a\b)
	(F\x\y X\a\b) -- (G\x\y X\a\b);
	
\end{scope}
}

\foreach \a in {0,...,5}
\foreach \b in {0,...,11}
\draw
	(F1-1X\a\b) -- (F11X\a\b) -- (E1-1X\a\b) -- (E11X\a\b) -- (A1-1X\a\b) -- (A11X\a\b) -- (A-1-1X\a\b) -- (A-11X\a\b) -- (E-1-1X\a\b) -- (E-11X\a\b) -- (F-1-1X\a\b) -- (F-11X\a\b)
	(B11X\a\b) -- (B-11X\a\b)
	(B1-1X\a\b) -- (B-1-1X\a\b);


\foreach \u / \v in {5/5,6/6}
\draw[dotted, shift={(2.2*\u cm,-1.6*\v cm)}]
	(-1.1,-0.7) -- (-1.1,0.7);


\node at (P00) {\tiny $\times 2$};
\node at (P10) {\tiny $\times 2$};
\node at (P20) {\tiny $\times 2$};
\node at (P30) {\tiny $\times 2$};
\node at (P40) {\tiny $\times 2$};
\node at (P50) {\tiny $\times 2$};
\node at (P01) {\tiny $\times 2$};
\node at (P11) {\tiny $\times 1$};
\node at (P21) {\tiny $\times 1$};
\node at (P31) {\tiny $\times 2$};
\node at (P41) {\tiny $\times 2$};
\node at (P51) {\tiny $\times 2$};
\node at (P02) {\tiny $\times 2$};
\node at (P12) {\tiny $\times 2$};
\node at (P22) {\tiny $\times 2$};
\node at (P32) {\tiny $\times 1$};
\node at (P42) {\tiny $\times 2$};
\node at (P52) {\tiny $\times 2$};
\node at (P03) {\tiny $\times 2$};
\node at (P13) {\tiny $\times 1$};
\node at (P23) {\tiny $\times 2$};
\node at (P33) {\tiny $\times 2$};
\node at (P43) {\tiny $\times 2$};
\node at (P53) {\tiny $\times 1$};
\node at (P04) {\tiny $\times 2$};
\node at (P14) {\tiny $\times 2$};
\node at (P24) {\tiny $\times 1$};
\node at (P34) {\tiny $\times 2$};
\node at (P44) {\tiny $\times 1$};
\node at (P54) {\tiny $\times 1$};
\node at (P05) {\tiny $\times 2$};
\node at (P15) {\tiny $\times 2$};
\node at (P25) {\tiny $\times 1$};
\node at (P35) {\tiny $\times 2$};
\node at (P45) {\tiny $\times 1$};

\node at (P55) {\tiny $\times 2$};
\node at (P06) {\tiny $\times 1$};
\node at (P16) {\tiny $\times 1$};
\node at (P26) {\tiny $\times 1$};
\node at (P36) {\tiny $\times 2$};
\node at (P46) {\tiny $\times 2$};
\node at (P56) {\tiny $\times 1$};

\node at (P011) {\tiny $\times 1$};

\foreach \a in {0,1,...,5}
\foreach \b in {7,8,9,10}
{
\node at (P\a\b) {\tiny $\times 1$};
}

\fill[white,shift={(2 cm, -17.6 cm)}] (-0.8,-0.6) rectangle (10,0.6);

\draw[line width=1.8]
	(D-11X00) -- (E-11X00)	
	(E-1-1X00) -- (D-1-1X00)	
	(D1-1X00) -- (B1-1X00)	
	(D11X00) -- (B11X00)
	(A-1-1X00) -- (A11X00)
	(F11X00) -- (E1-1X00)
	
	(D-11X10) -- (E-11X10)	
	(E-1-1X10) -- (D-1-1X10)	
	(D1-1X10) -- (C1-1X10)	
	(D11X10) -- (B11X10)
	(A-1-1X10) -- (A11X10)
	(E11X10) -- (E1-1X10)
	
	(D-11X20) -- (E-11X20)	
	(E-1-1X20) -- (D-1-1X20)	
	(D1-1X20) -- (C1-1X20)	
	(D11X20) -- (C11X20)
	(A-1-1X20) -- (A11X20)
	(E11X20) -- (A1-1X20)

	(F-1-1X30) -- (E-11X30)	
	(B-1-1X30) -- (D-1-1X30)	
	(F11X30) -- (E1-1X30)	
	(D-11X30) -- (B-11X30)
	(A-1-1X30) -- (A11X30)
	(E11X30) -- (A1-1X30)
	
	(F-1-1X40) -- (E-11X40)	
	(B11X40) -- (D11X40)	
	(E11X40) -- (E1-1X40)	
	(D1-1X40) -- (C1-1X40)
	(A-1-1X40) -- (A11X40)
	(E-1-1X40) -- (A-11X40)

	(F-1-1X50) -- (E-11X50)	
	(C11X50) -- (D11X50)	
	(E11X50) -- (A1-1X50)	
	(D1-1X50) -- (C1-1X50)
	(A-1-1X50) -- (A11X50)
	(E-1-1X50) -- (A-11X50)

	(E-11X01) -- (D-11X01)	
	(A-11X01) -- (A-1-1X01)
	(B11X01) -- (D11X01)		
	(D-1-1X01) -- (C-1-1X01)
	(B1-1X01) -- (A1-1X01)
	(F11X01) -- (E1-1X01)
	
	(E-11X11) -- (D-11X11)	
	(A-11X11) -- (A-1-1X11)
	(C11X11) -- (D11X11)
	(D-1-1X11) -- (C-1-1X11)
	(A11X11) -- (A1-1X11)
	(D1-1X11) -- (E1-1X11)
	
	(B-11X21) -- (D-11X21)	
	(A-11X21) -- (A-1-1X21)
	(E1-1X21) -- (F11X21)
	(F-1-1X21) -- (E-11X21)
	(A11X21) -- (A1-1X21)
	(D1-1X21) -- (B1-1X21)

	(A-11X31) -- (A-1-1X31)	
	(B1-1X31) -- (B-1-1X31)
	(E-11X31) -- (F-1-1X31)
	(A1-1X31) -- (E11X31)
	(B-11X31) -- (D-11X31)
	(E1-1X31) -- (F11X31)

	(A-11X41) -- (A-1-1X41)	
	(D1-1X41) -- (C1-1X41)
	(E-11X41) -- (F-1-1X41)
	(A1-1X41) -- (A11X41)
	(B-11X41) -- (D-11X41)
	(E1-1X41) -- (E11X41)
	
	(A-11X51) -- (A-1-1X51)	
	(B1-1X51) -- (B-1-1X51)
	(E-11X51) -- (F-1-1X51)
	(D11X51) -- (E11X51)
	(B-11X51) -- (D-11X51)
	(E1-1X51) -- (D1-1X51)

	(A-11X02) -- (A-1-1X02)	
	(B1-1X02) -- (A1-1X02)
	(E-11X02) -- (F-1-1X02)
	(D11X02) -- (E11X02)
	(B-11X02) -- (D-11X02)
	(C1-1X02) -- (D1-1X02)

	(A-11X12) -- (A-1-1X12)	
	(B1-1X12) -- (B-1-1X12)
	(E-11X12) -- (E-1-1X12)
	(D11X12) -- (E11X12)
	(C-11X12) -- (D-11X12)
	(E1-1X12) -- (D1-1X12)

	(A-11X22) -- (A-1-1X22)	
	(B1-1X22) -- (B-1-1X22)
	(E-11X22) -- (E-1-1X22)
	(A1-1X22) -- (E11X22)
	(C-11X22) -- (D-11X22)
	(E1-1X22) -- (F11X22)

	(A-11X32) -- (A-1-1X32)	
	(D1-1X32) -- (C1-1X32)
	(E-11X32) -- (E-1-1X32)
	(A1-1X32) -- (A11X32)
	(C-11X32) -- (D-11X32)
	(E1-1X32) -- (E11X32)

	(A-11X42) -- (A-1-1X42)	
	(D1-1X42) -- (C1-1X42)
	(E-11X42) -- (E-1-1X42)
	(A1-1X42) -- (B1-1X42)
	(C-11X42) -- (D-11X42)
	(D11X42) -- (E11X42)

	(B11X52) -- (A11X52)	
	(D1-1X52) -- (C1-1X52)
	(D-1-1X52) -- (E-1-1X52)
	(A1-1X52) -- (B1-1X52)
	(E-11X52) -- (D-11X52)
	(D11X52) -- (C11X52)

	(A11X03) -- (B11X03)	
	(B1-1X03) -- (B-1-1X03)
	(E-11X03) -- (D-11X03)
	(D11X03) -- (C11X03)
	(E-1-1X03) -- (D-1-1X03)
	(E1-1X03) -- (D1-1X03)

	(A11X13) -- (B11X13)	
	(B-1-1X13) -- (A-1-1X13)
	(E-11X13) -- (D-11X13)
	(D11X13) -- (C11X13)
	(C-1-1X13) -- (D-1-1X13)
	(E1-1X13) -- (D1-1X13)

	(A11X23) -- (B11X23)	
	(B1-1X23) -- (A1-1X23)
	(F-1-1X23) -- (E-11X23)
	(D11X23) -- (C11X23)
	(A-11X23) -- (E-1-1X23)
	(C1-1X23) -- (D1-1X23)

	(A11X33) -- (B11X33)	
	(B1-1X33) -- (B-1-1X33)
	(E-11X33) -- (F-1-1X33)
	(D11X33) -- (C11X33)
	(E-1-1X33) -- (A-11X33)
	(E1-1X33) -- (D1-1X33)

	(A11X43) -- (B11X43)	
	(B-1-1X43) -- (A-1-1X43)
	(F-1-1X43) -- (E-11X43)
	(C1-1X43) -- (D1-1X43)
	(B-11X43) -- (D-11X43)
	(E1-1X43) -- (E11X43)

	(A11X53) -- (B11X53)	
	(B-1-1X53) -- (A-1-1X53)
	(F-1-1X53) -- (E-11X53)
	(B1-1X53) -- (D1-1X53)
	(B-11X53) -- (D-11X53)
	(E1-1X53) -- (F11X53)

	(A11X04) -- (B11X04)	
	(B-1-1X04) -- (D-1-1X04)
	(F-1-1X04) -- (E-11X04)
	(B1-1X04) -- (A1-1X04)
	(B-11X04) -- (D-11X04)
	(E1-1X04) -- (F11X04)

	(A11X14) -- (B11X14)	
	(C-1-1X14) -- (D-1-1X14)
	(E-1-1X14) -- (A-11X14)
	(B1-1X14) -- (A1-1X14)
	(C-11X14) -- (D-11X14)
	(E1-1X14) -- (F11X14)

	(A11X24) -- (B11X24)	
	(B-1-1X24) -- (A-1-1X24)
	(E-1-1X24) -- (E-11X24)
	(C1-1X24) -- (D1-1X24)
	(C-11X24) -- (D-11X24)
	(E1-1X24) -- (E11X24)

	(A11X34) -- (B11X34)	
	(B-1-1X34) -- (D-1-1X34)
	(E-1-1X34) -- (E-11X34)
	(B1-1X34) -- (A1-1X34)
	(C-11X34) -- (D-11X34)
	(E1-1X34) -- (F11X34)

	(D11X44) -- (B11X44)	
	(B-1-1X44) -- (D-1-1X44)
	(F-1-1X44) -- (E-11X44)
	(B1-1X44) -- (A1-1X44)
	(A-11X44) -- (B-11X44)
	(E1-1X44) -- (F11X44)

	(D11X54) -- (E11X54)	
	(E-1-1X54) -- (D-1-1X54)
	(C-11X54) -- (D-11X54)
	(B1-1X54) -- (A1-1X54)
	(A-11X54) -- (B-11X54)
	(C1-1X54) -- (D1-1X54)

	(A1-1X05) -- (E11X05)	
	(E-1-1X05) -- (D-1-1X05)
	(F11X05) -- (E1-1X05)
	(B1-1X05) -- (B-1-1X05)
	(A-11X05) -- (B-11X05)
	(C-11X05) -- (D-11X05)

	(D11X15) -- (E11X15)	
	(E-1-1X15) -- (D-1-1X15)
	(D1-1X15) -- (E1-1X15)
	(B1-1X15) -- (B-1-1X15)
	(A-11X15) -- (B-11X15)
	(C-11X15) -- (D-11X15)

	(D11X25) -- (E11X25)	
	(E-1-1X25) -- (D-1-1X25)
	(D1-1X25) -- (E1-1X25)
	(B1-1X25) -- (B-1-1X25)
	(B11X25) -- (B-11X25)
	(E-11X25) -- (D-11X25)

	(A1-1X35) -- (E11X35)	
	(E-1-1X35) -- (D-1-1X35)
	(F11X35) -- (E1-1X35)
	(B1-1X35) -- (B-1-1X35)
	(B11X35) -- (B-11X35)
	(E-11X35) -- (D-11X35)
		
	(A1-1X45) -- (E11X45)	
	(E-1-1X45) -- (A-11X45)
	(F11X45) -- (E1-1X45)
	(B1-1X45) -- (B-1-1X45)
	(B11X45) -- (B-11X45)
	(E-11X45) -- (F-1-1X45);

\foreach \a in {0,1,...,5}
{
\draw[line width=1.8]
	(F-11X\a6) -- (G-11X\a6)
	(F-1-1X\a6) -- (G-1-1X\a6)
	(F11X\a6) -- (G11X\a6)
	(F1-1X\a6) -- (G1-1X\a6);
}

\draw[line width=1.8]
	(F-11X55) -- (G-11X55)
	(F-1-1X55) -- (G-1-1X55)
	(F11X55) -- (G11X55)
	(F1-1X55) -- (G1-1X55);

\draw[line width=1.8]	
	(B11X55) -- (B-11X55)	
	(E-1-1X55) -- (A-11X55)
	(D1-1X55) -- (B1-1X55)	
	(A11X55) -- (A1-1X55)

	(B11X06) -- (B-11X06)	
	(E-1-1X06) -- (A-11X06)
	(B-1-1X06) -- (B1-1X06)	
	(E11X06) -- (A1-1X06)

	(B-11X16) -- (D-11X16)	
	(A-1-1X16) -- (B-1-1X16)	
	(D1-1X16) -- (B1-1X16)	
	(A11X16) -- (B11X16)

	(B-11X26) -- (D-11X26)	
	(A-1-1X26) -- (A-11X26)	
	(D1-1X26) -- (B1-1X26)	
	(A11X26) -- (A1-1X26)
	
	(B-11X36) -- (D-11X36)	
	(B-1-1X36) -- (D-1-1X36)	
	(A1-1X36) -- (B1-1X36)	
	(A11X36) -- (B11X36)

	(B-11X46) -- (D-11X46)	
	(B-1-1X46) -- (D-1-1X46)	
	(A1-1X46) -- (E11X46)	
	(A11X46) -- (A-1-1X46)

	(B-11X56) -- (A-11X56)	
	(B-1-1X56) -- (D-1-1X56)	
	(A1-1X56) -- (B1-1X56)	
	(D11X56) -- (B11X56);
	
\foreach \a in {0,...,5}
\foreach \b in {7,...,10}
{
\draw[line width=1.8]
	(F-11X\a\b) -- (G-11X\a\b)
	(F-1-1X\a\b) -- (G-1-1X\a\b);
}

\draw[line width=1.8]
	(F-11X011) -- (G-11X011)
	(F-1-1X011) -- (G-1-1X011);

\draw[line width=1.8]
	(B11X07) -- (D11X07)	
	(E11X07) -- (E1-1X07)
	(D1-1X07) -- (C1-1X07)
	(A-1-1X07) -- (A11X07)
	(E-1-1X07) -- (A-11X07)

	(B11X17) -- (D11X17)	
	(F11X17) -- (E1-1X17)
	(D1-1X17) -- (B1-1X17)
	(A-1-1X17) -- (A11X17)
	(E-1-1X17) -- (A-11X17)

	(C11X27) -- (D11X27)	
	(E11X27) -- (A1-1X27)
	(D1-1X27) -- (C1-1X27)
	(A-1-1X27) -- (A11X27)
	(E-1-1X27) -- (A-11X27)

	(B-11X37) -- (D-11X37)	
	(B-1-1X37) -- (D-1-1X37)	
	(D11X37) -- (E11X37)	
	(A11X37) -- (A-1-1X37)
	(D1-1X37) -- (E1-1X37)

	(B-11X47) -- (D-11X47)	
	(B-1-1X47) -- (D-1-1X47)	
	(A1-1X47) -- (E11X47)
	(A11X47) -- (A-1-1X47)
	(F11X47) -- (E1-1X47)

	(B-11X57) -- (D-11X57)	
	(C1-1X57) -- (D1-1X57)	
	(A1-1X57) -- (A11X57)
	(A-11X57) -- (A-1-1X57)
	(E11X57) -- (E1-1X57)
	
	(B-11X08) -- (D-11X08)	
	(E1-1X08) -- (D1-1X08)
	(B1-1X08) -- (B-1-1X08)
	(A-11X08) -- (A-1-1X08)
	(E11X08) -- (D11X08)	

	(B-11X18) -- (D-11X18)	
	(C1-1X18) -- (D1-1X18)
	(B1-1X18) -- (A1-1X18)
	(A-11X18) -- (A-1-1X18)
	(E11X18) -- (D11X18)

	(B-11X28) -- (D-11X28)	
	(B1-1X28) -- (D1-1X28)
	(A11X28) -- (A1-1X28)
	(A-11X28) -- (A-1-1X28)
	(E1-1X28) -- (F11X28)

	(B-11X38) -- (D-11X38)	
	(B1-1X38) -- (B-1-1X38)
	(E11X38) -- (A1-1X38)
	(A-11X38) -- (A-1-1X38)
	(E1-1X38) -- (F11X38)
	
	(B11X48) -- (A11X48)	
	(B1-1X48) -- (B-1-1X48)
	(C11X48) -- (D11X48)
	(A-11X48) -- (E-1-1X48)
	(E1-1X48) -- (D1-1X48)

	(B11X58) -- (A11X58)	
	(B1-1X58) -- (A1-1X58)
	(C11X58) -- (D11X58)
	(A-11X58) -- (E-1-1X58)
	(C1-1X58) -- (D1-1X58)

	(B-11X09) -- (D-11X09)	
	(C1-1X09) -- (D1-1X09)	
	(B11X09) -- (A11X09)
	(B-1-1X09) -- (A-1-1X09)
	(E11X09) -- (E1-1X09)

	(B-11X19) -- (D-11X19)	
	(B1-1X19) -- (A1-1X19)
	(B11X19) -- (A11X19)
	(B-1-1X19) -- (D-1-1X19)
	(F11X19) -- (E1-1X19)

	(B-11X29) -- (D-11X29)	
	(B1-1X29) -- (D1-1X29)
	(B11X29) -- (A11X29)
	(B-1-1X29) -- (A-1-1X29)
	(F11X29) -- (E1-1X29)

	(B-11X39) -- (A-11X39)	
	(B1-1X39) -- (D1-1X39)
	(B11X39) -- (D11X39)
	(B-1-1X39) -- (A-1-1X39)
	(F11X39) -- (E1-1X39)

	(B-11X49) -- (A-11X49)	
	(C1-1X49) -- (D1-1X49)
	(B11X49) -- (D11X49)
	(B-1-1X49) -- (A-1-1X49)
	(E11X49) -- (E1-1X49)

	(B-11X59) -- (A-11X59)	
	(C1-1X59) -- (D1-1X59)
	(C11X59) -- (D11X59)
	(B-1-1X59) -- (A-1-1X59)
	(E11X59) -- (A1-1X59)

	(B-11X010) -- (A-11X010)	
	(E1-1X010) -- (F11X010)
	(B11X010) -- (D11X010)
	(B-1-1X010) -- (D-1-1X010)
	(B1-1X010) -- (A1-1X010)

	(B-11X110) -- (A-11X110)	
	(E1-1X110) -- (D1-1X110)
	(C11X110) -- (D11X110)
	(B-1-1X110) -- (D-1-1X110)
	(A11X110) -- (A1-1X110)

	(B-11X210) -- (B11X210)	
	(E1-1X210) -- (D1-1X210)
	(E11X210) -- (D11X210)
	(B-1-1X210) -- (B1-1X210)
	(A-11X210) -- (E-1-1X210)

	(B-11X310) -- (B11X310)	
	(E1-1X310) -- (F11X310)
	(E11X310) -- (A1-1X310)
	(B-1-1X310) -- (B1-1X310)
	(A-11X310) -- (E-1-1X310)

	(B-11X410) -- (B11X410)	
	(E1-1X410) -- (F11X410)
	(A11X410) -- (A1-1X410)
	(D1-1X410) -- (B1-1X410)
	(A-11X410) -- (E-1-1X410)

	(B-11X510) -- (B11X510)	
	(E1-1X510) -- (E11X510)
	(A11X510) -- (A1-1X510)
	(D1-1X510) -- (C1-1X510)
	(A-11X510) -- (E-1-1X510)

	(B-11X011) -- (B11X011)	
	(D11X011) -- (E11X011)
	(B1-1X011) -- (A1-1X011)
	(D1-1X011) -- (C1-1X011)
	(A-11X011) -- (E-1-1X011);
	
\end{tikzpicture}
\caption{Even family of timezone tiling, distance $3$, edge combination $a^4b$.}
\label{timezone4ab_d3even}
\end{figure}

Up to the symmetry of $180^{\circ}$ rotation, the complete list for the even family of timezone tilings is given by Figure \ref{timezone4ab_d3even}. If the rotation preserves the meridian boundaries of a tiling and changes the interior, then we indicate $\times 2$. If the rotation changes the meridian boundaries, or preserves the whole timezone tiling (including the boundaries as well as the interior), then we indicate $\times 1$.
\begin{itemize}
\item Nos.1-35 have the meridian $aaa$ as both boundaries. There are total of $60$ timezone tilings.
\item Nos.36-42 have the meridian $bab$ as both boundaries. There are total of $10$ timezone tilings.
\item Nos.43-67 have the meridians $bab$ and $aaa$ as the left and right boundaries. There are total of $25$ timezone tilings. Their $180^{\circ}$ rotations give $25$ timezone tilings with $aaa$ and $bab$ as boundaries.
\end{itemize}

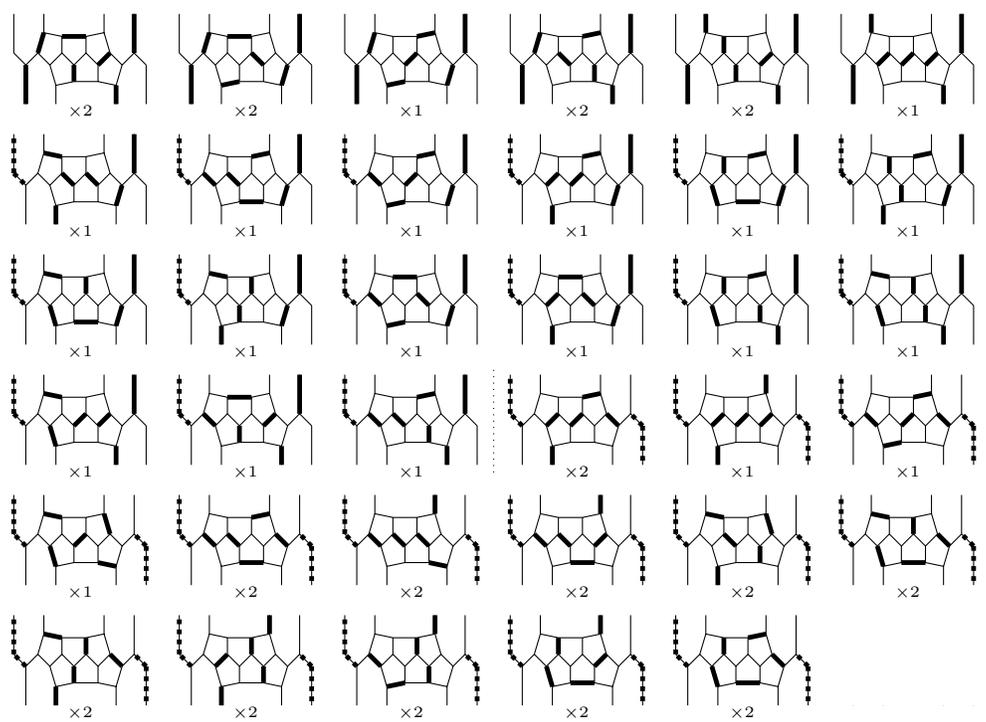
\begin{figure}[htp]
\centering
\begin{tikzpicture}[>=latex]

\foreach \x in {-1,1}
\foreach \y in {-1,1}
\foreach \a in {0,...,5}
\foreach \b in {0,...,5}
{
\begin{scope}[shift={(2.2*\a cm, -1.6*\b cm)}]

\coordinate  (A\x\y X\a\b) at (0.16*\x-0.08*\y,0.08*\y);
\coordinate  (B\x\y X\a\b) at (0.16*\x-0.08*\y,0.3*\y);
\coordinate  (C\x\y X\a\b) at (0.40*\x-0.08*\y,0.6*\y);
\coordinate  (D\x\y X\a\b) at (0.40*\x-0.08*\y,0.35*\y);
\coordinate  (E\x\y X\a\b) at (0.48*\x-0.08*\y,0.08*\y);
\coordinate  (F\x\y X\a\b) at (0.80*\x-0.08*\y,0.08*\y);
\coordinate  (G\x\y X\a\b) at (0.80*\x-0.08*\y,0.6*\y);

\coordinate  (P\a\b) at (0,-0.7);

\draw
	(A\x\y X\a\b) -- (B\x\y X\a\b)
	(C\x\y X\a\b) -- (D\x\y X\a\b) -- (B\x\y X\a\b)
	(D\x\y X\a\b) -- (E\x\y X\a\b)
	(F\x\y X\a\b) -- (G\x\y X\a\b);
	
\end{scope}
}

\foreach \a in {0,...,5}
\foreach \b in {0,...,5}
\draw
	(F1-1X\a\b) -- (F11X\a\b) -- (E1-1X\a\b) -- (E11X\a\b) -- (A1-1X\a\b) -- (A11X\a\b) -- (A-1-1X\a\b) -- (A-11X\a\b) -- (E-1-1X\a\b) -- (E-11X\a\b) -- (F-1-1X\a\b) -- (F-11X\a\b)
	(B11X\a\b) -- (B-11X\a\b)
	(B1-1X\a\b) -- (B-1-1X\a\b);


\foreach \u / \v in {3/3,6/0}
\draw[dotted, shift={(2.2*\u cm,-1.6*\v cm)}]
	(-1.1,-0.7) -- (-1.1,0.7);


\node at (P00) {\tiny $\times 2$};
\node at (P10) {\tiny $\times 2$};
\node at (P20) {\tiny $\times 1$};
\node at (P30) {\tiny $\times 2$};
\node at (P40) {\tiny $\times 2$};
\node at (P50) {\tiny $\times 1$};

\node at (P01) {\tiny $\times 1$};
\node at (P11) {\tiny $\times 1$};
\node at (P21) {\tiny $\times 1$};
\node at (P31) {\tiny $\times 1$};
\node at (P41) {\tiny $\times 1$};
\node at (P51) {\tiny $\times 1$};
\node at (P02) {\tiny $\times 1$};
\node at (P12) {\tiny $\times 1$};
\node at (P22) {\tiny $\times 1$};
\node at (P32) {\tiny $\times 1$};
\node at (P42) {\tiny $\times 1$};
\node at (P52) {\tiny $\times 1$};
\node at (P03) {\tiny $\times 1$};
\node at (P13) {\tiny $\times 1$};
\node at (P23) {\tiny $\times 1$};

\node at (P33) {\tiny $\times 2$};
\node at (P43) {\tiny $\times 1$};
\node at (P53) {\tiny $\times 1$};
\node at (P04) {\tiny $\times 1$};
\node at (P14) {\tiny $\times 2$};
\node at (P24) {\tiny $\times 2$};
\node at (P34) {\tiny $\times 2$};
\node at (P44) {\tiny $\times 2$};
\node at (P54) {\tiny $\times 2$};
\node at (P05) {\tiny $\times 2$};
\node at (P15) {\tiny $\times 2$};
\node at (P25) {\tiny $\times 2$};
\node at (P35) {\tiny $\times 2$};
\node at (P45) {\tiny $\times 2$};

\foreach \a in {0,...,5}
{
\draw[line width=1.8]
	(F-1-1X\a0) -- (G-1-1X\a0)
	(F11X\a0) -- (G11X\a0);
}

\draw[line width=1.8]
	(B-11X00) -- (B11X00)	
	(E11X00) -- (A1-1X00)
	(C1-1X00) -- (D1-1X00)	
	(B-1-1X00) -- (A-1-1X00)
	(E-11X00) -- (D-11X00)

	(B-11X10) -- (B11X10)	
	(A11X10) -- (A1-1X10)
	(E1-1X10) -- (D1-1X10)	
	(B-1-1X10) -- (D-1-1X10)
	(E-11X10) -- (D-11X10)

	(D11X20) -- (B11X20)	
	(A11X20) -- (A-1-1X20)
	(E1-1X20) -- (D1-1X20)	
	(B-1-1X20) -- (D-1-1X20)
	(E-11X20) -- (D-11X20)

	(D11X30) -- (B11X30)	
	(A-11X30) -- (A-1-1X30)
	(C1-1X30) -- (D1-1X30)	
	(B1-1X30) -- (A1-1X30)
	(E-11X30) -- (D-11X30)

	(A-1-1X40) -- (B-1-1X40)	
	(A-11X40) -- (B-11X40)
	(C1-1X40) -- (D1-1X40)	
	(E11X40) -- (A1-1X40)
	(C-11X40) -- (D-11X40)

	(A-1-1X50) -- (A11X50)	
	(A-11X50) -- (E-1-1X50)
	(C1-1X50) -- (D1-1X50)	
	(E11X50) -- (A1-1X50)
	(C-11X50) -- (D-11X50);
	

\foreach \a in {0,1,2,3,4,5}
\foreach \b in {1,2}
{
\draw[line width=1.8,dotted]
	(F-11X\a\b) -- (G-11X\a\b)
	(F-11X\a\b) -- (F-1-1X\a\b);
}

\foreach \a in {0,1,2}
{
\draw[line width=1.8,dotted]
	(F-11X\a3) -- (G-11X\a3)
	(F-11X\a3) -- (F-1-1X\a3);
}

\draw[line width=1.8]
	(F11X01) -- (G11X01)	
	(B-11X01) -- (D-11X01)
	(C-1-1X01) -- (D-1-1X01)
	(A-11X01) -- (A-1-1X01)
	(A11X01) -- (A1-1X01)	
	(D1-1X01) -- (E1-1X01)	

	(F11X11) -- (G11X11)	
	(A-11X11) -- (A-1-1X11)
	(D1-1X11) -- (E1-1X11)
	(E-11X11) -- (E-1-1X11)
	(B1-1X11) -- (B-1-1X11)
	(B11X11) -- (D11X11)
	
	(F11X21) -- (G11X21)	
	(B11X21) -- (D11X21)	
	(A11X21) -- (A-1-1X21)
	(D1-1X21) -- (E1-1X21)		
	(E-11X21) -- (E-1-1X21)
	(B-1-1X21) -- (D-1-1X21)

	(F11X31) -- (G11X31)	
	(B11X31) -- (D11X31)	
	(A11X31) -- (A-1-1X31)
	(D1-1X31) -- (E1-1X31)		
	(A-11X31) -- (E-1-1X31)
	(C-1-1X31) -- (D-1-1X31)
	
	(F11X41) -- (G11X41)	
	(B11X41) -- (D11X41)	
	(B1-1X41) -- (B-1-1X41)
	(D1-1X41) -- (E1-1X41)		
	(A-11X41) -- (B-11X41)
	(E-1-1X41) -- (D-1-1X41)

	(F11X51) -- (G11X51)	
	(B11X51) -- (D11X51)	
	(A-1-1X51) -- (B-1-1X51)
	(D1-1X51) -- (E1-1X51)		
	(A-11X51) -- (B-11X51)
	(C-1-1X51) -- (D-1-1X51)

	(F11X02) -- (G11X02)	
	(B-11X02) -- (D-11X02)	
	(B1-1X02) -- (B-1-1X02)	
	(D1-1X02) -- (E1-1X02)	
	(E-1-1X02) -- (D-1-1X02)
	(A11X02) -- (B11X02)

	(F11X12) -- (G11X12)	
	(B-11X12) -- (D-11X12)	
	(A-1-1X12) -- (B-1-1X12)	
	(D1-1X12) -- (E1-1X12)	
	(C-1-1X12) -- (D-1-1X12)
	(A11X12) -- (B11X12)
	
	(F11X22) -- (G11X22)	
	(B-11X22) -- (B11X22)	
	(D-1-1X22) -- (B-1-1X22)	
	(D1-1X22) -- (E1-1X22)	
	(A11X22) -- (A1-1X22)
	(E-11X22) -- (E-1-1X22)

	(F11X32) -- (G11X32)	
	(B-11X32) -- (B11X32)	
	(D-1-1X32) -- (C-1-1X32)	
	(D1-1X32) -- (E1-1X32)	
	(A11X32) -- (A1-1X32)
	(A-11X32) -- (E-1-1X32)

	(F11X42) -- (G11X42)	
	(B11X42) -- (D11X42)	
	(B1-1X42) -- (A1-1X42)
	(C1-1X42) -- (D1-1X42)	
	(A-11X42) -- (B-11X42)
	(D-1-1X42) -- (E-1-1X42)
	
	(F11X52) -- (G11X52)	
	(B-11X52) -- (D-11X52)	
	(B1-1X52) -- (A1-1X52)
	(C1-1X52) -- (D1-1X52)	
	(A11X52) -- (B11X52)
	(D-1-1X52) -- (E-1-1X52)

	(F11X03) -- (G11X03)	
	(B-11X03) -- (D-11X03)	
	(E11X03) -- (A1-1X03)
	(C1-1X03) -- (D1-1X03)	
	(A11X03) -- (A-1-1X03)
	(D-1-1X03) -- (E-1-1X03)

	(F11X13) -- (G11X13)	
	(B-11X13) -- (B11X13)	
	(E11X13) -- (A1-1X13)
	(C1-1X13) -- (D1-1X13)	
	(B-1-1X13) -- (A-1-1X13)
	(E-11X13) -- (E-1-1X13)

	(F11X23) -- (G11X23)	
	(D11X23) -- (B11X23)	
	(B1-1X23) -- (A1-1X23)
	(C1-1X23) -- (D1-1X23)	
	(A-11X23) -- (A-1-1X23)
	(E-11X23) -- (E-1-1X23);
	

\foreach \a in {0,1,2,3,4,5}
\foreach \b in {4,5}
{
\draw[line width=1.8,dotted]
	(F-1-1X\a\b) -- (F-11X\a\b)
	(F-11X\a\b) -- (G-11X\a\b)
	(F1-1X\a\b) -- (F11X\a\b)
	(F1-1X\a\b) -- (G1-1X\a\b);
}

\foreach \a in {3,4,5}
{
\draw[line width=1.8,dotted]
	(F-1-1X\a3) -- (F-11X\a3)
	(F-11X\a3) -- (G-11X\a3)
	(F1-1X\a3) -- (F11X\a3)
	(F1-1X\a3) -- (G1-1X\a3);
}

\draw[line width=1.8]
	(B11X33) -- (D11X33)	
	(A11X33) -- (A-1-1X33)
	(E11X33) -- (E1-1X33)		
	(A-11X33) -- (E-1-1X33)
	(C-1-1X33) -- (D-1-1X33)

	(C11X43) -- (D11X43)	
	(A11X43) -- (A-1-1X43)
	(E11X43) -- (A1-1X43)		
	(A-11X43) -- (E-1-1X43)
	(C-1-1X43) -- (D-1-1X43)
	
	(B11X53) -- (D11X53)	
	(A11X53) -- (A-1-1X53)
	(E11X53) -- (E1-1X53)		
	(E-11X53) -- (E-1-1X53)
	(B-1-1X53) -- (D-1-1X53)

	(B-11X04) -- (D-11X04)	
	(A11X04) -- (A-1-1X04)
	(E11X04) -- (D11X04)		
	(D-1-1X04) -- (E-1-1X04)
	(B1-1X04) -- (D1-1X04)

	(A-11X14) -- (A-1-1X14)	
	(E11X14) -- (E1-1X14)
	(E-11X14) -- (E-1-1X14)
	(B1-1X14) -- (B-1-1X14)
	(B11X14) -- (D11X14)

	(A-11X24) -- (A-1-1X24)	
	(B1-1X24) -- (D1-1X24)
	(E-11X24) -- (E-1-1X24)
	(A1-1X24) -- (A11X24)
	(C11X24) -- (D11X24)

	(A-11X34) -- (A-1-1X34)	
	(B1-1X34) -- (B-1-1X34)
	(E-11X34) -- (E-1-1X34)
	(A1-1X34) -- (E11X34)
	(C11X34) -- (D11X34)

	(B-11X44) -- (D-11X44)	
	(A-11X44) -- (A-1-1X44)
	(E11X44) -- (D11X44)		
	(D-1-1X44) -- (C-1-1X44)
	(B1-1X44) -- (A1-1X44)

	(B-11X54) -- (D-11X54)	
	(A11X54) -- (B11X54)
	(E11X54) -- (E1-1X54)		
	(D-1-1X54) -- (E-1-1X54)
	(B1-1X54) -- (B-1-1X54)

	(B-11X05) -- (D-11X05)	
	(A11X05) -- (B11X05)
	(E11X05) -- (E1-1X05)		
	(D-1-1X05) -- (C-1-1X05)
	(A-1-1X05) -- (B-1-1X05)

	(C11X15) -- (D11X15)	
	(A11X15) -- (B11X15)
	(B1-1X15) -- (A1-1X15)		
	(A-11X15) -- (E-1-1X15)
	(C-1-1X15) -- (D-1-1X15)

	(C11X25) -- (D11X25)	
	(A11X25) -- (B11X25)
	(B1-1X25) -- (A1-1X25)		
	(E-11X25) -- (E-1-1X25)
	(B-1-1X25) -- (D-1-1X25)

	(A-11X35) -- (B-11X35)	
	(B1-1X35) -- (B-1-1X35)
	(D-1-1X35) -- (E-1-1X35)
	(A1-1X35) -- (E11X35)
	(C11X35) -- (D11X35)

	(A-11X45) -- (B-11X45)	
	(B1-1X45) -- (B-1-1X45)
	(D-1-1X45) -- (E-1-1X45)
	(E1-1X45) -- (E11X45)
	(B11X45) -- (D11X45);

\fill[white,shift={(11 cm, -8 cm)}] (-1.2,-0.6) rectangle (1.2,0.6);
	
\end{tikzpicture}
\caption{Odd family of timezone tiling, distance $3$, edge combination $a^4b$.}
\label{timezone4ab_d3odd}
\end{figure}

Similarly, up to symmetry, the complete list for the odd family of timezone tilings is given by Figure \ref{timezone4ab_d3odd}. To simplify the drawings, we introduce dotted thick adjacent edges, which means that any one of the two edges can be the thick edge (i.e., the $b$-edge). For example, the left meridian of the $7$-th tiling can be either $baa$ or $aba$. 
\begin{itemize}
\item Nos.1-6 have $aab$ and $baa$ the left and right boundaries. There are total of $10$ timezone tilings.
\item Nos.7-21, with $baa$ as both boundaries, or with $aba$ and $baa$ as the left and right boundaries. The $180^{\circ}$ rotation changes the boundaries to $aab$ as both boundaries, or to $aab$ and $aba$ as the left and right boundaries. In each of the four cases, there are total of $15$ timezone tilings.
\item Nos.22-35, with choices of boundary meridians to be preserved by the $180^{\circ}$ rotation. The boundary meridians can be both $aba$, or have $baa$ and $aab$ as the left and right boundaries. In either case, there are total of $25$ timezone tilings. 
\item Nos.22-35, with choices of boundary meridians to be changed by the $180^{\circ}$ rotation. The boundary meridians can have $baa$ and $aba$ as the left and right boundaries, or have $aba$ and $baa$ as the left and right boundaries. The $14$ tilings with two possible choices of boundary meridians and their $180^{\circ}$ rotations give total of $50$ timezone tilings. The number $50$ is the total number for the two cases. Since the $180^{\circ}$ rotation switches between the two cases, and the two cases are disjoint, we get total of $25$ timezone tilings for each of the two cases. 
\end{itemize}

\begin{figure}[htp]
\centering
\begin{tikzpicture}[>=latex]


\node[rectangle,draw,inner sep=2] (A) at (-1.2,0) {\small $aaa$};
\node[rectangle,draw,inner sep=1] (C) at (1.2,0) {\small $bab$};	

\draw[->] (A) to[out=100,in=90] (-2.5,0) to[out=270,in=260] (A);
\node[draw,shape=circle,fill=white,inner sep=0] at (-2.5,0) {\small $60$};	

\draw[->] (C) to[out=280,in=270] (2.5,0) to[out=90,in=80] (C);
\node[draw,shape=circle,fill=white,inner sep=0] at (2.5,0) {\small $10$};	
	
\draw (A) edge[<->] node[draw,shape=circle,fill=white,inner sep=0] {\small $25$} (C);


\begin{scope}[xshift=7cm]

\node[rectangle,draw,inner sep=1] (A) at (90:1.7) {\small $aba$};	
\node[rectangle,draw,inner sep=1] (B) at (200:1.7) {\small $aab$};
\node[rectangle,draw,inner sep=1] (C) at (-20:1.7) {\small $baa$};

\draw[->] (A) to[out=0,in=0] (90:2.6) to[out=180,in=180] (A);	
\node[draw,shape=circle,fill=white,inner sep=0] at (90:2.6) {\small 25};

\draw[->] (B) to[out=100,in=100] (200:2.6) to[out=290,in=290] (B);	
\node[draw,shape=circle,fill=white,inner sep=0] at (200:2.6) {\small 15};

\draw[->] (C) to[out=250,in=250] (-20:2.6) to[out=70,in=70] (C);	
\node[draw,shape=circle,fill=white,inner sep=0] at (-20:2.6) {\small 15};

\draw[->] (A) to[out=215,in=75] (B);
\draw[->] (B) to[out=45,in=245] (A);
\node[draw,shape=circle,fill=white,inner sep=0] at (145:1.35) {\small 25};
\node[draw,shape=circle,fill=white,inner sep=0] at (145:0.75) {\small 15};

\draw[<-] (A) to[out=-35,in=105] (C);
\draw[<-] (C) to[out=135,in=295] (A);
\node[draw,shape=circle,fill=white,inner sep=0] at (35:1.35) {\small 25};
\node[draw,shape=circle,fill=white,inner sep=0] at (35:0.75) {\small 15};

\draw[->] (B) to[out=15,in=165] (C);
\draw[->] (C) to[out=195,in=-15] (B);
\node[draw,shape=circle,fill=white,inner sep=0] at (270:0.9) {\small 25};
\node[draw,shape=circle,fill=white,inner sep=0] at (270:0.3) {\small 10};
	
\end{scope}

\end{tikzpicture}
\caption{Earth map tiling, distance $3$, edge combination $a^4b$.}
\label{4ab_d3}
\end{figure}
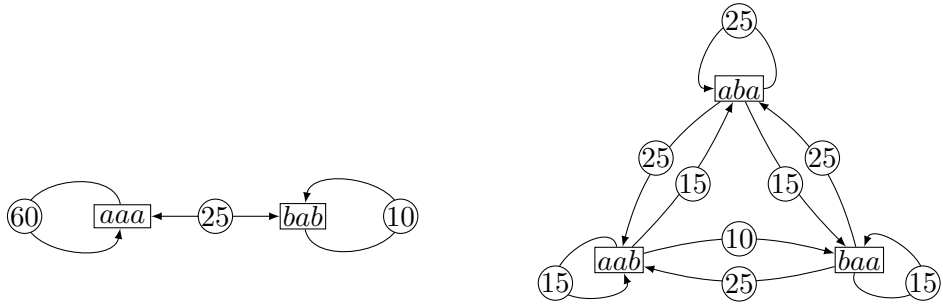

All the timezone tilings for the earth map tiling of distance $3$ are represented by the arrows in Figure \ref{4ab_d3}. The graph on the left consists of the even family, and the graph on the right consists of the odd family. An edge congruent earth map tiling is a directed cycle in the graph.

\subsection*{Distance 2}

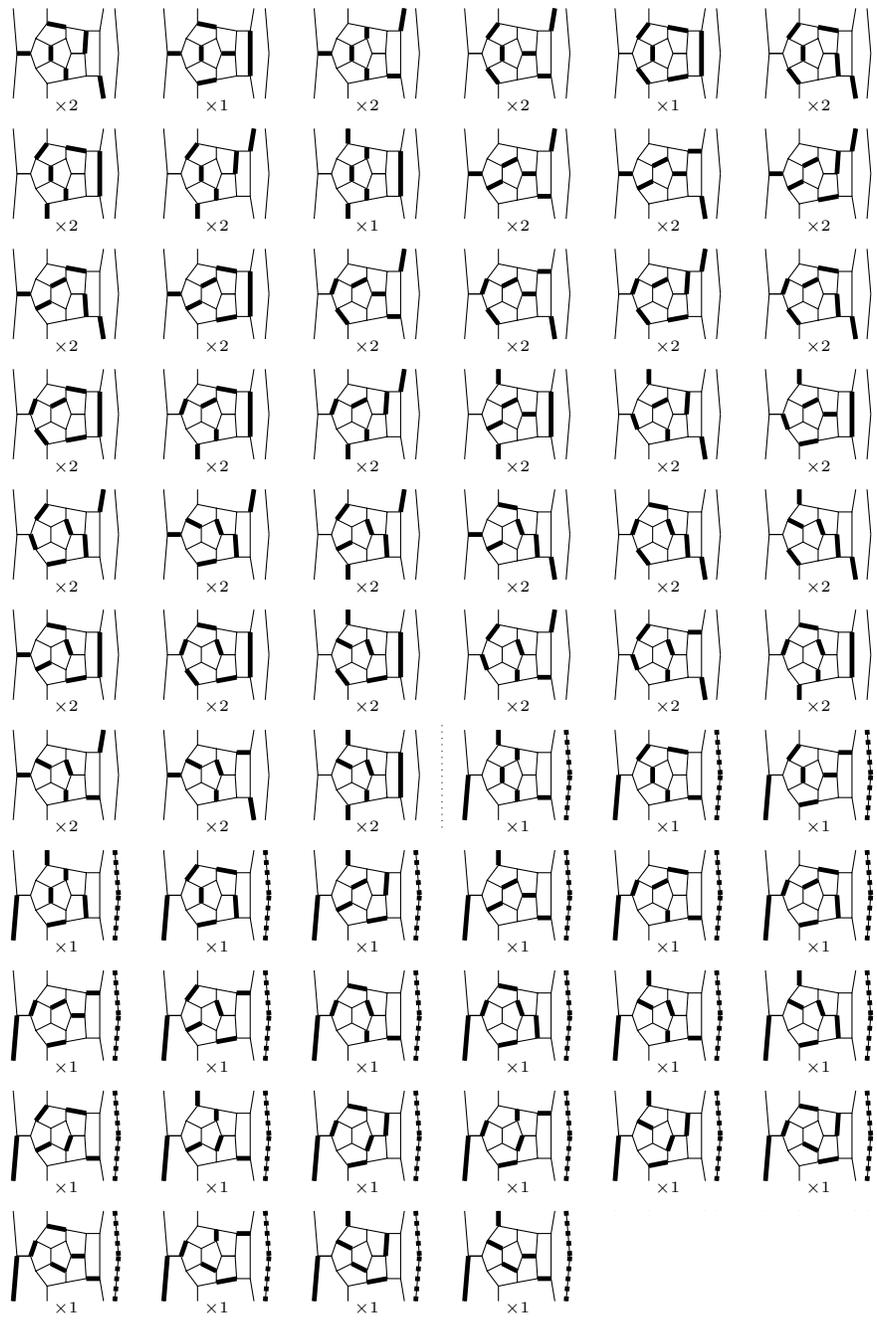
\begin{figure}[htp]
\centering
\begin{tikzpicture}[>=latex]

\foreach \y in {-1,1}
\foreach \a in {0,...,5}
\foreach \b in {0,...,10}
{
\begin{scope}[shift={(2*\a cm, -1.6*\b cm)}]

\coordinate  (A\y X\a\b) at (0,0.1*\y);
\coordinate  (B\y X\a\b) at (-0.2,0.2*\y);
\coordinate  (C\y X\a\b) at (0.2,0.2*\y);
\coordinate  (D\y X\a\b) at (-0.05,0.4*\y);
\coordinate  (E\y X\a\b) at (0.2,0.35*\y);
\coordinate  (F\y X\a\b) at (-0.05,0.6*\y);
\coordinate  (GX\a\b) at (0.27,0);
\coordinate  (HX\a\b) at (0.45,0);
\coordinate  (I\y X\a\b) at (0.47,0.3*\y);
\coordinate  (J\y X\a\b) at (0.65,0.3*\y);
\coordinate  (K\y X\a\b) at (0.7,0.6*\y);
\coordinate  (L\y X\a\b) at (0.85,0.6*\y);
\coordinate  (MX\a\b) at (0.9,0);
\coordinate  (NX\a\b) at (-0.27,0);
\coordinate  (OX\a\b) at (-0.45,0);
\coordinate  (P\y X\a\b) at (-0.5,0.6*\y);

\coordinate  (Q\a\b) at (0.2,-0.7);

\draw
	(A\y X\a\b) -- (B\y X\a\b) -- (D\y X\a\b) -- (E\y X\a\b) -- (C\y X\a\b) -- cycle
	(D\y X\a\b) -- (F\y X\a\b)
	(C\y X\a\b) -- (GX\a\b) -- (HX\a\b) -- (I\y X\a\b)
	(E\y X\a\b) -- (I\y X\a\b) -- (J\y X\a\b) -- (K\y X\a\b)
	(L\y X\a\b) -- (MX\a\b)
	(B\y X\a\b) -- (NX\a\b) -- (OX\a\b) -- (P\y X\a\b);
	
\end{scope}
}

\foreach \a in {0,...,5}
\foreach \b in {0,...,10}
{\draw
	(A1X\a\b) -- (A-1X\a\b)
	(J1X\a\b) -- (J-1X\a\b);
	}


\draw[dotted, shift={(2*3 cm,-1.6*6 cm)}]
	(-0.8,-0.7) -- (-0.8,0.7);


\foreach \a / \b in {1/0,4/0,2/1,3/6,4/6,5/6,0/10,1/10,2/10,3/10}
\node at (Q\a\b) {\tiny $\times 1$};

\foreach \a in {0,...,5}
\foreach \b in {7,8,9}
\node at (Q\a\b) {\tiny $\times 1$};

\foreach \a / \b in {0/0,2/0,3/0,5/0,0/1,1/1,3/1,4/1,5/1,0/6,1/6,2/6}
\node at (Q\a\b) {\tiny $\times 2$};

\foreach \a in {0,...,5}
\foreach \b in {2,...,5}
\node at (Q\a\b) {\tiny $\times 2$};


\draw[line width=1.8]

	(C-1X00) -- (E-1X00)	
	(E1X00) -- (D1X00)
	(I1X00) -- (HX00)
	(A-1X00) -- (A1X00)
	(J-1X00) -- (K-1X00)
	(OX00) -- (NX00)

	(D-1X10) -- (E-1X10)	
	(E1X10) -- (D1X10)
	(GX10) -- (HX10)
	(A-1X10) -- (A1X10)
	(J-1X10) -- (J1X10)
	(OX10) -- (NX10)

	(C-1X20) -- (E-1X20)	
	(E1X20) -- (C1X20)
	(K1X20) -- (J1X20)
	(A-1X20) -- (A1X20)
	(J-1X20) -- (I-1X20)
	(OX20) -- (NX20)
	
	(D-1X30) -- (B-1X30)	
	(B1X30) -- (D1X30)
	(GX30) -- (HX30)
	(A-1X30) -- (A1X30)
	(J-1X30) -- (I-1X30)
	(J1X30) -- (K1X30)
	
	(D-1X40) -- (B-1X40)	
	(B1X40) -- (D1X40)
	(I1X40) -- (E1X40)
	(A-1X40) -- (A1X40)
	(J-1X40) -- (J1X40)
	(E-1X40) -- (I-1X40)
	
	(D-1X50) -- (B-1X50)	
	(B1X50) -- (D1X50)
	(I1X50) -- (E1X50)
	(A-1X50) -- (A1X50)
	(J-1X50) -- (K-1X50)
	(HX50) -- (I-1X50)
	
	(D-1X01) -- (F-1X01)	
	(B1X01) -- (D1X01)
	(I1X01) -- (E1X01)
	(A-1X01) -- (A1X01)
	(J-1X01) -- (J1X01)
	(E-1X01) -- (C-1X01)

	(D-1X11) -- (F-1X11)	
	(B1X11) -- (D1X11)
	(I1X11) -- (HX11)
	(A-1X11) -- (A1X11)
	(K1X11) -- (J1X11)
	(E-1X11) -- (C-1X11)

	(D-1X21) -- (F-1X21)	
	(F1X21) -- (D1X21)
	(C1X21) -- (E1X21)
	(A-1X21) -- (A1X21)
	(J-1X21) -- (J1X21)
	(E-1X21) -- (C-1X21)

	(J-1X31) -- (I-1X31)	
	(A1X31) -- (C1X31)
	(GX31) -- (HX31)
	(A-1X31) -- (B-1X31)
	(K1X31) -- (J1X31)
	(OX31) -- (NX31)

	(J1X41) -- (I1X41)	
	(A1X41) -- (C1X41)
	(GX41) -- (HX41)
	(A-1X41) -- (B-1X41)
	(K-1X41) -- (J-1X41)
	(OX41) -- (NX41)

	(E-1X51) -- (I-1X51)	
	(A1X51) -- (C1X51)
	(I1X51) -- (HX51)
	(A-1X51) -- (B-1X51)
	(K1X51) -- (J1X51)
	(OX51) -- (NX51)

	(E1X02) -- (I1X02)	
	(A1X02) -- (C1X02)
	(I-1X02) -- (HX02)
	(A-1X02) -- (B-1X02)
	(K-1X02) -- (J-1X02)
	(OX02) -- (NX02)

	(E1X12) -- (I1X12)	
	(A1X12) -- (C1X12)
	(I-1X12) -- (E-1X12)
	(A-1X12) -- (B-1X12)
	(J1X12) -- (J-1X12)
	(OX12) -- (NX12)

	(GX22) -- (HX22)	
	(B1X22) -- (NX22)
	(C1X22) -- (A1X22)
	(B-1X22) -- (D-1X22)
	(I-1X22) -- (J-1X22)
	(K1X22) -- (J1X22)

	(GX32) -- (HX32)	
	(B1X32) -- (NX32)
	(C1X32) -- (A1X32)
	(B-1X32) -- (D-1X32)
	(I1X32) -- (J1X32)
	(K-1X32) -- (J-1X32)

	(I1X42) -- (HX42)	
	(B1X42) -- (NX42)
	(C1X42) -- (A1X42)
	(B-1X42) -- (D-1X42)
	(I-1X42) -- (E-1X42)
	(K1X42) -- (J1X42)

	(I-1X52) -- (HX52)	
	(B1X52) -- (NX52)
	(C1X52) -- (A1X52)
	(B-1X52) -- (D-1X52)
	(I1X52) -- (E1X52)
	(K-1X52) -- (J-1X52)

	(I-1X03) -- (E-1X03)	
	(B1X03) -- (NX03)
	(C1X03) -- (A1X03)
	(B-1X03) -- (D-1X03)
	(I1X03) -- (E1X03)
	(J1X03) -- (J-1X03)

	(C-1X13) -- (E-1X13)	
	(B1X13) -- (NX13)
	(C1X13) -- (A1X13)
	(F-1X13) -- (D-1X13)
	(I1X13) -- (E1X13)
	(J1X13) -- (J-1X13)

	(C-1X23) -- (E-1X23)	
	(B1X23) -- (NX23)
	(C1X23) -- (A1X23)
	(F-1X23) -- (D-1X23)
	(I1X23) -- (HX23)
	(J1X23) -- (K1X23)

	(D1X33) -- (F1X33)	
	(GX33) -- (HX33)
	(C1X33) -- (A1X33)
	(J1X33) -- (J-1X33)
	(B-1X33) -- (A-1X33)
	(D-1X33) -- (F-1X33)

	(D1X43) -- (F1X43)	
	(I1X43) -- (HX43)
	(C1X43) -- (A1X43)
	(C-1X43) -- (E-1X43)
	(B-1X43) -- (NX43)
	(J-1X43) -- (K-1X43)

	(D1X53) -- (F1X53)	
	(GX53) -- (HX53)
	(C1X53) -- (A1X53)
	(D-1X53) -- (E-1X53)
	(B-1X53) -- (NX53)
	(J-1X53) -- (J1X53)

	(D1X04) -- (B1X04)	
	(I-1X04) -- (HX04)
	(C1X04) -- (GX04)
	(D-1X04) -- (E-1X04)
	(B-1X04) -- (NX04)
	(K1X04) -- (J1X04)

	(A1X14) -- (B1X14)	
	(I-1X14) -- (HX14)
	(C1X14) -- (GX14)
	(D-1X14) -- (E-1X14)
	(OX14) -- (NX14)
	(K1X14) -- (J1X14)

	(D1X24) -- (B1X24)	
	(I-1X24) -- (HX24)
	(C1X24) -- (GX24)
	(D-1X24) -- (F-1X24)
	(B-1X24) -- (A-1X24)
	(K1X24) -- (J1X24)

	(A-1X34) -- (B-1X34)	
	(I-1X34) -- (HX34)
	(C1X34) -- (GX34)
	(D1X34) -- (E1X34)
	(OX34) -- (NX34)
	(K-1X34) -- (J-1X34)

	(D-1X44) -- (B-1X44)	
	(I-1X44) -- (HX44)
	(C1X44) -- (GX44)
	(D1X44) -- (E1X44)
	(B1X44) -- (NX44)
	(K-1X44) -- (J-1X44)

	(D-1X54) -- (B-1X54)	
	(I-1X54) -- (HX54)
	(C1X54) -- (GX54)
	(D1X54) -- (F1X54)
	(B1X54) -- (A1X54)
	(K-1X54) -- (J-1X54)

	(A-1X05) -- (B-1X05)	
	(I-1X05) -- (E-1X05)
	(C1X05) -- (GX05)
	(D1X05) -- (E1X05)
	(OX05) -- (NX05)
	(J1X05) -- (J-1X05)

	(D-1X15) -- (B-1X15)	
	(I-1X15) -- (E-1X15)
	(C1X15) -- (GX15)
	(D1X15) -- (E1X15)
	(B1X15) -- (NX15)
	(J1X15) -- (J-1X15)

	(D-1X25) -- (B-1X25)	
	(I-1X25) -- (E-1X25)
	(C1X25) -- (GX25)
	(D1X25) -- (F1X25)
	(B1X25) -- (A1X25)
	(J1X25) -- (J-1X25)

	(D1X35) -- (B1X35)	
	(I-1X35) -- (J-1X35)
	(C1X35) -- (GX35)
	(C-1X35) -- (E-1X35)
	(B-1X35) -- (NX35)
	(K1X35) -- (J1X35)

	(D1X45) -- (B1X45)	
	(I1X45) -- (J1X45)
	(C1X45) -- (GX45)
	(C-1X45) -- (E-1X45)
	(B-1X45) -- (NX45)
	(K-1X45) -- (J-1X45)

	(F-1X55) -- (D-1X55)	
	(C-1X55) -- (E-1X55)
	(C1X55) -- (GX55)
	(D1X55) -- (E1X55)
	(B1X55) -- (NX55)
	(J1X55) -- (J-1X55)

	(A1X06) -- (B1X06)	
	(I-1X06) -- (J-1X06)
	(C1X06) -- (GX06)
	(C-1X06) -- (E-1X06)
	(OX06) -- (NX06)
	(K1X06) -- (J1X06)

	(A1X16) -- (B1X16)	
	(I1X16) -- (J1X16)
	(C1X16) -- (GX16)
	(C-1X16) -- (E-1X16)
	(OX16) -- (NX16)
	(K-1X16) -- (J-1X16)

	(F-1X26) -- (D-1X26)	
	(C-1X26) -- (E-1X26)
	(C1X26) -- (GX26)
	(D1X26) -- (F1X26)
	(B1X26) -- (A1X26)
	(J1X26) -- (J-1X26);


\foreach \a / \b in {3/6,4/6,5/6,0/10,1/10,2/10,3/10}
{
\draw[line width=1.8,dotted]
	(L1X\a\b) -- (MX\a\b)
	(L-1X\a\b) -- (MX\a\b);

\draw[line width=1.8]
	(OX\a\b) -- (P-1X\a\b);
}

\foreach \a in {0,...,5}
\foreach \b in {7,8,9}
{
\draw[line width=1.8,dotted]
	(L1X\a\b) -- (MX\a\b)
	(L-1X\a\b) -- (MX\a\b);

\draw[line width=1.8]
	(OX\a\b) -- (P-1X\a\b);
}

\draw[line width=1.8]

	(D1X36) -- (F1X36)	
	(C1X36) -- (E1X36)	
	(A-1X36) -- (A1X36)	
	(C-1X36) -- (E-1X36)	
	(I-1X36) -- (J-1X36)

	(D1X46) -- (B1X46)	
	(I1X46) -- (E1X46)
	(A-1X46) -- (A1X46)	
	(C-1X46) -- (E-1X46)	
	(I-1X46) -- (J-1X46)
	
	(D1X56) -- (B1X56)	
	(I1X56) -- (J1X56)
	(A-1X56) -- (A1X56)	
	(D-1X56) -- (E-1X56)	
	(GX56) -- (HX56)

	(D1X07) -- (F1X07)	
	(C1X07) -- (E1X07)
	(A-1X07) -- (A1X07)	
	(D-1X07) -- (E-1X07)	
	(I-1X07) -- (HX07)

	(D1X17) -- (B1X17)	
	(I1X17) -- (E1X17)
	(A-1X17) -- (A1X17)	
	(D-1X17) -- (E-1X17)
	(I-1X17) -- (HX17)

	(D1X27) -- (F1X27)	
	(I1X27) -- (HX27)
	(C1X27) -- (A1X27)
	(I-1X27) -- (E-1X27)
	(B-1X27) -- (A-1X27)

	(D1X37) -- (F1X37)	
	(GX37) -- (HX37)
	(C1X37) -- (A1X37)
	(I-1X37) -- (J-1X37)
	(B-1X37) -- (A-1X37)

	(C-1X47) -- (E-1X47)	
	(B1X47) -- (NX47)
	(C1X47) -- (A1X47)
	(I-1X47) -- (J-1X47)
	(E1X47) -- (I1X47)

	(D-1X57) -- (E-1X57)	
	(B1X57) -- (NX57)
	(C1X57) -- (A1X57)
	(I-1X57) -- (HX57)
	(E1X57) -- (I1X57)

	(D-1X08) -- (E-1X08)	
	(B1X08) -- (NX08)
	(C1X08) -- (A1X08)
	(GX08) -- (HX08)
	(J1X08) -- (I1X08)

	(I-1X18) -- (E-1X18)	
	(B1X18) -- (D1X18)
	(C1X18) -- (GX18)
	(B-1X18) -- (A-1X18)
	(J1X18) -- (I1X18)

	(C-1X28) -- (E-1X28)	
	(E1X28) -- (D1X28)
	(C1X28) -- (GX28)
	(B1X28) -- (NX28)
	(J-1X28) -- (I-1X28)
	
	(D-1X38) -- (E-1X38)	
	(E1X38) -- (D1X38)
	(C1X38) -- (GX38)
	(B1X38) -- (NX38)
	(HX38) -- (I-1X38)

	(C-1X48) -- (E-1X48)	
	(F1X48) -- (D1X48)
	(C1X48) -- (GX48)
	(B1X48) -- (A1X48)
	(J-1X48) -- (I-1X48)

	(D-1X58) -- (E-1X58)	
	(F1X58) -- (D1X58)
	(C1X58) -- (GX58)
	(B1X58) -- (A1X58)
	(HX58) -- (I-1X58)
	
	(I1X09) -- (E1X09)	
	(B-1X09) -- (A-1X09)
	(C-1X09) -- (GX09)
	(B1X09) -- (D1X09)
	(J-1X09) -- (I-1X09)
	
	(C1X19) -- (E1X19)	
	(B-1X19) -- (A-1X19)
	(C-1X19) -- (GX19)
	(F1X19) -- (D1X19)
	(J-1X19) -- (I-1X19)

	(I1X29) -- (HX29)	
	(B1X29) -- (NX29)
	(C-1X29) -- (GX29)
	(E1X29) -- (D1X29)
	(D-1X29) -- (E-1X29)

	(I1X39) -- (J1X39)	
	(B1X39) -- (NX39)
	(C-1X39) -- (GX39)
	(E1X39) -- (C1X39)
	(D-1X39) -- (E-1X39)

	(I1X49) -- (HX49)	
	(B1X49) -- (A1X49)
	(C-1X49) -- (GX49)
	(F1X49) -- (D1X49)
	(D-1X49) -- (E-1X49)
	
	(I1X59) -- (HX59)	
	(B1X59) -- (NX59)
	(C-1X59) -- (A-1X59)
	(E1X59) -- (D1X59)
	(I-1X59) -- (E-1X59)
	
	(GX010) -- (HX010)	
	(B1X010) -- (NX010)
	(C-1X010) -- (A-1X010)
	(E1X010) -- (D1X010)
	(I-1X010) -- (J-1X010)
	
	(E-1X110) -- (I-1X110)	
	(B1X110) -- (NX110)
	(C-1X110) -- (A-1X110)
	(E1X110) -- (C1X110)
	(I1X110) -- (J1X110)
	
	(I1X210) -- (HX210)	
	(B1X210) -- (A1X210)
	(C-1X210) -- (A-1X210)
	(F1X210) -- (D1X210)
	(I-1X210) -- (E-1X210)
	
	(GX310) -- (HX310)	
	(B1X310) -- (A1X310)
	(C-1X310) -- (A-1X310)
	(F1X310) -- (D1X310)
	(I-1X310) -- (J-1X310);

\fill[white,shift={(8 cm, -16 cm)}] (-0.8,-0.6) rectangle (3,0.6);

\end{tikzpicture}
\caption{Timezone tiling, distance $2$, edge combination $a^4b$.}
\label{timezone4ab_d2}
\end{figure}

Up to the symmetry of vertical flipping, the complete list for the even family (i.e., both boundaries are $aa$) of timezone tilings is given by the first $39$ tilings in Figure \ref{timezone4ab_d2}. By applying the vertical flipping, those labeled $\times 2$ give new timezone tilings. So there are total of $75$ timezone tilings with $aa$ as both boundary meridians.

Similarly, up to symmetry, the complete list for the odd family (i.e., the boundaries are $ab$ or $ba$) of timezone tilings is given by the last $25$ tilings in Figure \ref{timezone4ab_d2}. In case the right boundary is $ab$, we get total of $25$ timezone tilings with $ab$ as both boundaries. The vertical flippings of these give total of $25$ timezone tilings with $ba$ as both boundaries. On the other hand, in case the right boundary is $ba$, we get total of $25$ timezone tilings with $ab$ and $ba$ as left and right boundaries. The vertical flippings of these give total of $25$ timezone tilings with $ba$ and $ab$ as left and right boundaries.

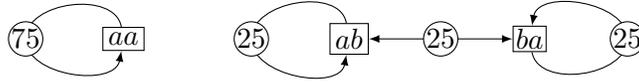
\begin{figure}[htp]
\centering
\begin{tikzpicture}[>=latex]


\node[rectangle,draw,inner sep=2] (A) at (-1.2,0) {\small $aa$};
\draw[->] (A) to[out=100,in=90] (-2.5,0) to[out=270,in=260] (A);
\node[draw,shape=circle,fill=white,inner sep=0] at (-2.5,0) {\small $75$};	


\begin{scope}[xshift=3cm]

\node[rectangle,draw,inner sep=2] (A) at (-1.2,0) {\small $ab$};
\node[rectangle,draw,inner sep=1] (C) at (1.2,0) {\small $ba$};	

\draw[->] (A) to[out=100,in=90] (-2.5,0) to[out=270,in=260] (A);
\node[draw,shape=circle,fill=white,inner sep=0] at (-2.5,0) {\small $25$};	

\draw[->] (C) to[out=280,in=270] (2.5,0) to[out=90,in=80] (C);
\node[draw,shape=circle,fill=white,inner sep=0] at (2.5,0) {\small $25$};	
	
\draw (A) edge[<->] node[draw,shape=circle,fill=white,inner sep=0] {\small $25$} (C);

\end{scope}
	
\end{tikzpicture}
\caption{Earth map tiling, distance $2$, edge combination $a^4b$.}
\label{4ab_d2}
\end{figure}

All the timezone tilings for earth map tiling of distance $2$ are represented by the arrows in Figure \ref{4ab_d2}. An edge congruent earth map tiling is a directed cycle in the graph. In particular, the even family is obtained by choosing any sequence of $75$ timezone tilings and glueing them together.

\subsection*{Distance 1}

\begin{figure}[htp]
\centering
\begin{tikzpicture}[>=latex]

\foreach \x in {-1,1}
\foreach \y in {-1,1}
\foreach \a in {0,...,5}
\foreach \b in {0,...,6}
{
\begin{scope}[shift={(2*\a cm, -1.6*\b cm)}]

\coordinate  (A\y X\a\b) at (0,0.1*\y);
\coordinate  (B\x\y X\a\b) at (0.2*\x,0.2*\y);
\coordinate  (C\x\y X\a\b) at (0.2*\x,0.35*\y);
\coordinate  (D\x X\a\b) at (0.3*\x,0);
\coordinate  (E\x X\a\b) at (0.45*\x,0);
\coordinate  (F\x\y X\a\b) at (0.5*\x,0.4*\y);
\coordinate  (G\x\y X\a\b) at (0.55*\x,0.6*\y);
\coordinate  (H\x\y X\a\b) at (0.7*\x,0.6*\y);

\coordinate  (Q\a\b) at (0,-0.7);

\draw
	(A\y X\a\b) -- (B\x\y X\a\b) -- (D\x X\a\b) -- (E\x X\a\b)  -- (F\x\y X\a\b) 
	(B\x\y X\a\b) -- (C\x\y X\a\b) -- (F\x\y X\a\b)  -- (G\x\y X\a\b);
	
\end{scope}
}

\foreach \x in {-1,1}
\foreach \a in {0,...,5}
\foreach \b in {0,...,6}
\draw
	(A1X\a\b) -- (A-1X\a\b)
	(C1\x X\a\b) -- (C-1\x X\a\b)
	(H\x 1X\a\b) -- (H\x -1X\a\b);


\draw[dotted, shift={(2*5 cm,-1.6*4 cm)}]
	(-1,-0.7) -- (-1,0.7);


\foreach \a / \b in {0/0,1/0,2/0,3/0,4/0,1/1,2/1,0/2,3/2,4/2,5/2,0/3,2/3,3/3,4/3,5/3,0/4,1/4,2/4,3/4,4/4,2/5,3/5,4/5,5/5,0/6}
\node at (Q\a\b) {\tiny $\times 4$};

\foreach \a / \b in {5/0,0/1,3/1,4/1,5/1,1/2,2/2,1/3,0/5,1/5}
\node at (Q\a\b) {\tiny $\times 2$};

\node at (Q54) {\tiny $\times 1$};


\draw[line width=1.8]
	(A1X00) -- (A-1X00)		
	(D-1X00) -- (E-1X00)
	(C11X00) -- (C-11X00)	
	(F11X00) -- (E1X00)	
	(B1-1X00) -- (C1-1X00)
	(F-1-1X00) -- (G-1-1X00)
	
	(A1X10) -- (A-1X10)		
	(D-1X10) -- (E-1X10)	
	(B11X10) -- (C11X10)
	(B1-1X10) -- (C1-1X10)
	(F-11X10) -- (G-11X10)
	(F1-1X10) -- (G1-1X10)
	
	(A1X20) -- (A-1X20)	
	(C-11X20) -- (C11X20)
	(A1X20) -- (A-1X20)	
	(E-1X20) -- (F-11X20)
	(E1X20) -- (F11X20)
	(B-1-1X20) -- (C-1-1X20)
	(C1-1X20) -- (F1-1X20)	
	
	(A1X30) -- (A-1X30)	
	(F-11X30) -- (G-11X30)
	(B-11X30) -- (C-11X30)	
	(E1X30) -- (F11X30)	
	(B-1-1X30) -- (C-1-1X30)	
	(C1-1X30) -- (F1-1X30)	
	
	(A1X40) -- (A-1X40)	
	(F-11X40) -- (G-11X40)
	(B-11X40) -- (C-11X40)	
	(E1X40) -- (F11X40)	
	(B1-1X40) -- (C1-1X40)	
	(C-1-1X40) -- (F-1-1X40)
	
	(B-11X50) -- (A1X50)	
	(A-1X50) -- (B1-1X50)
	(E-1X50) -- (D-1X50)	
	(E1X50) -- (D1X50)	
	(G-11X50) -- (F-11X50)
	(G1-1X50) -- (F1-1X50)
	
	(B-11X01) -- (A1X01)	
	(A-1X01) -- (B1-1X01)
	(E-1X01) -- (D-1X01)	
	(E1X01) -- (D1X01)	
	(G-1-1X01) -- (F-1-1X01)
	(G11X01) -- (F11X01)
	
	(B-11X11) -- (A1X11)	
	(A-1X11) -- (B1-1X11)
	(E-1X11) -- (D-1X11)	
	(E1X11) -- (F11X11)	
	(C1-1X11) -- (F1-1X11)
	(G-11X11) -- (F-11X11)
	
	(B-11X21) -- (A1X21)	
	(A-1X21) -- (B1-1X21)
	(E-1X21) -- (D-1X21)	
	(C11X21) -- (F11X21)	
	(E1X21) -- (F1-1X21)
	(G-1-1X21) -- (F-1-1X21)
	
	(B-11X31) -- (A1X31)	
	(A-1X31) -- (B1-1X31)
	(C-11X31) -- (F-11X31)	
	(C1-1X31) -- (F1-1X31)	
	(E-1X31) -- (F-1-1X31)
	(E1X31) -- (F11X31)	
	
	(B-11X41) -- (A1X41)	
	(A-1X41) -- (B1-1X41)
	(C11X41) -- (F11X41)	
	(C-1-1X41) -- (F-1-1X41)	
	(E1X41) -- (F1-1X41)
	(E-1X41) -- (F-11X41)
	
	(B-11X51) -- (D-1X51)	
	(D1X51) -- (B1-1X51)
	(C-11X51) -- (C11X51)	
	(C-1-1X51) -- (C1-1X51)	
	(E-1X51) -- (F-1-1X51)
	(E1X51) -- (F11X51)	
	
	(B-11X02) -- (D-1X02)	
	(D1X02) -- (B1-1X02)
	(C-11X02) -- (C11X02)	
	(E1X02) -- (F11X02)		
	(B-1-1X02) -- (C-1-1X02)	
	(G-1-1X02) -- (F-1-1X02)	
	
	(B-11X12) -- (D-1X12)	
	(D1X12) -- (B1-1X12)
	(C11X12) -- (B11X12)	
	(C-1-1X12) -- (B-1-1X12)		
	(G-11X12) -- (F-11X12)	
	(G1-1X12) -- (F1-1X12)	
	
	(B-11X22) -- (D-1X22)	
	(D1X22) -- (B1-1X22)
	(C11X22) -- (B11X22)	
	(C-1-1X22) -- (B-1-1X22)		
	(G11X22) -- (F11X22)	
	(G-1-1X22) -- (F-1-1X22)	
	
	(B-11X32) -- (D-1X32)	
	(A-1X32) -- (B1-1X32)
	(C11X32) -- (C-11X32)	
	(F-1-1X32) -- (E-1X32)		
	(D1X32) -- (E1X32)	
	(G1-1X32) -- (F1-1X32)	
	
	(B-11X42) -- (D-1X42)	
	(A-1X42) -- (B1-1X42)
	(C11X42) -- (C-11X42)	
	(F-1-1X42) -- (E-1X42)		
	(F11X42) -- (E1X42)	
	(C1-1X42) -- (F1-1X42)	
	
	(B-11X52) -- (D-1X52)	
	(A-1X52) -- (B1-1X52)
	(C11X52) -- (B11X52)	
	(G-11X52) -- (F-11X52)	
	(F-1-1X52) -- (C-1-1X52)		
	(F1-1X52) -- (E1X52)	
	
	(B-11X03) -- (D-1X03)	
	(A-1X03) -- (B1-1X03)
	(C11X03) -- (B11X03)	
	(F-1-1X03) -- (E-1X03)		
	(F11X03) -- (G11X03)	
	(C1-1X03) -- (F1-1X03)	
	
	(B-11X13) -- (D-1X13)	
	(B11X13) -- (D1X13)
	(C11X13) -- (C-11X13)	
	(C-1-1X13) -- (C1-1X13)
	(F1-1X13) -- (E1X13)		
	(F-1-1X13) -- (E-1X13)	
	
	(B-11X23) -- (D-1X23)	
	(B11X23) -- (D1X23)
	(C11X23) -- (C-11X23)	
	(F1-1X23) -- (E1X23)		
	(F-1-1X23) -- (G-1-1X23)	
	(B-1-1X23) -- (C-1-1X23)
	
	(B-11X33) -- (D-1X33)	
	(B11X33) -- (A1X33)
	(C1-1X33) -- (C-1-1X33)	
	(F-1-1X33) -- (E-1X33)		
	(C11X33) -- (F11X33)	
	(E1X33) -- (F1-1X33)	
	
	(B-11X43) -- (D-1X43)	
	(B11X43) -- (A1X43)
	(C1-1X43) -- (C-1-1X43)	
	(F-1-1X43) -- (E-1X43)		
	(D1X43) -- (E1X43)	
	(G11X43) -- (F11X43)	
	
	(B-11X53) -- (D-1X53)	
	(B11X53) -- (A1X53)
	(D1X53) -- (E1X53)
	(B-1-1X53) -- (C-1-1X53)	
	(F-11X53) -- (G-11X53)	
	(G1-1X53) -- (F1-1X53)	
	
	(B-11X04) -- (D-1X04)	
	(B11X04) -- (A1X04)
	(D1X04) -- (E1X04)
	(B-1-1X04) -- (C-1-1X04)	
	(F11X04) -- (G11X04)	
	(G-1-1X04) -- (F-1-1X04)
	
	(B-11X14) -- (D-1X14)	
	(B11X14) -- (A1X14)
	(B-1-1X14) -- (C-1-1X14)
	(F11X14) -- (E1X14)	
	(C1-1X14) -- (F1-1X14)	
	(G-11X14) -- (F-11X14)	
		
	(B-11X24) -- (D-1X24)	
	(B11X24) -- (A1X24)
	(B-1-1X24) -- (C-1-1X24)
	(F11X24) -- (C11X24)	
	(E1X24) -- (F1-1X24)	
	(G-1-1X24) -- (F-1-1X24)
	
	(B-11X34) -- (D-1X34)	
	(B11X34) -- (A1X34)
	(B1-1X34) -- (C1-1X34)
	(F11X34) -- (E1X34)	
	(C-1-1X34) -- (F-1-1X34)	
	(F-11X34) -- (G-11X34)	
	
	(B-11X44) -- (D-1X44)	
	(B11X44) -- (A1X44)
	(B1-1X44) -- (C1-1X44)
	(F11X44) -- (C11X44)	
	(G1-1X44) -- (F1-1X44)	
	(F-1-1X44) -- (E-1X44)					
	
	;
			

\draw[line width=1.8]
	(H11X54) -- (H1-1X54)	
	(H-11X54) -- (H-1-1X54)
	(A1X54) -- (A-1X54)
	(D-1X54) -- (E-1X54)	
	(D1X54) -- (E1X54)
	(C11X54) -- (C-11X54)		
	(C-1-1X54) -- (C1-1X54)
	
	(H11X05) -- (H1-1X05)	
	(H-11X05) -- (H-1-1X05)
	(A1X05) -- (A-1X05)	
	(C-11X05) -- (F-11X05)
	(C11X05) -- (B11X05)	
	(B-1-1X05) -- (C-1-1X05)
	(C1-1X05) -- (F1-1X05)
	
	(H11X15) -- (H1-1X15)	
	(H-11X15) -- (H-1-1X15)
	(A1X15) -- (A-1X15)	
	(C-11X15) -- (F-11X15)
	(C11X15) -- (B11X15)
	(B1-1X15) -- (C1-1X15)
	(C-1-1X15) -- (F-1-1X15)	
	
	(H11X25) -- (H1-1X25)	
	(H-11X25) -- (H-1-1X25)
	(B-11X25) -- (A1X25)	
	(B1-1X25) -- (A-1X25)	
	(D-1X25) -- (E-1X25)	
	(C11X25) -- (F11X25)	
	(C1-1X25) -- (F1-1X25)	
	
	(H11X35) -- (H1-1X35)	
	(H-11X35) -- (H-1-1X35)
	(B-11X35) -- (D-1X35)
	(B1-1X35) -- (A-1X35)
	(D1X35)	-- (E1X35)
	(C-11X35) -- (C11X35)
	(C-1-1X35) -- (F-1-1X35)
	
	(H11X45) -- (H1-1X45)	
	(H-11X45) -- (H-1-1X45)
	(B-11X45) -- (D-1X45)	
	(B11X45) -- (D1X45)	
	(C-11X45) -- (C11X45)
	(B-1-1X45) -- (C-1-1X45)
	(C1-1X45) -- (F1-1X45)

	(H11X55) -- (H1-1X55)	
	(H-11X55) -- (H-1-1X55)
	(B-11X55) -- (D-1X55)	
	(B11X55) -- (A1X55)	
	(C11X55) -- (F11X55)
	(B-1-1X55) -- (C-1-1X55)
	(C1-1X55) -- (F1-1X55)

	(H11X06) -- (H1-1X06)	
	(H-11X06) -- (H-1-1X06)
	(B-11X06) -- (D-1X06)	
	(B11X06) -- (A1X06)	
	(C11X06) -- (F11X06)
	(B1-1X06) -- (C1-1X06)
	(C-1-1X06) -- (F-1-1X06)
			
	;

\fill[white,shift={(2 cm, -9.6 cm)}] (-0.8,-0.6) rectangle (8.8,0.6);

\end{tikzpicture}
\caption{Timezone tiling, distance $1$, edge combination $a^4b$.}
\label{timezone4ab_d1}
\end{figure}
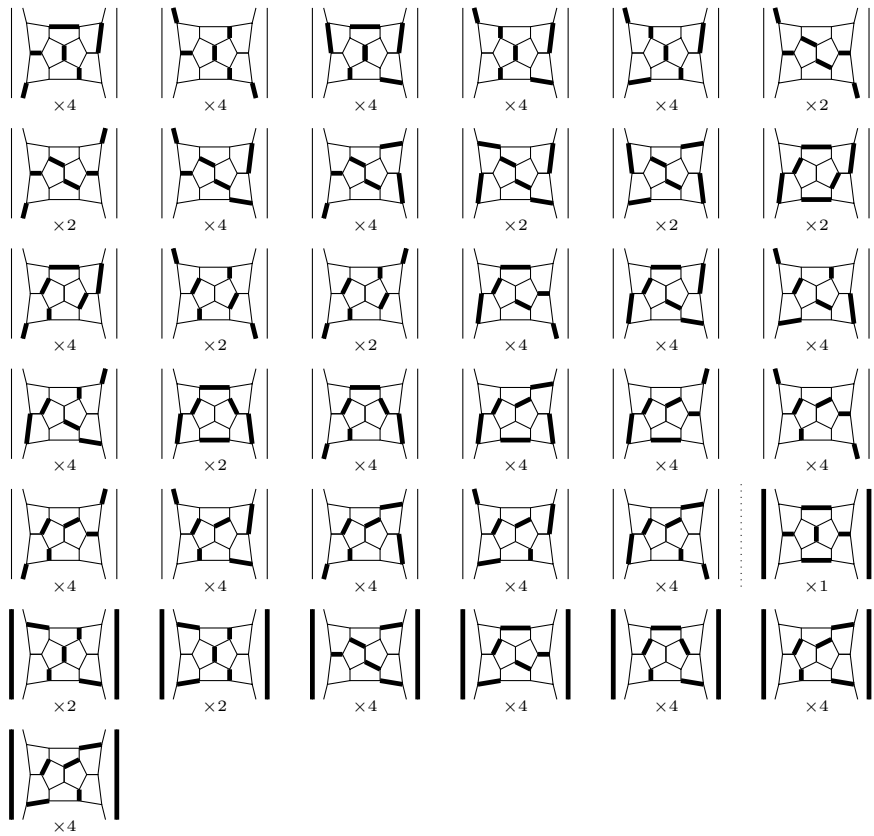

For earth map tilings of distance $1$, the even family simply means the meridians are $a$, and the odd family means the meridians are $b$. Up to symmetry, Figure \ref{timezone4ab_d1} gives all the timezone tilings. We have total of $100$ timezone tilings with $a$ as both boundary meridians. We also have total of $25$ timezone tilings with $b$ as both boundary meridians. By glueing any sequence of timezone tilings in any one of the two families, we get all the edge congruent earth map tilings of distance $1$.

\end{document}